\documentclass[a4paper,11pt]{article}

\usepackage[pages=all, color=black, position={current page.south}, placement=bottom, scale=1, opacity=1, vshift=5mm]{background}

 \SetBgContents{
 }      % copyright

\usepackage[margin=1in]{geometry} % full-width

\usepackage{amsmath}
\usepackage{amsthm}
\usepackage{subcaption} % Per i subfigure
\usepackage{graphicx}
\usepackage{amsmath}
\usepackage{wrapfig}
\numberwithin{equation}{section}

\usepackage{amsfonts}
\usepackage{amssymb}
\usepackage{graphicx}
\usepackage{enumitem}   % per personalizzare enumerate
\usepackage{booktabs}
\usepackage{longtable} % Per tabelle che si estendono su più pagine, se necessario
\usepackage{ragged2e} % Per controllo dell'allineamento del testo nelle celle
\usepackage[absolute,overlay]{textpos} % Per posizionamento assoluto
\usepackage{overpic} % Per sovrapporre immagini
\usepackage{graphicx}
\usepackage{subcaption} % Per i subfigure
\usepackage{amsopn}
\usepackage{multirow}
\usepackage[normalem]{ulem}

\usepackage[linesnumbered, ruled, vlined]{algorithm2e}
\SetAlgoLined % Disegna linee verticali per i blocchi While/If
\DontPrintSemicolon % Non stampa i punti e virgola alla fine di ogni riga
\SetKwFunction{KwDefine}{Define} % Per definire la funzione T_tau(x)

\usepackage[utf8]{inputenc}
\usepackage{hyperref}
\hypersetup{
	unicode,
	pdfauthor={Author One, Author Two, Author Three},
	pdftitle={Deep Equilibrium models for Poisson Imaging Inverse problems via Mirror Descent},
	pdfsubject={A simple article template},
	pdfkeywords={article, template, simple},
	pdfproducer={LaTeX},
	pdfcreator={pdflatex}
}

\usepackage{mathtools}
\mathtoolsset{showonlyrefs}
\usepackage[sort&compress,numbers,square]{natbib}

\theoremstyle{plain}
\newtheorem{theorem}{Theorem}[section]
\newtheorem{corollary}[theorem]{Corollary}
\newtheorem{lemma}[theorem]{Lemma}

\newtheorem{proposition}[theorem]{Proposition}

\newtheorem{definition}[theorem]{Definition}
\newtheorem{example}[theorem]{Example}
\newtheorem{remark}[theorem]{Remark}

\newtheorem{notation}{Notation}[section] % numerazione opzionale per sezione

\usepackage{graphicx, color}
\graphicspath{{fig/}}

\usepackage{mathrsfs} % for \mathscr command

\usepackage{lipsum}
\usepackage{etoolbox} % utile per definire ambienti personalizzati

\newenvironment{keywords}{
    \noindent\textbf{Keywords:} 
}{\par}

\newenvironment{MSCcodes}{
    \noindent\textbf{MSC codes:} 
}{\par}

\title{Deep Equilibrium models for Poisson imaging inverse problems via Mirror Descent}
\author{Christian Daniele$^1$ \and Silvia Villa$^2$ \and Samuel Vaiter$^3$ \and Luca Calatroni$^4$}

\date{
	$^1$ MaLGa Center, DIBRIS, University of Genoa, Italy   \\ \texttt{christian.daniele@edu.unige.it}\\%
	$^2$MaLGa Center, DIIMA, University of Genoa, Italy \\ \texttt{silvia.villa@unige.it}\\
    $^3$ CNRS \& Universit\'e C\^ote d’Azur, Laboratoire J. A. Dieudonn\'e. Nice, France. \\ \texttt{samuel.vaiter@cnrs.fr}\\
    MaLGa Center, DIBRIS, University of Genoa \& MMS, Istituto Italiano di Tecnologia (IIT), Italy\\
    \texttt{luca.calatroni@unige.it} \\[2ex]%
}

\begin{document}
	\maketitle
	
\begin{abstract}
Deep Equilibrium Models (DEQs) are implicit neural networks with fixed points, which have recently gained attention for learning image regularization functionals, particularly in settings involving Gaussian fidelities, where assumptions on the forward operator ensure contractiveness of standard (proximal) Gradient Descent operators. In this work, we extend the application of DEQs to Poisson inverse problems, where the data fidelity term is more appropriately modeled by the Kullback–Leibler divergence. To this end, we introduce a novel DEQ formulation based on Mirror Descent defined in terms of a tailored non-Euclidean geometry that naturally adapts with the structure of the data term. This enables the learning of neural regularizers within a principled training framework. We derive sufficient conditions and establish refined convergence results based on the Kurdyka–Łojasiewicz framework for subanalytic functions with non-closed domains to guarantee the convergence of the learned reconstruction scheme and propose computational strategies that enable both efficient training and parameter-free inference. Numerical experiments show that our method outperforms traditional model-based approaches and it is comparable to the performance of Bregman Plug-and-Play methods, while mitigating their typical drawbacks, such as time-consuming tuning of hyper-parameters.\\ The code is publicly available at \url{https://github.com/christiandaniele/DEQ-MD}.
\end{abstract}

% REQUIRED

% REQUIRED

\section{Introduction}

Image reconstruction problems in the presence of Poisson noise arise naturally in various scientific imaging domains where the measurements are governed by photon-limited observations. Prominent examples include fluorescence microscopy \cite{Bertero_2009}, astronomical imaging \cite{Dupe2009}, and Positron Emission Tomography \cite{Shepp1982}, all of which involve scenarios where photon counts are inherently low and follow discrete statistics.

Mathematically, given dimensions $m, n \in \mathbb{N}_*$, the task is to recover an unknown non-negative image $x \in \mathbb{R}^n_{\geq 0}$ from Poisson-degraded measurements $y \in \mathbb{R}^m_{\geq 0 }$, modeled as:
\begin{equation}
y \sim \text{Poiss}(Ax), \label{PIP}
\end{equation}
where $A \in \mathbb{R}^{m \times n}_{\geq 0}$ 
is a known, typically ill-posed, forward operator. The operator $\text{Poiss}(\cdot)$ denotes a multidimensional Poisson distribution, where each component $y_j$ satisfies $\mathbb{E}(y_j) = \text{Var}(y_j) = (Ax)_j$ for $j = 1, \ldots, m$. We use the convention\footnote{This choice is motivated by the property that for a sequence of random variables $\{X_n\}_{n \in \mathbb{N}}$ with \\$X_n \sim \text{Poiss}(\lambda_n)$, if $\lambda_n \to 0$ as $n \to \infty$, then $X_n \xrightarrow{d} \mathbf{0}$, where $\xrightarrow{d}$ denotes convergence in distribution.} Poiss(0)=$\mathbf{0}$, where $\mathbf{0}$ is the distribution s.t. $\mathbb{E}(\mathbf{0})=\text{Var}({\mathbf{0}})=0$. This indicates the case where the Poisson distribution's parameter is 0, resulting in a degenerate Poisson distribution. 
Note that the formulation in \eqref{PIP} does not include an explicit background term $b \in \mathbb{R}^m_{>0}$, which is sometimes introduced in practice to avoid degeneracy in regions of low expected intensity \cite{KLsmoothed}.

In comparison with standard Gaussian noise modeling, the Poisson assumption can be easily shown to correspond to a non-quadratic log-likelihood data term which is the Kullback-Leibler (KL) divergence function $\text{KL}: \mathbb{R}^m_{\geq 0} \times \mathbb{R}^m_{>0} \to \mathbb{R}_{\geq 0}$ that reads:
\begin{equation}
	\text{KL}(y,Ax) := \sum_{i=1}^{m} y_i \log\left( \frac{y_i}{(Ax)_i} \right) + (Ax)_i - y_i, \label{KL equation}
\end{equation}
which is well-defined also for vanishing values $y_j$ upon the convention $0\log 0 = 0$. In this framework, to ensure that $(Ax)_i > 0$ for all $i \in \{1, \dots, m\}$, we consider a forward operator $A$ such that its $i$-th row vector, $a_i$, satisfies $a_i \in \mathbb{R}_{\geq 0}^n \setminus \{0\}$ for all $i \in \{1, \dots, m\}$, and $x \in \mathbb{R}_{>0}^n$. These conditions are not overly restrictive, for instance, a standard convolution operator typically satisfies them.

\medskip

Early contributions to image reconstruction based on Kullback–Leibler (KL) divergence minimization date back to the foundational works of \cite{Richardson72, Lucy1974, Shepp1982}, which introduced iterative regularization schemes grounded in the Expectation-Maximization (EM) framework. A notable instance is the Richardson–Lucy algorithm, developed independently in \cite{Richardson72, Lucy1974}, which applies EM updates to linear image deconvolution under Poisson noise. 
Building on this framework, Byrne proposed in \cite{Byrne1993} the Multiplicative Algebraic Reconstruction Technique (MART), a class of iterative methods derived from cross-entropy (Kullback–Leibler) minimization using multiplicative updates under non-negativity constraints. While MART (and its simultaneous variant SMART) aim at minimizing $\mathrm{KL}(A x, y)$, this objective differs from the more common Poisson-likelihood divergence $\mathrm{KL}(y, A x)$. In this respect, the expectation-maximisation maximum-likelihood (EMML) algorithm — and in particular Byrne’s accelerated block-iterative variant (RBI-EMML) as described in \cite{ByrneAccEMML1998} — align more closely with the standard Poisson noise model.

While these methods remain effective and are still used in modern applications (see, e.g., \cite{Zunino2023_ISMinverse}), they typically require early stopping to prevent overfitting to noise, which introduces a reliance on heuristic tuning to obtain satisfactory results.

To address these limitations, a wide range of model-based regularization techniques have been developed to incorporate explicit prior knowledge about the image, often in combination with smoothed variants of the KL divergence in \eqref{KL equation} to enable standard gradient-based optimization. A common strategy involves introducing a small background term \cite{KLsmoothed}, resulting in the modified fidelity $\mathrm{KL}(y, Ax + b)$, which is twice continuously differentiable on the non-negative orthant. To ensure non-negativity of the reconstruction, projected Gradient Descent and its scaled or preconditioned variants have been explored \cite{Bonettini_2009}. In parallel, more structured priors such as sparsity-promoting $\ell_1$-norm penalties \cite{Dupe2009, Figuereido2010}, total variation (TV) \cite{Sawatzky2009}, and higher-order variants \cite{DiSerafino2021} have been successfully employed. These models are typically optimized using (accelerated) proximal gradient algorithms \cite{Bonettini2018a, RebegoldiCalatroni2022} or primal–dual methods \cite{Chambolle2018}, often incorporating advanced features such as variable metric scaling, inexact proximal steps, and stochastic updates to improve convergence and computational efficiency.

\medskip

Directly minimizing the KL divergence in \eqref{KL equation} without adding a background smoothing term poses significant algorithmic challenges. Although the functional $\mathrm{KL}(y, A \cdot)$ is differentiable on $\mathbb{R}^n_{>0}$, its gradient is not Lipschitz continuous therein due to the singularity near zero. This prevents the straightforward application of standard (proximal) gradient algorithms. To overcome this, Bregman-type optimization methods—such as those developed in \cite{nolip, Bolte2018,Teboulle2018}—relax the need for Lipschitz-smoothness and naturally accommodate domain constraints, for instance, non-negativity. A prominent example is Mirror Descent (MD) \cite{nemirovskii_yudin_1983}, whose convergence properties can be studied under a non-Euclidean analogue of $L$-smoothness, defined relative to a strictly convex function $h$ that specifies the underlying geometry. 
%In the Euclidean setting, this corresponds to the familiar case $h(x) = \frac{1}{2}\|x\|^2$.
In the context of Poisson inverse problems, this notion of relative smoothness amounts to the existence of a constant $L > 0$ such that the function $Lh - \mathrm{KL}(y, A\cdot)$ is convex \cite{nolip, Bolte2018, RelativeSmoothness}. For the KL divergence, this condition is satisfied by choosing $h$ as Burg's entropy, defined by $h(x) = -\sum_{i=1}^n \log(x_i)$. Based on this choice, one can construct Bregman (proximal) gradient operators using the associated Bregman distance, leading to optimization schemes that share strong analogies with the multiplicative updates found in EM and MART algorithms.

\medskip

Leveraging the connection between optimization algorithms and deep learning, we propose a novel approach to Poisson inverse problems based on the MD iteration interpreted as the recurrent core of a deep equilibrium model (DEQ). 
%In this architecture, each MD iteration corresponds to a network layer, w. 
Our goal is to learn image regularization functionals in a principled and convergent manner—without requiring proximal operators or contractiveness of classical descent schemes. This builds upon the broader idea of using optimization-inspired architectures in image reconstruction. Among them, algorithm unrolling \cite{Unrolling,AdlerOktem2018} treats each iteration of an algorithm as a trainable neural layer, while Plug-and-Play (PnP) methods \cite{Venkatakrishnan2013,meinhardt2017learning,Kamilov2023} replace proximal steps with learned Gaussian denoisers, often interpreted through Tweedie's formula \cite{laumont2022bayesian}. Variants such as RED \cite{Romano2017}, contractive network models \cite{sun2018online,sun2021scalable,Terris}, and Bregman-PnP approaches \cite{Shabili2022,Hurault} offer convergence guarantees under various assumptions. In particular, \cite{Hurault} analyzes convergence in the Poisson setting using relative smoothness and suitable training upon inverse Gamma noise statistics, while \cite{klatzer2025} introduces a mirror-Langevin sampler tailored to Poisson likelihoods, handling non-negativity and lack of $L$-smoothness without any Gaussian smoothing.

\medskip

Our approach departs from these by embedding Mirror Descent within a DEQ framework. While DEQs have been previously studied under Gaussian fidelities using contractive operators such as (proximal) Gradient Descent or ADMM \cite{Willett}, their extension to the non-Euclidean and non-contractive MD case with a deep, non-convex regularization functional remains unexplored. This is precisely the focus of our work: providing a convergent and trainable formulation of DEQs for Poisson inverse problems grounded in Bregman geometry, enabling a fully implicit and less sensitive to parameters tuning scheme w.r.t. to PnP approaches.

\medskip
\paragraph{Contribution}

This work introduces a novel DEQ framework for solving linear Poisson inverse problems in imaging by learning regularization functionals directly from data. We begin in Section~\ref{sec:DEQ} by recalling the general DEQ formulation and its interpretation as a fixed-point iteration. In Section~\ref{sec:Poisson_DEQ}, we propose a DEQ architecture built upon Mirror Descent (MD), specifically tailored to minimize the Kullback–Leibler divergence \eqref{KL equation}. This choice allows us to bypass the standard \( L \)-smoothness assumption by instead relying on a relative smoothness condition with respect to a non-Euclidean geometry.

Section~\ref{eq:DEQ-MD_convergence} presents our main theoretical contribution: despite the lack of strong convexity of the KL divergence and without any convexity assumption on the learnable regularization functional, we derive sufficient conditions to guarantee convergence of the proposed fixed-point DEQ-MD reconstruction scheme. To this end, we extend the result \cite[Theorem 3.1]{Bolte2007} on the validity of the Kurdyka–Łojasiewicz property to subanalytic functions with not-necessarily closed domain. This is indeed a crucial property which allows us to extend existing DEQ analyses beyond strongly convex data fidelities and contractive operators.

Implementation-related aspects are discussed in Sections~\ref{sec:fw_bw} and~\ref{sec:training}, where we address the computational challenges associated with both the forward and backward passes of the fixed-point layer. These include, for instance, backtracking strategies for step-size selection and Jacobian-free strategies for efficient and memory-saving training.

Finally, in Section~\ref{numerical experiments}, we validate our approach on several imaging tasks. Numerical results show improved reconstruction quality over classical model-based methods and competitive performance compared to recent Bregman-PnP algorithms, while offering a less complex architecture, reduced sensitivity to hyper-parameter tuning, and full parameter-freedom at inference time, provided the level of noise in the image is known at test time. As an outlook, some encouraging numerical tests addressing generalization properties (to different tasks and to different noise models than the one use during training) are reported. Section \ref{sec:conclusions} contains some concluding remarks and perspectives.

\smallskip

For the sake of readability, detailed technical definitions and results are deferred to Appendix \ref{appendix: definitions and results}, whereas the main statements and principal results are presented in the main body of the paper.
%Moreover, by directly incorporating the forward operator \( A \) into the learning loop, our method provides increased interpretability compared to black-box denoisers used in PnP approaches.

\section{Deep Equilibrium Models (DEQs)} \label{sec:DEQ}

In this section, we provide a self-contained introduction to DEQs \cite{Primi_DEQ}. Roughly speaking, they can be defined as infinite-depth neural networks with fixed points. More precisely, a DEQ is defined by considering a single operator layer $ f_{\theta} : X \subseteq \mathbb{R}^n \to X $. A DEQ can be viewed as the infinite composition of  $f_{\theta} $, where the same parameters $ \theta \in \Theta \subset \mathbb{R}^D $ are reused at each step (a property known as weight tying).
The operator $f_{\theta}$ thus serves as the fundamental building block of a DEQ and can itself be implemented as a neural network.
% by considering an operator} \sout{single layer} \( f_{\theta} : \mathbb{R}^n \to \mathbb{R}^n \), \textcolor{red}{a DEQ is the infinite composition of $f_{\theta}$} considering the same weights \( \theta \in \Theta \subset \mathbb{R}^D \) (weight tying) \textcolor{red}{for each term in the composition. The DEQ layer $f_{\theta}$ is thus the fundamental building block of a DEQ and, potentially, it can itself be defined in terms of a neural network. 
Equivalently, a DEQ can be interpreted as the fixed point operator that maps each input $x^0$ to one of the fixed points of the DEQ-layer $f_{\theta}$. To make these notions more precise, we will study the convergence properties of DEQs in terms of existence of fixed points for $f_\theta$. To this end, we introduce the following notion of well-posedness.

\begin{definition}[Well-posedness of a DEQ] \label{well posedness of deq layer}
    Let $X \subset \mathbb{R}^n$ be an open set and let $f_{\theta}: X \to X$ with $\theta \in \Theta \subset \mathbb{R}^D$ be a DEQ layer. We say that the DEQ associated with $f_\theta$ is well-posed if:
    \begin{enumerate}
        \item The set of fixed points $\operatorname{Fix}(f_{\theta})=\{x^\infty \in X : f_{\theta}(x^\infty)=x^\infty\} \neq \emptyset$;
        \item For all initialization points $x^0 \in X$, the sequence $\{x^k\}_{k \in \mathbb{N}}$, defined by
\[
x^{k} :=
\begin{cases}
x^0 & \text{if } k=0 \\
f_{\theta}(x^{k-1}) & \text{if } k > 0,
\end{cases}
\]
converges to $x^{\infty} \in \operatorname{Fix}(f_{\theta})$.
    \end{enumerate}

\end{definition}
Designing a well-posed DEQ means constructing an operator \( f_\theta : X \to X \) that admits at least one fixed point and ensures convergence of its associated iteration. A class of operators which defines naturally a well-posed DEQ is the set of contractive operators, which ensures both existence and uniqueness of the fixed point as well as convergence of fixed point iteration schemes.
 \subsection{DEQs for image reconstruction}
In the context of image reconstruction, given a noisy measurement $ y \in \mathbb{R}^m $, the goal is to design a DEQ layer  $f_\theta(\cdot; y)$ such that its fixed points $x^\infty = f_\theta(x^\infty; y) $ provides a good approximation of the unknown ground truth $x$. 

When $ f_\theta $ is parameterized by a neural network, training a DEQ thus amounts to learning the weights $ \theta$ that minimize the reconstruction error. Note that since $f_\theta(\cdot; y) $ maps the input space to itself, computing one of its fixed point requires an initial guess $x^0 \in \mathbb{R}^n$, which can typically be chosen as a simple function of $y$, such as $x^0 = A^\top y $. While convergence is ensured under strong assumptions like contractiveness (see \cite{Willett} and Section \ref{gaussian case}), the choice of $x^0$ becomes critical in more general settings where multiple fixed points may exist. Differently from contractive networks, which are independent on the choice of the initialization, this choice becomes critical in settings where multiple fixed points may exist. Note, however, that enforcing a network to be a contraction may limit its expressivity. %Thus drawing a proper DEQ may be challenging

% For an operator $f_{\theta}:\mathbb{R}^n \to \mathbb{R}^n$ having fixed point and its associated fixed-point iteration being convergent is an hard condition to satisfy. Thus drawing these types of models is challenging. To understand how to get such models we focus on the image reconstruction task, where the goal is to build a DEQ to reconstruct an input source, given its noisy measurement (as in our problem of interest \eqref{PIP}). In particular given the pair $(x,y)$, where $x\in \mathbb{R}^n$ is the input and $y \in \mathbb{R}^m$ its associated measurement, we want to draw an operator $f_{\theta}(\cdot,y)$ such that, one of its fixed points $x^{\infty}=f_\theta(x^{\infty};y)$ is a good approximation of $x$ (meaning that $\|x^{\infty}-x\|\leq \epsilon$ for some $\epsilon > 0$). Firstly we underline that $f_{\theta}(\cdot;y)$ is a mapping from the input space to the input space, thus given the measurement $y$, to get an output $x^{\infty}$, we need an initial point $x^0$ belonging to the input space (e.g. $x^0=A^{\top}y$) to iteratively apply our map. The choice of $x^0$ is not so important if $f_{\theta}(\cdot;y)$ satisfies some assumptions (e.g. being a contraction, that has unique fixed point and each starting point converges to the same limit), but it is crucial if there exists more than one fixed point. In addition, once the structure of $f_{\theta}(\cdot;y)$ is fixed (e.g. a CNN, MLP...), we need to find the optimal weights $\theta$ for our task of interest. 

\medskip 
To formalize the DEQ training procedure, consider a dataset of image pairs \( \{(x_i^*, y_i)\}_{i=1}^n \), where each \( x_i^* \in \mathbb{R}^n \) is a ground truth image and \( y_i \in \mathbb{R}^m \) is its corresponding degraded and noisy observation. The goal is to learn parameters \( \hat{\theta} \in \Theta \) such that
\begin{equation}
    \hat{\theta} \in \underset{\theta \in \Theta}{\operatorname{argmin}}~ \mathcal{L}(\theta):=\frac{1}{n}\sum_{i=1}^{n} \ell(x^{\infty}(\theta;y_i,x^0_i),x_i^*), \label{Training eq.}
\end{equation}
where \( x^0_i \in \mathbb{R}^n \) is the $i$-th initialization, \( x^{\infty}(\theta;y_i,x^0_i) \in \mathrm{Fix}(f_{\theta}(\cdot; y_i)) \) is the limit of the fixed point iteration in Definition~\ref{well posedness of deq layer}, and \( \ell : \mathbb{R}^n \times \mathbb{R}^n \to \mathbb{R}_{\geq 0} \) is a differentiable loss function (e.g., the square loss). In the following, we describe how the existence of such fixed points can be guaranteed, and how they can be computed in practice.
 We emphasize that the DEQ model is used here in a specific way. The main objective is to learn the solution map of an inverse problem, mapping noisy measurements to the ground truth. This is achieved through the composition: $y \mapsto f_{\theta}(\cdot; y) \mapsto x^{\infty}(\theta; y, x^0)$.

    %The goal of this training phase is to find optimal weights such that, given $y$, a fixed point $x^{\infty}(\theta;y,x^0)$ of $f_{\theta}(\cdot;y)$ matches $x^*$.

\subsubsection{Gaussian case, DE-GRAD}
\label{gaussian case}

% \begin{equation}
%   \text{find }x\in\mathbb{R}^n\quad\text{such that}\quad  y= Ax  + \eta, \label{GIP}
% \end{equation}
% where $\eta \sim N(0,\sigma^2I)$ and $\sigma^2 >0$. 
%To retrieve $x \in \mathbb{R}^n$, given $y \in \mathbb{R}^m$, it is often considered the following minimization problem:
A simple example of \( f_\theta \) arises from the formulation of a classical regularized inverse problem under additive white Gaussian noise, which reads:
\begin{equation}
    \underset{x \in \mathbb{R}^n}{\operatorname{argmin}}~ \frac{1}{2}\|Ax - y\|_2^2 + \lambda R(x), \label{gauss ip}
\end{equation}
where \( A \in \mathbb{R}^{n \times n} \) is a square matrix, \( R : \mathbb{R}^n \to \mathbb{R}_{\geq 0} \) is a $L$-smooth regularization functional, and \( \lambda > 0 \) is a regularization parameter. In \cite{Willett}, the authors considered a DEQ model, based on applying the gradient descent algorithm to solve \eqref{gauss ip}. The iteration reads:
\begin{equation}
    x^{k+1} = x^k - \tau A^{\top}(Ax^k - y) - \tau \lambda \nabla R(x^k),
\end{equation}
for a suitable step size \( \tau > 0 \). 

To define the operator \( f_\theta(\cdot; y): \mathbb{R}^n \to \mathbb{R}^n \), the authors replaced the regularization term \( \lambda \nabla R \) with a neural network \( r_\theta : \mathbb{R}^n \to \mathbb{R}^n \), yielding:
\begin{equation}  \label{eq:DEQ_GP}
    f_{\theta}(x; y) := x - \tau A^{\top}(Ax - y) - \tau r_{\theta}(x).
\end{equation}
Here, the regularization parameter is implicitly encoded in \( r_\theta \), and the weights of \( f_\theta \) are shared with those of \( r_\theta \).

To ensure well-posedness of \eqref{eq:DEQ_GP}, the authors imposed conditions on \( r_\theta \) and $A$ to guarantee that \( f_\theta \) is a contraction, which ensures existence and uniqueness of a fixed point. In particular, if the operator \( r_\theta - \mathrm{Id} \) is Lipschitz-continuous with sufficiently small Lipschitz constant, then \( f_\theta \) becomes a contraction \cite[Theorem 1]{Willett}. This result relies on the fact that the data fidelity term \( \|Ax - y\|_2^2 \) is strongly convex when \( A \) is injective. However, injectivity is often not satisfied in practical inverse problems, such as downsampling or inpainting, where \( A \) is typically rank-deficient. Moreover, when the data fidelity term is no longer quadratic -- as in the case of the Kullback–Leibler divergence \eqref{KL equation} -- standard strong convexity arguments no longer apply, and more sophisticated algorithmic strategies and theoretical analyses are required to establish convergence guarantees. In this spirit, the authors in \cite{DEQKamilov} introduced a DEQ framework for Gaussian problems where convergence guarantees are established without relying on contraction properties, but based on fully optimization-based arguments. Building on this idea, we consider, in the following, optimization algorithms tailored to minimize the Kullback–Leibler divergence, and we derive fixed-point convergence results by leveraging the theoretical foundations of the underlying optimization algorithm.

\section{Poisson case, Deep Equilibrium Mirror Descent (DEQ-MD)}  \label{sec:Poisson_DEQ}

Following \cite{Willett}, we leverage an optimization algorithm well-suited to the minimization of KL \eqref{KL equation} to be used as a building block of a DEQ strategy. Note that such functional is constrained over the non-negative orthant and is not $L$-smooth, thus classical optimization schemes (e.g. Proximal Gradient Descent) do not apply in a straightforward manner.

\subsection{Bregman Proximal Gradient Descent} \label{sec:BPG}
In this section and for the following analysis, given $S \subset \mathbb{R}^n$, we will use the notation $\operatorname{int}(S)$ and $\overline{S}$ to indicate, respectively, its interior and its closure. 
The Bregman Proximal Gradient (BPG) algorithm \cite{nemirovskii_yudin_1983,nolip,Bolte2018} combined with generalized notions of smoothness, has been widely used in the literature to go beyond standard convergence assumptions of first-order optimization schemes. This algorithm can be seen as a generalization of Proximal Gradient Descent under a different geometry, potentially non-Euclidean, naturally encoding domain constraints and allowing for a weaker notion of smoothness.
In more detail, given:
\begin{itemize}
    \item a convex function $f:\mathbb{R}^n \to \mathbb{R} \cup \{+\infty\}$, differentiable on the interior of its domain,
    \item a proper, convex and lower semicontinuous (l.s.c.) function $\mathcal{R}:\mathbb{R}^n \to \mathbb{R}\cup \{+\infty\}$,
    \item a Legendre function \ref{def:Legendre function} $h:\mathbb{R}^n \to \mathbb{R}\cup \{+\infty\}$,
\end{itemize}
for $k \geq 0$, and $x^0 \in \operatorname{int (dom}(h))$, the BPG iteration can be written as:
\begin{equation}
  x^{k+1}=\operatorname*{argmin}_{x\in \mathbb{R}^n} ~ \mathcal{R}(x) + \langle x - x^k, \nabla f(x^k) \rangle + \frac{1}{\tau} D_h(x,x^k) , \label{implicit BPG}
\end{equation}
where $D_{h}$ is the Bregman divergence associated with $h$ defined in \ref{bregman distance}
and $\tau >0$ is the step-size.
Under certain compatibility conditions between $f,\mathcal{R}$ and $h$ \cite[Lemma 2]{nolip}), it is possible to get an explicit formulation of \eqref{implicit BPG}, which reads:
\begin{equation}
    x^{k+1} = \operatorname{Prox}^{h}_{\tau \mathcal{R}}(\nabla h^{*}(\nabla h(x^{k}) - \tau \nabla f(x^k))), \label{explicit MD}
\end{equation}
where $h^*$ denotes the Fenchel conjugate~\cite{Rockafellar} of the function $h$ for which there holds $\nabla h^*=(\nabla h)^{-1}$ on $\operatorname{int}(\operatorname{dom}(h^*))$ and $\operatorname{Prox}_{\tau \mathcal{R}}^h$ is the Bregman proximal operator of $\mathcal{R}$ w.r.t. $h$, defined in \ref{BProx}. 
The choice of the Bregman potential function $h$ is crucial and it is related to the functional we aim to minimize. In the case of the KL divergence \eqref{KL equation}, a tailored choice for $h$ is the so-called Burg's entropy \cite{nolip}:
\begin{equation}
    h(x)=-\sum_{i=1}^n \log (x_i), \label{burg's entropy}
\end{equation}
which enforces naturally the positivity of the desired solution.
% Equipped with BPG we aim to build a DEQ building block $f_{\theta}$. The need of sticking with such algorithm is inherently related to the particular choice of KL as data fidelity term. 
In this particular case and assuming that $R_{\theta}: \mathbb{R}^n \to \mathbb{R}$ is a regularization functional parameterized in terms of a differentiable neural network with weights $\theta \in \Theta$, we consider the following operator as the underlying DEQ layer:
\begin{equation}
    f_{\theta}(x;y):=\Pi_{[0,a]^n}(\nabla h^{*}(\nabla h (x)-\tau( \nabla \text{KL}(y,Ax) + \nabla R_{\theta}(x))),\label{f_theta of DEMD} \tag{DEQ-MD}
\end{equation}
where $h$ is the Burg's entropy defined in \eqref{burg's entropy}, $a>0$ and $\Pi_{[0,a]^n}$ denotes the projection operator on $[0,a]^n$. Note that \eqref{f_theta of DEMD} corresponds to a step of BPG minimizing a composite functional $\Psi:\mathbb{R}^n\to \mathbb{R}\cup \left\{ + \infty\right\}$ defined by:
\begin{equation}
    \Psi(x):= \text{KL}(y, Ax) + R_{\theta}(x) + \iota_{[0,a]^n}(x). \label{Psi fun}
\end{equation}
Due to the presence of $R_{\theta}$ and the indicator function  $\iota_{[0,a]^n}$, the functional is \eqref{Psi fun} is non-convex and non-smooth. Such an indicator function is required for technical reasons which will be clarified in the following.
In the case where non-smoothness is present only in the form of an indicator function $\iota_C$ of a closed and convex set $C\subset \overline{\text{dom}(h)}$, BPG corresponds to a Projected Mirror Descent \footnote{For a general $h$, we have that  $\operatorname{Prox}_{\iota_C}^h = \Pi^h_{C}$, the generalized Bregman projection.
%Given $C \subset \overline{\operatorname{dom}(h)}$ and $y \in \operatorname{int(dom}(h)), \Pi^h_{C}(y):= \underset{x \in C}{\operatorname{argmin}} \  D_{h}(x,y)$. 
In the particular case of $h$ being Burg's entropy, it can be proved that the Bregman projection on the cube coincides with the Euclidean projection (see Appendix \ref{Example}).}. For this reason, with some abuse of language, we name the proposed approach Deep Equilibrium Mirror Descent and refer to it as DEQ-MD in the following. 

Note, that we decided to parameterize directly the regularization function and not its gradient as in Section \eqref{gaussian case}, since this choice guarantees the existence of an explicit functional \eqref{Psi fun} being minimized, which allows for using tools from non-convex optimization theory to derive the well-posedness of $f_\theta$ as we do in the following Section. Furthermore, parameterizing directly $r_\theta:\mathbb{R}^n \to \mathbb{R}^n$, instead of considering $\nabla R_{\theta}$ in \eqref{f_theta of DEMD}, would require to find contractive-type bounds making use of the norm of $f_{\theta}(\cdot;y)$, which presents further challenges than in the Gaussian case due to the presence of $h$. 

\section{Convergence theory} \label{eq:DEQ-MD_convergence}

%In this section, we consider $A \in \mathbb{R}^{m \times n}_{\geq 0}$ with rows $a_i$ such that $a_i \in \mathbb{R}^n_{\geq 0} \setminus \{0_{\mathbb{R}^n}\} \ \forall i =1,...,m.$ 
In this section, we prove  the well-posedness of the DEQ-MD operator $f_{\theta}(\cdot;y)$ defined in \eqref{f_theta of DEMD} in the sense of Definition \ref{well posedness of deq layer}. To this end, we introduce a non-trivial generalization of \cite[Theorem 2]{Bolte2007}, which, to the best of our knowledge, has not been explicitly addressed in the existing literature. %Given the deep parameterization of $R_{\theta}$, we focus in particular on which conditions need to hold to have a well-defined DEQ layer. 
Having in mind the structure of the functional  $\Psi$ in \eqref{Psi fun} to be minimized,  we first study conditions upon which the MD sequence converges to some point $x^{\infty}$, then we prove that $x^{\infty} \in \operatorname{Fix}(f_{\theta}(\cdot;y)).$

Standard convergence results for first-order optimization schemes rely on the Lipschitz continuity of the gradient of the smooth part (KL$(y,A\cdot)+R_{\theta}$ in our case). This condition is not met in our framework, due to the presence of $\text{KL}(y,A\cdot)$. We thus consider a more general smoothness condition for KL expressed in terms of the Bregman potential $h$, which was first considered in \cite{nolip} in the convex setting  and then extended in \cite{Bolte2018} to the non-convex case to prove convergence of BPG: 
%\footnote{This condition is also known as $L$-SMAD \cite{Bolte2018} and Relative-Smoothness \cite{RelativeSmoothness}.},  

\begin{definition}[NoLip, $L$-SMAD, Relative smoothness \cite{nolip},\cite{Bolte2018},\cite{RelativeSmoothness}] \label{NoLip}
Let $f:\mathbb{R}^n \rightarrow \mathbb{R} \cup \{+\infty\}$ be proper, lower semicontinuous with $\operatorname{dom}(f) \supset \operatorname{dom}(h)$ and differentiable on $\operatorname{int}(\operatorname{dom}(h))$.
We say that $f$ satisfies the (NoLip) condition with constant $L>0$, if the following holds:
\begin{equation}
   \text{There exists } L>0 \text{ such that } Lh-f \text{ is convex on $\operatorname{int}(\operatorname{dom}(h))$} \label{eq nolip}   \tag{NoLip}
\end{equation}
\end{definition}

Such condition generalizes the standard notion of $L$-smoothness to the geometry induced by $h$. 
Note indeed that if $h=\frac{1}{2}\|\cdot\|_2^2$ then \eqref{eq nolip} is equivalent to the standard $L$-smoothness assumption. We have indeed that the following equivalence holds true \cite{Villa}:
\[ \frac{L}{2}\|x\|_2^2-f \text{ is convex on } \mathbb{R}^n \iff \|\nabla f(x)-\nabla f (y)\|_2 \leq L\|x-y\|_2 \quad \forall x,y \in \mathbb{R}^n.\]

We now recall for convenience  \cite[Lemma 7]{nolip}. It provides the existence of such  constant $L$ for the pair $(\text{KL}(y,A\cdot),h)$, justifying the choice of the specific Bregman potential considered.

\begin{lemma} [\cite{nolip} Lemma~7]\label{NoLip constant of KL}
Let $f: x \mapsto \text{KL}(y,Ax)$ and $h$ be the Burg's entropy, defined in \eqref{burg's entropy}. Then, for any $L$ satisfying
$
L \geq \|y\|_1 = \sum_{i=1}^{m} |y_i|,
$
the function $Lh - f$ is convex on $\mathbb{R}_{>0}^n=\operatorname{int}(\operatorname{dom}(h))$, i.e \eqref{eq nolip} holds.
\end{lemma}

By Lemma~\ref{NoLip constant of KL}, to establish the relative smoothness condition \eqref{eq nolip} for the full objective \(\mathrm{KL}(y, A\cdot) + R_{\theta}\), it thus suffices to verify condition \eqref{eq nolip} for the regularization term \( R_{\theta} \). Specifically, if \eqref{eq nolip} holds for the pair \((R_{\theta}, h)\) with constant \( L' \), then it holds for \((\mathrm{KL}(y, A\cdot) + R_{\theta}, h)\) with any \( L \geq L' + \|y\|_1 \), since the sum of convex functions is still convex (see Appendix \ref{proof on nolip} for details).

As noted by Hurault et al.~\cite{Hurault}, verifying condition \eqref{eq nolip} on the entire interior of \(\operatorname{dom}(h)\), which is unbounded, can be challenging. However, this becomes significantly easier when restricting the analysis to a compact subset. The key insight is that the function \( h \) becomes strongly convex when restricted to bounded subsets of its domain. If the Hessian \(\nabla^2 R_{\theta}\) is also bounded on such a compact set, one can derive explicit bounds that ensure \eqref{eq nolip} to hold (see again Appendix \ref{proof on nolip}).
For this reason, we include an indicator function in the definition of the objective \eqref{Psi fun}, effectively restricting the domain to a compact set. This assumption is not restrictive in practice as normalized image data naturally lie in compact sets such as, e.g., \([0,1]^n\).

The other key ingredient is the Kurdyka–Łojasiewicz property \ref{def:KŁ}, a classical and powerful tool used to prove algorithmic convergence in non-convex optimization.
Although this property may be difficult to be checked directly, it turns out that is satisfied by a wide class of functions. For instance, subanalytic functions continuous on a closed domain satisfy this property \cite{Bolte2007}. However, the domain of $\Psi$ \eqref{Psi fun} is in general neither closed nor open (see Appendix \ref{proof of prop 4.4}). We overcome this difficulty by generalizing \cite[Theorem 3.1]{Bolte2007}. For more details we refer to the Appendix \ref{proof of prop 4.4}).
\subsection{Kurdyka-Łojasiewicz property for subanalytic functions with non-closed domain}
In this subsection, we show that property \ref{def:KŁ} still holds for subanalytic and coercive functions defined on a non-closed domain. The intuition behind the result is that even when the domain is not closed, if it is possible to separate the critical points of the functional being minimized from the set where the functional is infinite, then \ref{def:KŁ} can still be recovered.
\begin{proposition}
[KŁ with non-closed domain] \label{KL with not-closed domain}
    Let $\Psi:\mathbb{R}^n \to \mathbb{R} \cup \{+\infty\} $ be a lower semicontinuous and subanalytic function. If $\Psi$ is coercive, $\Psi \!\mid_{\operatorname{dom} \Psi}$ is continuous and such that  
    \begin{equation} \label{eq:S}
    S:=\{x \in \overline{\operatorname{dom}(\Psi)} \mid \Psi(x)= +\infty\}
    \end{equation}
    is closed, then $\Psi$ satisfies thr Kurdyka–Lojasiewicz property \ref{def:KŁ} at any $x^*\in\text{crit}(\Psi)$.
%     there exists an exponent $\theta \in [0,1)$ such that the function
% \begin{equation} \label{eq:8}
%     \frac{|\Psi - \Psi(x^*)|^{\theta}}{m_{\Psi}}
% \end{equation}
% is bounded around $x^*$.

\end{proposition}

\begin{proof}
The proof is available in Appendix \ref{proof of prop 4.4}.
\end{proof}
Thanks to Proposition \ref{KL with not-closed domain}, we can then state a key result for the regularized Kullback-Leibler functional with regularization.
\begin{corollary} \label{Kl corollary}
    Given $\Psi(\cdot):= \textrm{KL}(y,A\cdot) + \iota_{[0,a]^n}+ R_{\theta}$, if $R_{\theta}:\mathbb{R}^n \to \mathbb{R}$ is analytic, then $\Psi$ satisfies the Kurdyka-Łojasiewicz property \eqref{def:KŁ}.
\end{corollary}
\begin{proof}
The proof is available in Appendix \ref{proof of corollary 4.5}.
\end{proof}

\subsection{Convergence results of DEQ-MD} \label{sec: conv result}
 Thanks to the technical results reported in the previous subsection,  
we are now ready to state our main theoretical result on the convergence of DEQ-MD.

\begin{proposition} [Convergence of DEQ-MD with learnable regularization] \label{prop: convergence}
\leavevmode\\ 
Let $R_{\theta}: \mathbb{R}^n \to \mathbb{R}$ be an analytic neural network parameterized by parameters $\theta\in\Theta$ and let $h$ be defined as in \eqref{burg's entropy}.  Let $L>0$ be a constant satisfying \eqref{eq nolip} for the pair $(\text{KL}(y,A\cdot)+R_{\theta},h)$, let $a > 0$, and let  $x^0 \in \operatorname{int(dom}(h))$. Then, for every $\tau \in (0, \frac{1}{L})$, the sequence given by:
    \begin{equation}
    x^{k+1}\in\operatorname*{argmin}_{x\in \mathbb{R}^n} ~ \iota_{[0,a]^n}(x)+ \langle x - x^k, \nabla \left(\text{KL}(y,A\cdot)+R_{\theta}(\cdot)\right)(x^k) \rangle + \frac{1}{\tau} D_{h}(x,x^k) \label{PGD for KL} 
\end{equation}
    converges to $x^{\infty}(\theta;y,x^0)$, a critical point of the functional $\Psi$ defined in \eqref{Psi fun}. In addition, \eqref{PGD for KL} can be expressed in closed-form as:
\begin{align}
x^{k+1} &= \Pi_{[0,a]^n}\Big(\nabla h^{*}(\nabla h (x^k)-\tau( \nabla \text{KL}(y,Ax^k) + \nabla R_{\theta}(x^k)) \Big) \notag \\
&= \Pi_{[0,a]^n}\left(\frac{x^k}{1+\tau x^k (\nabla \text{KL}(y,Ax^k)+\nabla R_{\theta}(x^k))}\right).
\label{f_theta as iteration}
\end{align}

\end{proposition}

\begin{proof}
The proof is available in Appendix \ref{proof of prop 4.6}.
\end{proof}

Note that the operator in \eqref{f_theta as iteration} corresponds to the MD operator \( f_{\theta}(\cdot; y) \) defined in \\ \eqref{f_theta of DEMD}.

\begin{remark}[Analyticity is not restrictive]
   Although analyticity might appear to be a strong assumption, it is often satisfied in practice. Neural networks typically require differentiable activation functions to enable gradient-based training, and many commonly used smooth activations (e.g., SoftPlus, tanh, sigmoid) are analytic. As a result, the entire network -- hence the learnable regularizer \( R_{\theta} \) -- often inherits analyticity as a natural consequence of architectural choices (see Appendix \ref{proof of corollary 4.5}, paragraph 1). 
\end{remark}

\begin{remark}[On the choice of $h$]
  The choice of $h$ as the Burg's entropy is a standard choice in the literature \cite{nolip}. The rationale behind this choice lies on the compatibility between the domain of $h$ and the one of  $\text{KL}(y,A\cdot)$, that is $\operatorname{dom}(h) \subset \operatorname{dom}(\text{KL}(y,A\cdot ))$. Moreover, upon this choice, it is possible to find an explicit constant $L>0$ satisfying \eqref{NoLip} (Lemma \ref{NoLip constant of KL}). However, certainly this is not the only possible choice guaranteeing the conclusion of Proposition \ref{prop: convergence} as the two main properties required to $h$ are the compatibility condition as well as the strong convexity of $h$ on bounded sets. 
  %So, any $h$ satisfying these property could be used.
\end{remark}

We conclude this section proving the well-posedness of the DEQ-MD  $f_{\theta}(\cdot;y)$ in \eqref{f_theta of DEMD}.

\begin{proposition} [Well-posedness of DEQ-MD]\label{prop:well_posedness_DEQ_MD}

Under the same assumptions of Proposition \ref{prop: convergence}, the DEQ-MD model \eqref{f_theta of DEMD} is well-posed in the sense of Definition \ref{well posedness of deq layer}.
\end{proposition}
\begin{proof}

By Proposition \ref{prop: convergence}, for all $x^0 \in \operatorname{int(dom(h))}$, the sequence $\{x^k\}_{k \in \mathbb{N}}$ converges to $x^{\infty}\in\text{crit}(\Psi)\subset \operatorname{dom}(\Psi).$ %\cite{rockafellar1998variational}
Recalling the explicit expression for $f_{\theta}(\cdot;y)$ in \eqref{f_theta as iteration},
% $$
%  f_{\theta}(x;y)=\Pi_{[0,a]^n}\left(\frac{x}{1+\tau x (\nabla \text{KL}(y,Ax)+\nabla R_{\theta}(x))}\right),
% $$
and observing that $x^k\in \operatorname{int}(\operatorname{dom} (\text{KL}(y,A\cdot)))$, we get that:
\[
x^{\infty}=\lim_{k\to\infty}~ x^{k+1} = \lim_{k\to\infty}~ f_{\theta}(x^k;y)=f_{\theta}(\lim_{k\to\infty} x^k;y)=f_{\theta}(x^{\infty};y).
\] 
\end{proof}

\section{DEQ-MD implementation}  \label{sec:fw_bw}

Training a DEQ poses both several theoretical and practical challenges: first a grounded approach for computing fixed points is required; moreover, computing gradients through infinitely many applications of \( f_\theta \) is memory-intensive, making standard backpropagation infeasible.
More precisely, we focus in this Section on these two crucial steps:
\begin{itemize}
    \item \textbf{Forward pass}: it refers to how a fixed point $x^{\infty}(\theta;y,x^0)$ of $f_{\theta}(\cdot;y)$, given an initial point $x^0$, is computed.
    \item \textbf{Backward pass}: it refers to the computation of the gradients of $L(\theta)$ \eqref{Training eq.} with respect to its weights $\theta$.
\end{itemize}

\subsection{Forward pass}
Having defined the DEQ-MD operator $f_{\theta}(\cdot;y)$ as in \eqref{f_theta of DEMD}, we need an effective way to compute a fixed point $x^{\infty}(\theta;y,x^0)$ given an image $y$ and weights $\theta$. 

The fixed points of $f_{\theta}(\cdot;y)$ are in general not unique since they are critical points of a non convex function. On the other hand, Proposition \ref{prop: convergence} states that convergence holds for any choice of $x^0 \in \operatorname{int(dom}(h))$. Thus, for any such $x^{0}$, the solution map $\theta \mapsto x^{\infty}(\theta;y,x^0)$ is single-valued. 

Another key aspect is the choice of the step-size $\tau>0$ guaranteeing convergence. From Proposition \ref{prop: convergence}, convergence holds for $\tau \in (0,\frac{1}{L})$, with $L\geq \|y\|_1+L'$ chosen to satisfy \eqref{eq nolip} for the pair  $(\text{KL}(y,A\cdot) + R_{\theta},h)$ and where $L'$ is the NoLip constant for  $(R_{\theta},h)$. Such $L'$ exists (see Appendix \ref{proof on nolip}), but its computation might be challenging as it requires the computation of the maximum eigenvalue of the Hessian $\nabla^2 R_{\theta}$, which may be hard to compute. During training, the update of the weights $\theta$ would lead indeed to a different instance  of $R_{\theta}$, leading to the need of computing the maximum eigenvalue of $\nabla^2 R_{\theta}$ each time the weights $\theta$ change.
Furthermore, for images of standard size (in Section \ref{numerical experiments} we used  $256 \times 256$) $\|y\|_1$ is big enough to make $\tau$ very small. One may argue that with a normalization of the data this may be avoided, but this would lead in fact to smaller magnitude of the gradients involved in \eqref{f_theta of DEMD}, which would balance a smaller $\|y\|_1$.  As a consequence,  we decided to perform the forward pass  not computing such $L$ explicitly, but, rather, by estimating a suitable step-size by means of a backtracking strategy, similarly as done in \cite{Hurault} for PnP approaches.  Recalling the definition of $\Psi$ in \eqref{Psi fun}, we describe such strategy in Algorithm~\ref{BT}.

\begin{algorithm}[H]
\caption{Backtracking for DEQ-MD}
\label{BT}

\KwIn{Initial $\tau_0 > 0$, parameters $\gamma \in (0,1)$, $\eta \in (0,1)$}

\textbf{Define}: $T_\tau(x) :=\Pi_{[0,a]^n}(\nabla h^{*}(\nabla h (x)-\tau( \nabla \text{KL}(y,Ax) + \nabla R_{\theta}(x))))$ \;

\textbf{Initialize}: $\tau = \tau_0$ and $x^0 \in \operatorname{int(dom(}h))$ \;

\While{$\Psi(x^k) - \Psi(T_\tau(x^k)) < \frac{\gamma}{\tau} D_h(T_\tau(x^k), x^k)$}
{$\tau \gets \eta \tau$ \;}
\textbf{Update}: $x^{k+1} \gets T_\tau(x^k)$ \;
\end{algorithm}

The convergence result stated in Proposition Proposition \ref{prop: convergence} still holds upon this choice:
\begin{corollary}
    The Backtracking procedure described in Algorithm \eqref{BT} terminates in a finite number of steps and provides a sequence of strictly positive step-sizes $\{\tau_k\}_{k\in \mathbb{N}}$. In addition, the conclusion of Proposition \ref{prop: convergence} still holds under the choice $\tau_k$ in \eqref{PGD for KL}.
\end{corollary}
\begin{proof}
    Proof in Appendix D.2 of  \cite{Hurault}.
\end{proof}

\begin{remark}
In contrast to the case of standard DEQs employed in Gaussian cases (see, e.g., \cite{Willett}), we preferred a Backtracking strategy over standard acceleration strategies for fixed-point iterations, such as Anderson acceleration \cite{Anderson}. %Such strategy linearly combines past iterates to find promising search directions for the following iteration. 
To employ Anderson acceleration, we need to take $\tau \in (0,\frac{1}{L})$ and usually this theoretical step-size is much smaller than the ones given by Backtracking (see Figure \ref{backtracking step-sizes}). As a consequence, employing backtracking typically leads to faster convergence and, in addition, cheaper computations. Note, moreover, Anderson acceleration is typically applied to accelerate fixed-point iterations exhibiting linear convergence (e.g., contractive mappings), while the iterations \eqref{f_theta as iteration} do not correspond to a fixed point iteration of a contraction map.
\end{remark}

\subsection{Backward pass} \label{Back pass}

Given the structure of $f_{\theta}(\cdot;y)$ in \eqref{f_theta of DEMD}, the minimization problem \eqref{Training eq.} needs to be solved to compute optimal weights $\hat\theta$.
This problem is usually addressed  by means of first-order optimization algorithms (e.g., ADAM \cite{Adam}), hence a tractable way to compute gradients of the loss functional $\mathcal{L}(\cdot)$ needs to be found. We provide in the following a formal computation of such gradients. For $\bar{\theta}\in\Theta$ we want to compute:
\begin{equation}  \label{eq:grad_loss_theta}
    \nabla_{\theta} \mathcal{L}(\bar{\theta})=\sum_{i=1}^{n} \nabla_{\theta} \ell(x^{\infty}(\bar{\theta};y_i,x^0_i),x_i^*).
\end{equation}

%Note, that $f_{\theta}(\cdot;y)$ in \eqref{f_theta of DEMD} is not differentiable (both w.r.t. to $x$ and $\theta$) due to the presence of the projection operator $\Pi_{[0,a]^n}$. To compute gradients and make the following derivation formally correct, we thus replace $\Pi_{[0,a]^n}$ in \eqref{f_theta of DEMD} with its smoothed version $\Pi_{(- \infty,a]^n}^{\epsilon}:\mathbb{R}^n \to \mathbb{R}^n$ defined by $\Pi_{(- \infty,a]^n}^{\epsilon}(x)= \sum_{i=1}^n \Pi_{(- \infty,a]}^{\epsilon}(x_i)$, which, given $\epsilon>0$, writes as: 
% \begin{equation}
%     \Pi_{(- \infty,a]}^{\epsilon}(x_i)= x_i - \epsilon \log \left(1 + \text{exp}\left(\frac{1}{\epsilon}(x_i-a)\right)\right). \label{smooth proj}
% \end{equation}
% Note that for all $ x \in \mathbb{R}$, $\Pi_{(- \infty,a]}^{\epsilon}(x) \to \Pi_{(- \infty,a]}(x)$ whenever $\epsilon \to 0^+$. 
% Further, note that we do not need to project explicitly onto $[0,a]^n$, since the non negativity constraint is naturally enforced by the the use of the Bregman potential $h$. With some abuse of language, we will refer in the following to $f_{\theta}(\cdot,y)$, as the operator given by \eqref{f_theta of DEMD} where the projection operator $\Pi_{(-\infty,a]^n}$ is replaced by $\Pi^\epsilon_{(-\infty,a]^n}$ above.

Omitting the dependence on $i=1,\ldots,n$ in \eqref{eq:grad_loss_theta}, by the chain rule we get:
\begin{equation}
    \nabla_{\theta} \ell(x^\infty(\bar{\theta};y,x^0),x^*)=\frac{\partial x^\infty(
\bar{\theta};y,x^0)}{\partial \theta}^{\top} \frac{\partial \ell(x^\infty(\bar{\theta};y,x^0))}{\partial x}, \label{gradient of ell}
\end{equation}
where by $\frac{\partial \cdot}{\partial \theta}$ we denote the  formal Jacobian matrix of the function $\theta\mapsto x^\infty(\theta;y,x^0)$ w.r.t. $\theta$ and $\frac{\partial \cdot}{\partial x}$ is the derivative of $\ell$ w.r.t..~its argument. We refer to Remark \ref{rmk:diff} for more details on how to make such derivation precise in a non-smooth setting.

%{\color{red} Note that the 
% since the differentiability of $x^{\infty}$ w.r.t.~$\theta$ makes the derivation formal (see Remark \ref{rmk:diff}).} 
In \eqref{gradient of ell}, the quantity  $\frac{\partial x^\infty(\theta;y,x^0)}{\partial \theta}$ is difficult to deal with as it is the output of the infinitely many iterations of the mapping $f_{\theta}(\cdot;y)$.
To make such computation possible, for a fixed $\bar{\theta}\in\Theta$, the quantity $x^\infty(\bar{\theta};y,x^0)$ can be computed by exploiting  the fixed point equation \\$x^\infty(\bar{\theta};y,x^0)=f_{\bar \theta}(x^\infty(\bar{\theta};y,x^0);y)$.  We detail this in the following, using the shorthand notation $x^{\infty}(\theta)$ to indicate $x^{\infty}(\theta;y,x^0)$ and $f_{\theta}(\cdot)$ to indicate $f_{\theta}(\cdot;y).$ To avoid ambiguities, we will make precise the points where the functions are evaluated explicitly. For all $\theta\in\Theta$, we define the function $F: \Theta \to \mathbb{R}^n$ by
$$
F(\theta) 
%F(x^\infty(\theta),f_{\theta}(x^\infty(\theta))) 
:=x^\infty(\theta)-f_{\theta}(x^\infty(\theta)), $$
which is identically equal to zero by Proposition \ref{prop:well_posedness_DEQ_MD}. 
We can thus write:
\begin{equation}
   \frac{\text{d} F}{\text{d} \theta}=\frac{\partial x^\infty}{\partial \theta} -  \frac{\text{d} f_{\theta}(x^\infty)}{\text{d} \theta}=0, \label{chain role on F}
\end{equation}
where $\text{d}\cdot/\text{d}\theta$ denotes the total derivative w.r.t..~$\theta$. We get:
\begin{equation}  \label{eq:chain_rule2}
\frac{\text{d} f_{\theta}(x^\infty(\theta;y,x^0))}{\text{d}  \theta}\biggr|_{\theta=\bar{\theta}}=\frac{\partial f_{\bar{\theta}}(x)}{\partial x} \biggr|_{x=x^\infty} \frac{\partial x^\infty (\theta;y,x^0)}{\partial \theta}\biggr|_{\theta=\bar{\theta}} + \frac{\partial f_{\theta}(x^\infty(\bar{\theta};y,x^0))}{\partial \theta} \biggr|_{\theta=\bar{\theta}}.
\end{equation}
From \eqref{eq:chain_rule2} and \eqref{chain role on F}, it is thus possible to get the general linear equation solved by the Jacobian of the function $\theta\mapsto x^\infty(\theta)$, which reads:
\begin{equation}
    \left(I-\frac{\partial f_{\bar{\theta}}(x)}{\partial x}\biggr|_{x=x^\infty}\right)\frac{\partial x^\infty(\bar{\theta};y,x^0)}{\partial \theta}
    % =
    % \left(I-\frac{\partial f_{\bar{\theta}}(x^\infty)}{\partial x}\right)\frac{\partial x^\infty}{\partial \theta}(\bar{\theta};y,x^0)
    = \frac{\partial f_{\bar{\theta}}(x^\infty(\bar{\theta};y,x^0))}{\partial \theta}. \label{system solved by the jacobian of x inf}
\end{equation}
Denoting now by $\|\cdot\|$ the spectral norm, if \( \left\lVert  
\frac{\partial f_{\bar{\theta}}(x^\infty)}{\partial x}\right\rVert < 1 \), then $\left(I-\frac{\partial f_{\bar{\theta}}(x^\infty)}{\partial x}\right)$ is invertible and there exists a unique solution to \eqref{system solved by the jacobian of x inf} given by:
%$\frac{\partial x^{\infty}(\theta)}{\partial \theta}$:
\begin{equation}
    \frac{\partial x^{\infty}(\bar{\theta};y,x^0)}{\partial \theta}=\left(I-\frac{\partial f_{\bar{\theta}}(x)}{\partial x}\biggr|_{x=x^\infty}\right)^{-1}\frac{\partial f_{\bar{\theta}}(x^{\infty}(\bar{\theta};y,x^0))}{\partial \theta},\quad \forall \bar{\theta}\in\Theta,
    \label{solution of the systen with the inverse}
\end{equation}
which can be plugged in \eqref{gradient of ell} to get the desired expression.

In several contexts, inverting the matrix arising in \eqref{solution of the systen with the inverse} can be computationally expensive or even infeasible, as is often the case with implicit neural networks. In the following section, we present an approximation strategy that circumvents this limitation and enables efficient computation.

\begin{remark}[Implicit differentation in the non-smooth setting] \label{rmk:diff}
Note that due to the presence of the non-smooth projection operator $\Pi_{[0,a]^n}$ in the DEQ-MD layer \eqref{f_theta of DEMD}, the previous computations relying on \eqref{gradient of ell} are only formal.
%
% The explicit calculation of gradients in \eqref{gradient of ell} requires differentiability, although this may not hold whenever non-smooth parameterizations of the underlying DEQ layer are used, see e.g.,~\cite{Primi_DEQ}. However, the validity of the chain rule is often empirically justified by the practical effectiveness of the algorithms employed, and further supported by the robustness of modern automatic differentiation routines (e.g., \texttt{autograd}) which can effectively handle piece-wise differentiable functions. 
%
% In our case, we verified numerically that choosing the smoothed projection as above for
% $\epsilon$ sufficiently small, does not alter the the fixed points of the original $f_{\theta}$ in \eqref{f_theta of DEMD}, hence we decided to replace it as described above.
% Alternatively, one could avoid the use of an explicit smoothing by  considering a different Bregman potential $h$ enforcing naturally the projection onto $[0,a]$, such as, e.g., $h(x)= - \log(x) - \log(1-x)$. However, in order to compare our results with existing approaches relying on the standard choice of the Burg's entropy (see \cite{Hurault}) we preferred to not consider this alternative, although we believe that it could be an elegant idea to be considered in future work.
%This choice of $h$ as Legendre function in the BPG algorithm would force solution to stay in $(0,a)$.
%
The appropriate framework to compute derivatives and chain rules in the non-smooth case is the one of conservative derivatives (see, e.g., \cite{BolteLePauwelsSilveti2021,conservativederivatives}), which we did not consider in our work given the more practical focus on the use of DEQ-MD in the context of image reconstruction.

%However this type of analysis might be tricky. 

%To conclude the remark we comment our choice of replacing the non-smooth projection with its smooth approximation \eqref{smooth proj}. To be completely precise we would need to show that the operator associated with the smooth projection is well-posed in the sense of Definition \eqref{well posedness of deq layer}. This is non trivial, since Proposition \eqref{prop: convergence} does not apply directly. However, we found numerically that for 

\end{remark}

\subsubsection{Jacobian-Free Backpropagation} \label{JFB}
In recent years, research on implicit networks has led to the use of a technique called Jacobian-Free Backpropagation (JFB) \cite{JFB} where the inverse of the Jacobian matrix appearing in \eqref{solution of the systen with the inverse} is replaced by the identity matrix.  A possible interpretation of JFB is a zeroth-order approximation to the Neumann series of $\frac{\partial f_{\theta}(x^{\infty})}{\partial x}$ since there holds:
\begin{equation}
    \left(I-\frac{\partial f_{\theta}(x^{\infty})}{\partial x}\right)^{-1}=\sum_{k=0}^{+\infty} \left(\frac{\partial f_{\theta}(x^{\infty})}{\partial x}\right)^k,
\end{equation}
which holds whenever $  \left\|\frac{\partial f_{\theta}(x^{\infty})}{\partial x}\right\|<1$.
Although this approximation may seem overly simplistic, it has been shown to work effectively in many practical situations. Specifically it was shown in \cite{JFB} that, upon suitable assumptions, while the gradients computed using this technique differ from the theoretical gradients, they  point towards actual descent directions for the functional $L(\cdot)$, thus allowing its decrease and enabling successful training. In \cite{Samuel}, estimates in operator norm of the difference between the true Jacobian and the JFB-Jacobian are given. 
Using this approach simplifies the calculation of the gradients of $L(\cdot)$ to:
\begin{equation}
    \nabla_{\theta} \mathcal{L}(\bar{\theta}) \approx \sum_{i=1}^{n} \frac{\partial f_{\bar{\theta}}(x^{\infty}(\bar{\theta};y_i,x^0_i);y_i)}{\partial \theta}^{T} \frac{\partial \ell(x^{\infty}(\bar{\theta};y_i,x^0_i),x_i^{*})} {\partial x},\quad \bar{\theta}\in\Theta. \label{gradient with JFB}
\end{equation}

In the case of $f_\theta$ defined \eqref{f_theta of DEMD}, we remark that there is no theoretical guarantee that \(  \|\frac{\partial f_{\theta}(x^{\infty})}{\partial x}\| < 1 \) due the lack of convexity of the functional $\Psi$. However, we found out numerically that the relation holds with a non-strict (that is, with $\leq$) inequality.
Notice that the expression in \eqref{gradient with JFB} can be computed formally in this case, although there is no theoretical guarantee about the effectiveness of the approximation in this case. Nevertheless, as shown in the next Section \ref{sec:training}, this approach still works well in practice.

\begin{remark}[DEQs in practice]
It is useful to remind here that  the infinite depth aspect of a DEQ is primarily a theoretical concept. In practice, a finite number of iterations $N>0$  in the forward pass is computed to obtain a sufficiently accurate approximation of a fixed point within a certain tolerance.
%Indeed, all the equations in this section hold when $x^{\infty}$ represents a fixed point. Thus, 
Hence, in practice, $x^{\infty}$ shall be thought of as the output of an $N$-layered network with $N$ sufficiently large. This might suggest that standard backpropagation could be used to calculate the gradients of interest in \eqref{gradient of ell}. Note,  however, that even with a small $N$, memory issues quickly arise, since $f_{\theta}$ can be arbitrarily wide. This highlights one of the advantages of the gradient calculation method presented above: regardless of the number of iterations required in the forward pass, the gradients are indeed computed with \emph{constant memory} \cite{Primi_DEQ}, meaning it is only necessary to differentiate $f_{\theta}(\cdot;y)$ at the fixed point as shown in \eqref{gradient with JFB}.
\end{remark}

\section{Numerical implementation}   \label{sec:training}
In this section, we discuss the details concerned with the image datasets, the network architecture and the  implementation of the proposed DEQ-MD model. In particular, we discuss the parameterization of the deep regularization function, for which we propose two different possibilities, and their corresponding architectures. Following this, we detail both the forward-pass and some computational details of the proposed training methodology. For reproducibility purposes, the code used for the following experiments is available at \url{https://github.com/christiandaniele/DEQ-MD}.

\subsection{Image and training details}

In this section we present the image dataset considered, the data generation procedure and the training details. For image assessment, we will use the following three standard image quality metrics: PSNR, SSIM and LPIPS (see Appendix \ref{evaluation metrics}.)

\subsubsection{Image Dataset} \label{dataset}
For all experiments, we used the Berkeley image Segmentation Dataset 500 (BSDS500), a dataset of 500 natural RGB images resized to 256x256 by cropping a window of the desired size at the center of the original image. We used 406 images for training, 48 for validation and 46 for test.
The dataset is available at \url{https://www.kaggle.com/datasets/balraj98/berkeley-segmentation-dataset-500-bsds500}.
\subsubsection{Data Generation} \label{data generation}
\begin{figure}[h!]
    \centering
    \begin{subfigure}[b]{0.28\textwidth}
        \includegraphics[width=\linewidth]{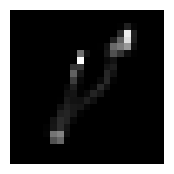}
        \caption{Motion blur kernel}
        \label{motion blur}
    \end{subfigure}
    \hfill
    \begin{subfigure}[b]{0.28\textwidth}
   
     \includegraphics[width=\linewidth]{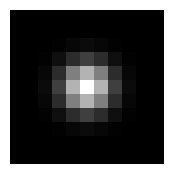}
        \caption{Gaussian Kernel} \label{gaussian kernel}
    \end{subfigure}
    \hfill
    \begin{subfigure}[b]{0.28\textwidth}
        \includegraphics[width=\linewidth]{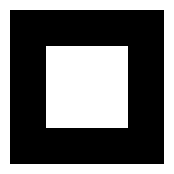}
        \caption{Uniform kernel}
        \label{uniform kernel}
    \end{subfigure}
    
    \caption{The 3 blur kernels used in the experiments. (a) is a real-world camera shake kernel, see \cite{Levin}. (b) is a \( 11 \times 11 \) Gaussian kernel with standard deviations $\sigma$ =1.2. (c) is a \( 9 \times 9 \) uniform kernel. }
    \label{fig:kernels}
\end{figure}

We generate data using the model
\begin{equation}
    y=\mathrm{Poiss}(\alpha Ax),
\end{equation}
where $A\in \mathbb{R}^{n \times n}_+$ is a convolution operator and $\alpha>0$ allows to module the intensity of Poisson noise. 
In Figure \ref{fig:kernels} we present a selection of kernels associated with the operator $A$ available within the DeepInverse library\footnote{\url{https://deepinv.github.io/deepinv/}}. We considered three levels of noise $\alpha\in\left\{100,60,40\right\}$ corresponding  to low, medium and high amount of Poisson noise and trained a DEQ-MD model for each level of noise and each kernel. 

\subsubsection{Training loss and optimizers} \label{loss function}

As a training loss we used the supervised MSE loss:
\begin{equation}
    \ell(x^*,x^{\infty})=\|x^*-x^{\infty}\|^2_2, \label{mixed loss}
\end{equation}
see \ref{training results} for more details.
As optimizer, we considered ADAM \cite{Adam} initialized with random weights and with piecewise constant step-size $\gamma=5 \times 10^{-4}$, that is halved after 25 iterations. The JFB strategy (section \ref{JFB}) was employed to approximate stochastic gradients. Training was performed over 50 epochs and optimal weights were chosen so as to maximize the PSNR on the validation set.

\subsection{Deep regularization: parameterization} \label{param of reg funz}

\begin{figure}[htbp]
    \centering
    \begin{subfigure}[b]{0.3\textwidth} % Adjust width as needed
        \includegraphics[width=\textwidth]{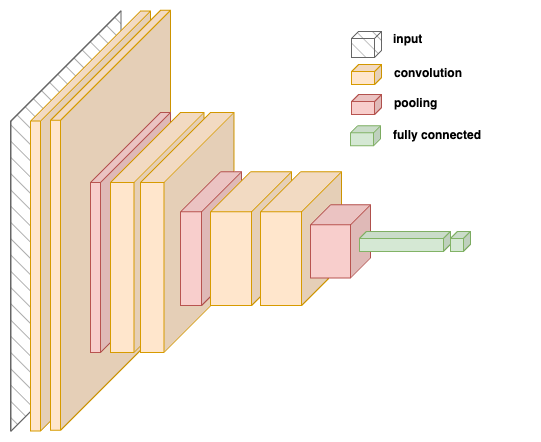} % Replace figure1.png with your first image file
        \caption{$R_{\theta}$ for DEQ-S.}
        \label{DEQ-S}
    \end{subfigure}
    \hfill % Puts space between the figures
    \begin{subfigure}[b]{0.6\textwidth} % Adjust width as needed, ensure total width is less than 1.0\textwidth
        \includegraphics[width=\textwidth]{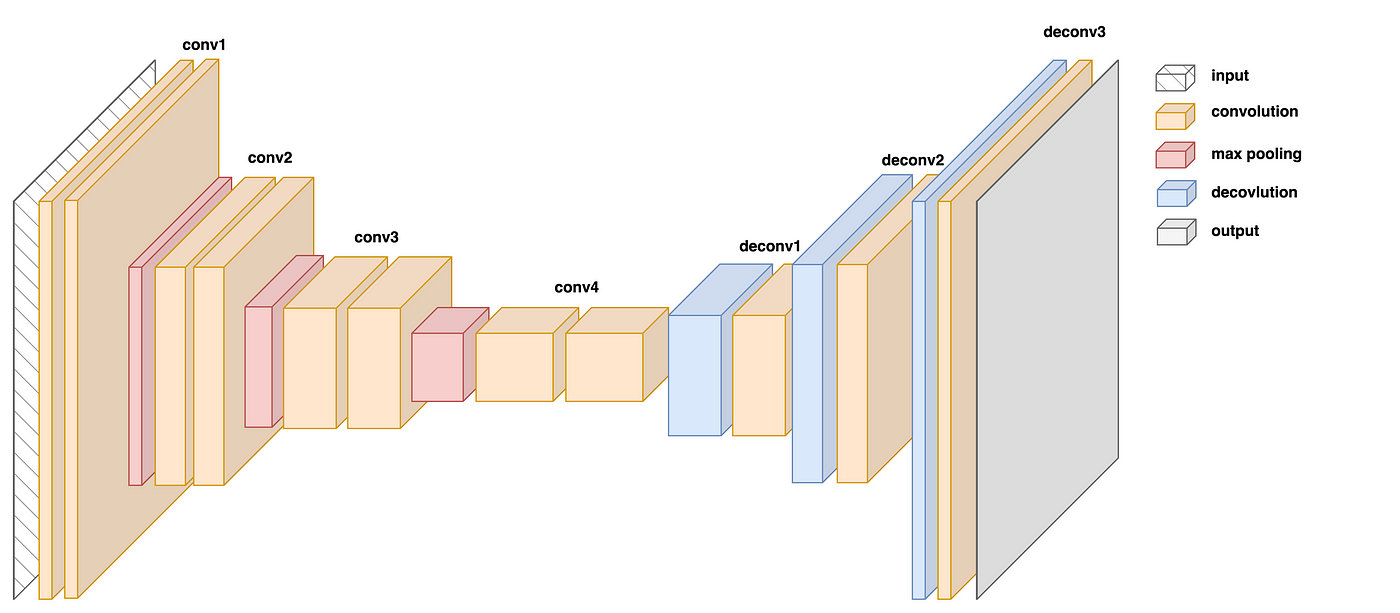} % Replace figure2.png with your second image file
        \caption{$N_{\theta}$ for DEQ-RED.}
        \label{DEQ-REDed}
    \end{subfigure}
\caption{Illustration of the two different networks employed for the parameterization  of the regularization function.}
\end{figure}

We consider the following two parameterizations of the regularization function $R_{\theta}:\mathbb{R}^n \to \mathbb{R}$ in \eqref{f_theta of DEMD}:
\begin{itemize}
    \item \textbf{Scalar Neural network (DEQ-S)} where $R_{\theta}$ is a neural network which maps an image input into a scalar value, e.g. a convolutional neural network (Figure \ref{DEQ-S}); 
    \item \textbf{RED regularization (DEQ-RED)} where $R_\theta$
   follows a Regularization by Denoising (RED) parameterization,  firstly proposed in \cite[Section 5.2]{Romano2017} and subsequently employed, e.g.,  in \cite{Hurault1} in the form
   $R_{\theta}(x)=\frac{1}{2}\|x-N_{\theta}(x)\|_2^2 $, with $N_{\theta}: \mathbb{R}^n \to \mathbb{R}^n$. The network $N_{\theta}$  can be chosen, for instance, as a Denoising CNN (DnCNN) \cite{DnCNN}, see Figure \ref{DEQ-REDed} for an illustration of the architecture. Note that contrary to plug-and-play approaches, $N_{\theta}$ is trained directly on the task of interest, that is image reconstruction under Poisson noise. 
   %This allows for the utilization of neural networks specifically designed for image reconstruction, such as DRUNet (\cite{DRUNet}) or DnCNN (\cite{DnCNN}). 
\end{itemize}

%Under the first choice, we parameterize a general regularization function to be learned through data-driven training. On the other hand, the second approach employs a regularization function with a fixed structure (RED-like).
Employing the latter strategy allows for pre-training $N_{\theta}$ on, e.g. a denoising task, which generally leads to a better initialization of the network weights. This has the advantage to accelerate convergence during training and, potentially, to achieve better critical points of the loss function (see Appendix \ref{pre-training}). 
As observed \cite{Salimans},  parameterizing directly $R_{\theta}$ via a scalar neural network may lead to  poorer performance.

%\textcolor{red}{Either case, in order for $R_{\theta}$} to be differentiable, smooth activation functions are used, see Section \ref{sec:network_choice} for more details.

\subsubsection{Network architecture}  \label{sec:network_choice}

% For the choice of the network we have two different cases: DEQ-S, $R_{\theta}$ is itself a scalar neural network or DEQ-RED, $R_{\theta}(x)=\frac{1}{2}\|x-N_{\theta}\|_2^2$ (section \ref{param of reg funz}).
For the DEQ-S regularization we considered a Convolutional neural network with residual connections as in \cite{ICNNpaper}. 

For DEQ-RED the choice was the DnCNN architecture introduced in \cite{DnCNN}: it is composed of a series of convolutional layers with ReLU activation functions which we replaced by Softplus activations defined by:
\begin{equation}
    \operatorname{Softplus}_{\beta}(x)=\frac{1}{\beta} \log (1+ \exp (\beta x)), \label{softplus}
\end{equation}
which makes the corresponding $R_{\theta}$  differentiable as required. We recall that the greater the parameter $\beta>0$ is, the more similar the Softplus function becomes to ReLU.

In the following table we report the number of trainable parameters.

\begin{longtable}{|p{5cm}|p{4cm}|p{5cm}|}
\hline
\textbf{Network} & \textbf{CNN} (for DEQ-S) & \textbf{DnCNN} (for DEQ-RED)\\
\hline
\endhead
\hline
\endfoot
\hline
\endlastfoot

\textbf{Activation function} &  
$\operatorname{Softplus}_{\beta}, \quad \beta=100$& $\operatorname{Softplus}_{\beta}, \quad \beta=100$ \\
\hline

\textbf{Trainable parameters $\theta$} &$\approx 5 \times 10^6 $ & $\approx 3\times 10^5 $  \\
\hline

\end{longtable}

\subsection{Forward pass}
We discuss here the details of the implementation regarding the forward pass. Recalling the analytic formula of \eqref{f_theta of DEMD}, 
we chose here $a=1$, which is reasonable for natural images  normalized in the range $[0,1]$. 
%As discussed in Section \ref{Back pass}, to have differentiability, we replaced the hard projection with $\Pi^{\epsilon}_{(-\infty,1)}$, with $\epsilon=10^{-5}$. 
In the spirit of DEQs, we stopped the algorithm when accurate (up to numerical precision) convergence is achieved, that is when the following stopping condition is met:
\begin{equation}
    \frac{\|x^{k+1}-x^k\|}{\|x^{k+1}\|} < \epsilon=2.5\times 10^{-5}. \label{convergence criterion}
\end{equation}
This value of this tolerance was selected based on empirical observation during numerical experimentation, representing an optimal trade-off: sufficiently small to ensure accuracy, yet large enough to maintain computational efficiency by preventing an excessive number of iterations.

\begin{figure}[h!] % Opzioni per il posizionamento della figura
    \centering
    \begin{subfigure}[t]{0.48\textwidth} % Larghezza della seconda subfigura
        \includegraphics[width=\textwidth]{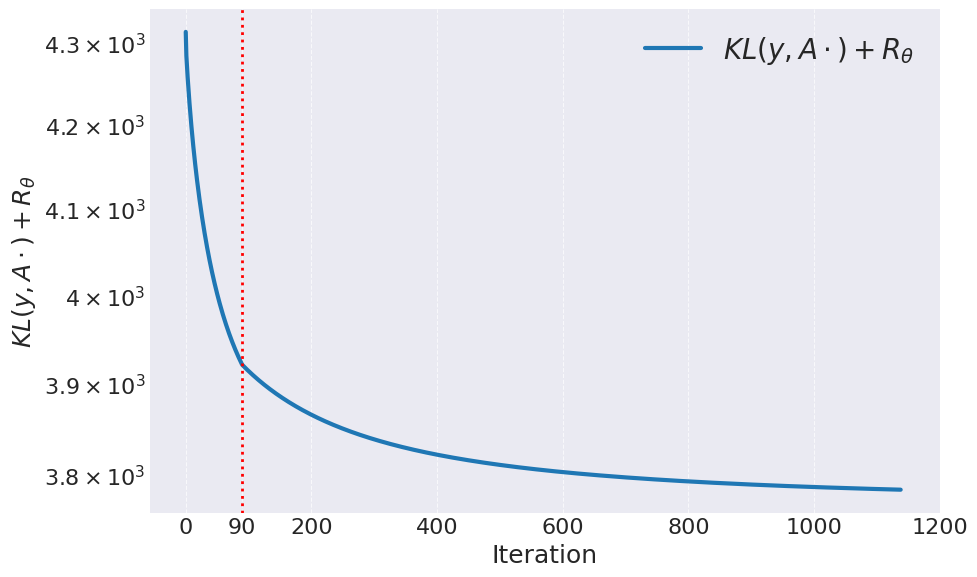} % Sostituisci con il nome del tuo file
        \caption{Functional values evaluated on the MD \\ sequence.}
        \label{fun val}
    \end{subfigure}
    \hfill % Spazio orizzontale tra le subfigure (riempie lo spazio rimanente)
    \begin{subfigure}[t]{0.48\textwidth} % Larghezza della prima subfigura (circa metà della larghezza del testo)
        \includegraphics[width=\textwidth]{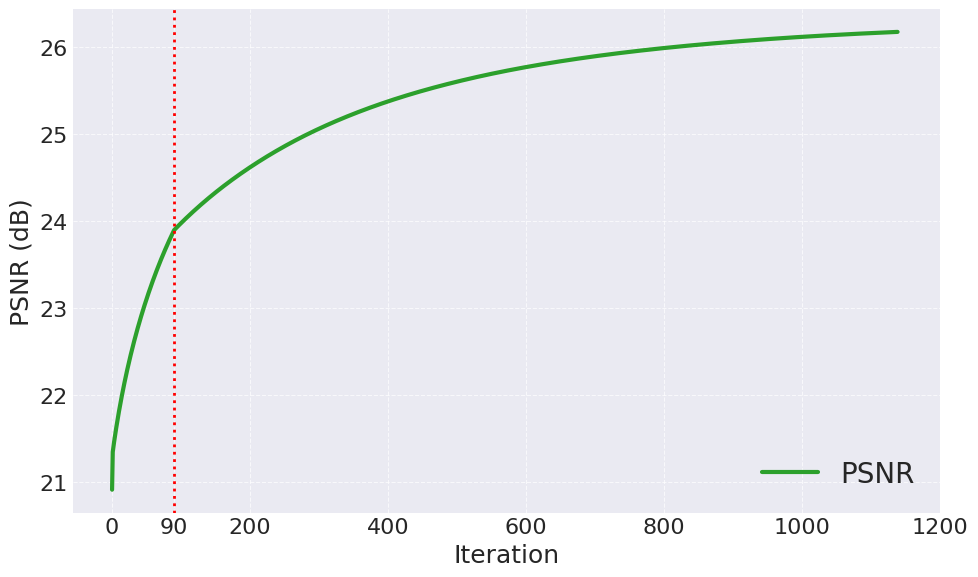} % Sostituisci con il nome del tuo file
        \caption{PSNR improvement along MD iterations.}
        \label{psnr w.r.t. gt}
    \end{subfigure}

    \caption{DEQ-MD forward pass at test time: functional values decrease and PSNR improvement along iterations.}
\end{figure}

We now report and comment some numerical results obtained at test time by running the forward pass of the proposed DEQ-RED approach on  a Gaussian deblurring task with $\alpha=40$ (see Section \ref{data generation} and Figure \ref{fig:kernels}). Figure \ref{fun val} illustrates the interpretability of the approach: running the DEQ-MD forward pass corresponds to minimize $\text{KL}(y,A\cdot)+R_{\theta}(\cdot)$, as detailed in Section \ref{eq:DEQ-MD_convergence}.
%We have theoretical guarantees about this fact. 
On the other hand, Figure \ref{psnr w.r.t. gt} shows that PSNR metric is increasing along iterates, with the final iterate maximizing it, which shows the adequacy of the minimization functional to the task considered in terms of reconstruction quality. %In fact it suffices to run the algorithm until convergence, in the sense of \eqref{convergence criterion}. We highlight that there are no theoretical guarantees about the maximization of PSNR along the iterates, but is \
We stress that this behavior was consistently observed on all the images considered for testing.

%While we currently lack a formal theoretical guarantee for this observed behavior, it is consistently encountered in practical applications.

\begin{figure}[h!] % Opzioni per il posizionamento della figura
    \centering
    \begin{subfigure}[t]{0.49\textwidth} % Larghezza della seconda subfigura
        \includegraphics[width=\textwidth]{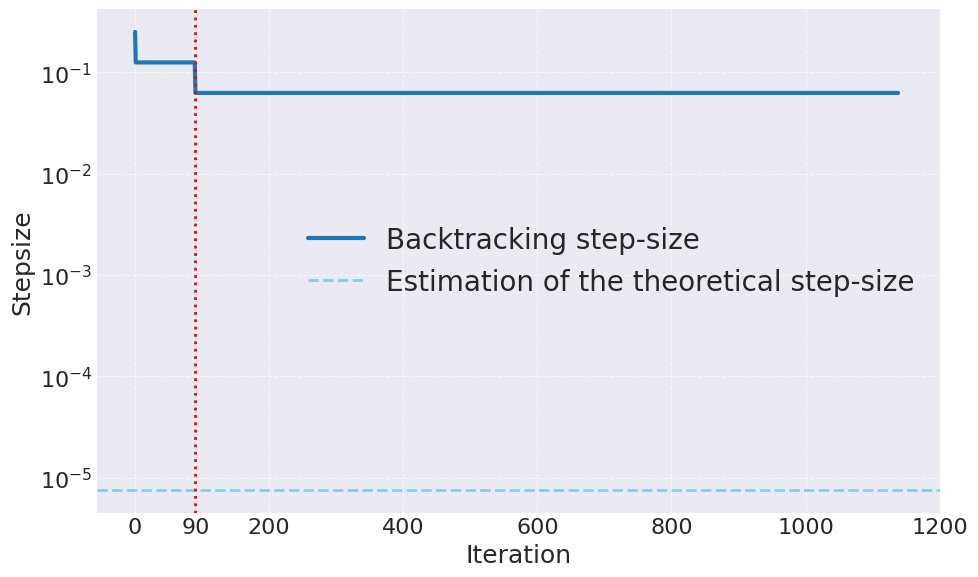} % Sostituisci con il nome del tuo file
        \caption{Step-size evolution given by Backtracking.}
        \label{backtracking step-sizes} % La label della subfigura va dopo la sua caption
    \end{subfigure}
    \hfill % Spazio orizzontale tra le subfigure (riempie lo spazio rimanente)
    \begin{subfigure}[t]{0.49\textwidth} % Larghezza della prima subfigura (circa metà della larghezza del testo)
        \includegraphics[width=\textwidth]{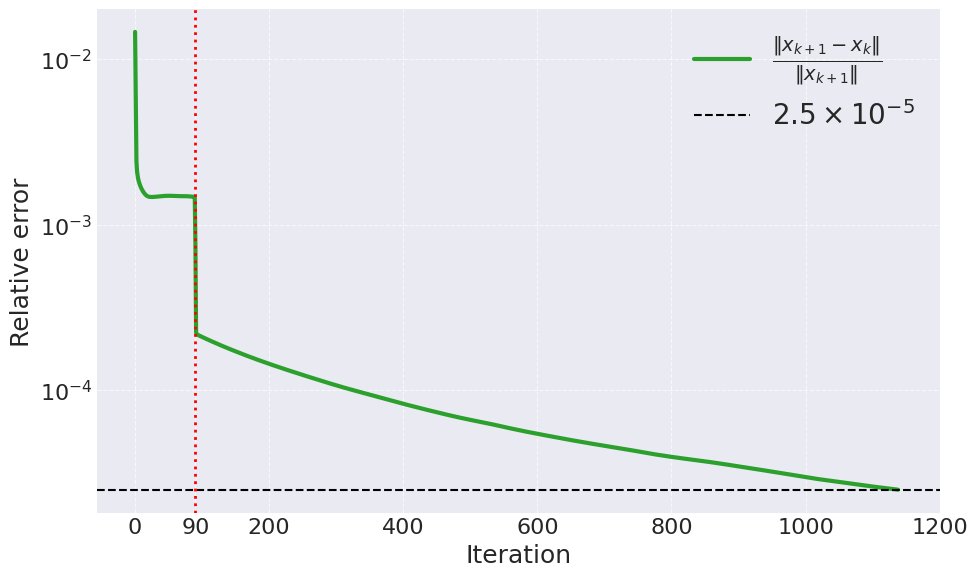} % Sostituisci con il nome del tuo file
        \caption{Relative error along the MD iterations.}
        \label{relative error} % La label della subfigura va dopo rla sua caption
    \end{subfigure}
    \caption{DEQ-MD forward pass at test time: step-size variations due to backtracking and  relative error behavior. The estimation of the theoretical step-size is computed as $\overline \tau =\frac{1}{L}$, where $L=\|y\|_1$ is the constant satisfying \eqref{eq nolip} for the pair $(\text{KL}(y,A\cdot),h)$. This value provides an upper bound to the theoretical step-size $\tau<\frac{1}{L+L'}<\frac{1}{L}= \overline{\tau}$, where $L'$ is the constant satisfying \eqref{eq nolip} for the pair $(R_{\theta},h)$.}
\end{figure}

Figure \ref{backtracking step-sizes} highlights the effectiveness of the backtracking strategy used to estimate admissible step-sizes guaranteeing convergence. Note that the such step-sizes are way bigger than theoretical one, which allows for faster algorithmic convergence. Figure \ref{relative error} demonstrates the decrease in relative error across iterations. We observe a correlation between changes in the step-size and variations in the relative error, as indicated by the vertical dashed red line.
%\textcolor{red}{We conclude these section observing that the consideration of these sections hold in general boht for DEQ-RED and DEQ-S and for different kernel and level of noise.}

\subsubsection{On the choice of the initialization}
Given the non-convex optimization regime where we are considering, the operator $f_{\theta}(\cdot;y)$ \eqref{f_theta of DEMD} may have multiple fixed points associated to the multiple critical points of the functional  we are minimizing (Proposition \ref{prop: convergence}). The choice of the initialization may thus affect the point the algorithm converges to, and, as such, it can be considered as an hyper-parameter affecting the quality of the results. Furthermore, in practice, the choice of initialization inputs algorithm efficiency as it impacts the number of iterations needed to get convergence. We thus analyzed four different choices of initializations $x^0$ in the forward pass:
\begin{itemize}
    \item \textbf{Adjoint-based}:  $x^0=A^*y$.
    \item \textbf{Random}: the components of $x^0$ are drawn independently from a  uniform distribution, that is $x^0_i \sim \mathcal{U}((0,1))$ for all $i=1,...,n$. 
    \item \textbf{Total Variation reconstruction}: $x^0$ is chosen as the numerical solution to: \begin{equation}
        x^0 \in
    \underset{x \geq 0}{\operatorname{argmin}}~\text{KL}(y,Ax)+ \lambda \textrm{TV}_{\epsilon}(x), \label{TV minimization}
    \end{equation}
    with $\lambda>0$, where $\textrm{TV}_\epsilon$ denotes a smoothed version of the Total Variation (TV) regularization functional.
    TV regularization is a widely used technique, particularly effective for noise reduction and image reconstruction tasks \cite{TV}. It operates by promoting solutions that have sparse gradients, which translates to preserving sharp edges while smoothing out noise in homogeneous regions. This approach is particularly beneficial for images containing piecewise constant or piecewise smooth structures. The TV functional in its natural expression is convex but not differentiable in $x=0$. For this reason, we smoothed it by  adding a small perturbation depending on $\epsilon>0$, that is by considering:
    \begin{equation}
        \textrm{TV}_{\epsilon}(x) = \sum_{i=1}^{N-1} \sum_{j=1}^{N-1} \sqrt{(x_{i,j+1} - x_{i,j})^2 + (x_{i+1,j} - x_{i,j})^2 + \epsilon},
        \label{TV}
    \end{equation}
   where here $x\in \mathbb{R}^{N \times N}$ with $n=N\times N$. The regularizer TV$_{\epsilon}$ is differentiable on $\mathbb{R}^{N \times N}$. We solved problem  \eqref{TV minimization} by running Mirror Descent iterations
    \begin{equation}  \label{eq:MD_TV}
       x^{k+1}=\nabla h^*(\nabla h(x^k)-\tau (\nabla \text{KL}(y,Ax^k))+\lambda\nabla \text{TV}_{\epsilon}(x^k))
    \end{equation}
     till convergence
   for a suitable choice of the step-size and for an optimal $\lambda$ chosen so as to maximize the PSNR w.r.t.~the ground truth.
    \item \textbf{Richardson-Lucy} \cite{Richardson72,Lucy1974}: we consider
    $x^0=x^K_{RL}$,  that is the $K$-th iteration of Richardson Lucy (RL) algorithm with $K>0$. The latter is an EM algorithm designed to minimize the Kullback-Leibler functional \eqref{KL equation} and particularly used in the context of image deblurring, particularly effective for data corrupted by Poisson noise. The multiplicative iteration of the algorithm reads:
\begin{equation}
x^{k+1}=\frac{x^k}{A^* \textbf{1}} \left( A^{*}\frac{y}{Ax^k} \right), \label{RL}
\end{equation}
where $\textbf{1}$ is a vector of ones, i.e., $\textbf{1}=(1,...,1) \in \mathbb{R}^m$, and division is intended element-wise.
%A key characteristic of this multiplicative scheme is that if the initial guess $x^0$ has strictly positive components ($x^0 \in \mathbb{R}^n_{++}$), then all subsequent iterates $\{x^k\}_k$ will also remain strictly positive.
%However, the RL algorithm in its basic form \eqref{RL} lacks an explicit regularization term. This can lead to overfitting noise if the iterations are continued until convergence, resulting in suboptimal reconstructions. Therefore, a robust stopping criterion is crucial for achieving high-quality deblurred images. 
Note that in such iterative regularization setting, the choice of an optimal early stopping $K$ corresponds to the choice of an optimal regularization parameter similarly as for $\lambda$ above.
\end{itemize}

\begin{figure}[h]
    \centering
    \includegraphics[width=0.75\linewidth]{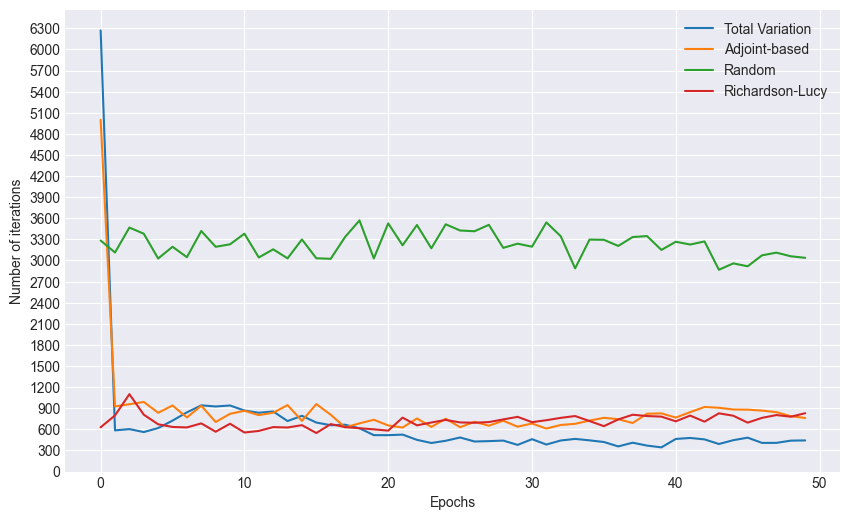}
    \caption{Number of iterations performed in the forward pass to reach the stopping criterion threshold for four different initializations, averaged on the test set, varying with epochs.}
    \label{number iter for different init}
\end{figure}

We compared the choice of the four initializations above by training four different DEQ-RED models following the modalities specified in Section  \ref{loss function} and using a subset of the image dataset described in Section \ref{dataset}, with training set size of 50, a validation set size of 20 and a test set size of 20. We considered an image deblurring problem defined in terms of a Gaussian kernel (see Figure \ref{fig:kernels}) and a level of Poisson noise $\alpha=100$.  
As shown in Figure \ref{number iter for different init}, in terms of iterations required to achieve convergence for the forward pass, Adjoint-based, TV and RL are comparable, while the Random initialization always requires many more iterations.

In Table \ref{tab:confronto_metodi}, we report the computational times  needed to get all such initializations and the PSNR values  obtained after evaluating each of the four DEQ-MD trained models on the test set. We notice that PSNR values do not change significantly except for $x^0$ being a random initialization , although we observe that computational time required to get $x^0$ in the case of TV is way more expensive than the others.  As a tradeoff, we thus decided to use for the following tests the Adjoint-based initialization.

\begin{table}[h]
\centering
\begin{tabular}{|c|c|c|c|c|}
\hline
 & Adjoint-based & $\textrm{TV}_\epsilon$ & RL & Random \\
\hline
PSNR (dB) &  25.52 & \textbf{25.57 } & 25.54  &  24.89  \\
CPU time (s) & 1e-3& 3.05& 5e-2& 1e-4\\
\hline
\end{tabular}
\caption{Average PSNR  and CPU times of the solution to the  four DEQ-RED models applied on the test set and initialized by choosing the four different initializations considered.}
\label{tab:confronto_metodi}
\end{table}

We conclude this section by highlighting that we also explored a warm restart initialization technique during the training phase by taking as initialization of the forward pass at epoch $i+1$ the reconstructions obtained at epoch $i$. This procedure, which theoretically should accelerate training, did not yield significant effects in terms of either computational time nor reconstruction quality, hence we did not consider it in our numerical results. 
%Note, moreover, that another drawback of this approach is that the same initializations used during training cannot be applied during inference.

% \subsection{Training details} \label{training details}

% We discuss here further training aspects, concerned with the architecture of the neural networks employed, the image dataset used and the loss function considered.

\subsection{Training results} \label{training results}

\begin{figure}[h]
    \centering
    \includegraphics[width=0.65\linewidth]{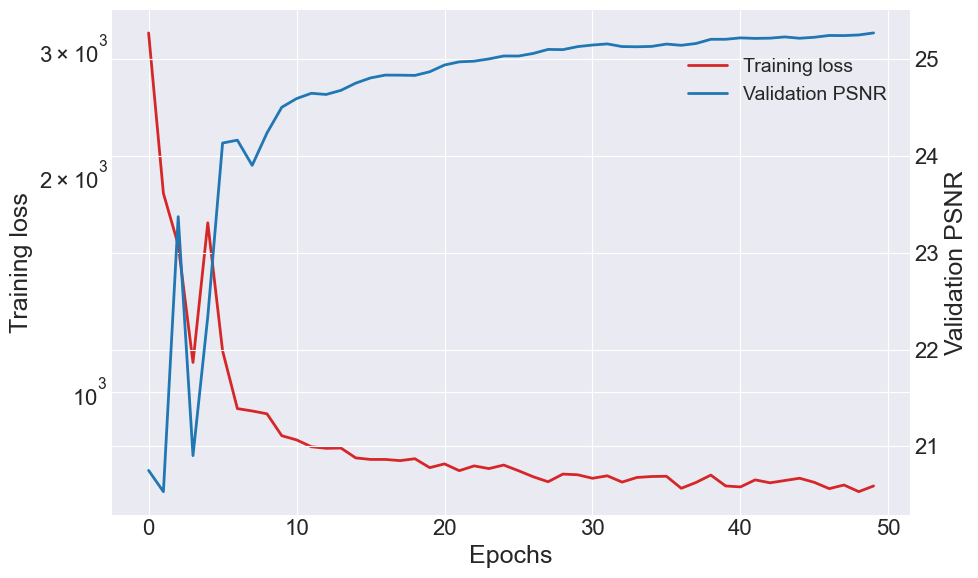}
    \caption{Training curves along epochs for DEQ-RED with $\alpha=60$ and Gaussian kernel. Training loss evaluated on training set (red), PSNR evaluated on validation set (blue).}
    \label{fig: training metrics}
\end{figure}

Training was performed by following the details described in Section \ref{loss function} and using the image dataset in Section \ref{dataset}. We trained a different DEQ-RED and DEQ-S for each level of noise $\alpha$ and tasks described in Section \ref{data generation}.
 As an illustration, we report in Figure \ref{fig: training metrics} the training curves obtained using the DEQ-RED model trained on the Gaussian kernel with $\alpha=60$. We can notice that after few epochs the metrics stabilize with no significative improvements after 50 epochs at the price of a significant computational overhead. Although theoretical guarantees for JFB are lacking, the empirical results presented in Figure \ref{fig: training metrics} show a positive trend in the validation loss and PSNR. Analogous considerations hold for both DEQ-RED and DEQ-S w.r.t.~all level of noise and convolution kernels considered.

To conclude this section, we present a study on the the impact of the choice of loss used for training the DEQ-MD models presented. For this purpose, we compare in the following the MSE loss \eqref{mixed loss} with the choices:
\begin{itemize}
    \item[-] $\ell(x^*,x^{\infty})=\text{KL}(x^*,x^\infty)$;
    \item[-]  $\ell(x^*,x^{\infty})=\|x^*-x^{\infty}\|_1$.
\end{itemize}
This comparison has been conducted training three DEQ-RED models on the same Gaussian deblurring task using the highest level of noise presented, i.e. $\alpha=40$ on a test dataset. From Table \ref{tab:comparison loss funz} we can see that all the three metrics provides comparable results with MSE being slightly better. Since, in addition, training $\text{KL}$ turned out to be less stable during the training iterations, the MSE  loss was employed throughout the experiments.

\begin{table}[h] 
\centering

\caption{Comparison between different loss functions used in training. Average results computed  on the test set with $\alpha=40$ and Gaussian blur.}
\label{tab:comparison loss funz}
\begin{tabular}{lccc}
\toprule
\textbf{Loss Function} & \textbf{PSNR $\uparrow$} & \textbf{SSIM $\uparrow$} & \textbf{LPIPS $\downarrow$}  \\
\midrule
MSE & \textbf{25.9250} & \textbf{0.7375} & \textbf{0.2921} \\
KL& \underline{25.8152} & 0.7313 & 0.3198 \\
L1 & 25.7876 & \underline{0.7341} & \underline{0.3144} \\
\bottomrule
\end{tabular}%

\end{table}

% \begin{table}[H]
%     \centering
%     \caption{Average Computational Times}
%     \label{tab:computational_times}
%     \begin{tabular}{|c|c|c|}
%         \hline
%         \textbf{Method} & \textbf{Full Training} & \textbf{Single epoch} \\
%         \hline
%         DEQ-S           & $\approx 105$ \text{ hours} & $\approx{130} $ \text{ mins} \\
%         \hline
%         DEQ-RED         & $\approx 50$ \text{ hours} & $\approx 60$ \text{ mins} \\
%         \hline
%     \end{tabular}
% \end{table}

\newpage
\section{Numerical experiments} \label{numerical experiments}

In this section, we present several numerical experiments validating the  DEQ-MD framework. We tested our method   under the settings described in Section \ref{data generation}. The image reconstructions obtained are compared with the following model-based and data-driven approaches tailored to Poisson noise removal:
%\textcolor{red}{At the end of this Section we study the generalization property of our approach. For example we tested DEQ-MD on a Super Resolution task, on which the model was not trained.}

% \begin{figure}[h!]
%     \centering
%     \begin{subfigure}[b]{0.28\textwidth}
%         \includegraphics[width=\linewidth]{images/motion_kernel_1.png}
%         \caption{Motion blur kernel}
%         \label{motion blur}
%     \end{subfigure}
%     \hfill
%     \begin{subfigure}[b]{0.28\textwidth}
   
%      \includegraphics[width=\linewidth]{images/gaussian_kernel.png}
%         \caption{Gaussian Kernel} \label{gaussian kernel}
%     \end{subfigure}
%     \hfill
%     \begin{subfigure}[b]{0.28\textwidth}
%         \includegraphics[width=\linewidth]{images/uniform_kernel.png}
%         \caption{Uniform kernel}
%         \label{uniform kernel}
%     \end{subfigure}
    
%     \caption{The 3 blur kernels used in the experiments. (a) is a real-world camera shake kernel, see \cite{Levin}. (b) is a \( 11 \times 11 \) Gaussian kernel with standard deviations $\sigma$ =1.2. (c) is a \( 9 \times 9 \) uniform kernel. }
%     \label{fig:kernels}
% \end{figure}

% We now present in detail the competing approaches considered:
\begin{itemize}
    \item[--] \textbf{Model-based}:
    \begin{itemize}
   \item[$\bullet$] \textbf{Richardson Lucy (RL)} \cite{Richardson72,Lucy1974} algorithm \eqref{RL}: we ran the algorithm for a fixed maximum number of $100$ iterations and chose as optimal early stopping the one whose corresponding iterate maximizes the PSNR w.r.t. ground truth.
    \item[$\bullet$] \textbf{\text{KL}+smoothed Total Variation (KL+TV)} \eqref{TV minimization}: we ran the instance \eqref{eq:MD_TV} of MD  for different values of $\lambda$, choosing as optimal value the one for which the corresponding solution maximizes the PSNR w.r.t.~ground truth.
    \end{itemize}
    \item[-]  \textbf{Data-driven}: 
    PnP-based methods tailored for Poisson Inverse Problems in the form \eqref{PIP} were considered in \cite{Hurault} and addressed using MD as an underlying optimization algorithm. In this framework, the denoiser considered takes the form of a generalization of a gradient step denoiser, and it is trained on data corrupted by Inverse Gamma noise. By choosing $h$ as the Burg's entropy \eqref{burg's entropy}, the Bregman proximal operator associated turns out to be the MAP estimator of an Inverse Gamma image denoising problem. To be more precise, considering a potential $g_{\gamma}: \mathbb{R}^n \to \mathbb{R}$ defined by: 
    \begin{equation}
        g_\gamma(x)=\frac{1}{2}\|x-N_{\gamma}(x)\|_2^2,
    \end{equation} 
   the corresponding Bregman denoiser takes the form:
    \begin{equation}
        \mathcal{B}_{\gamma}(x)=x- (\nabla^2h(x))^{-1}\nabla g_{\gamma}(x) . \label{bregman score denoiser}
    \end{equation}
    Such denoiser is then trained on different levels $\gamma$ of Inverse Gamma noise. %\textcolor{red}{One might argue that the training should have been performed on Poisson noise. However, according to the PnP philosophy, the denoising task is not dependent on the specific problem we aim to solve, but rather on the optimization algorithm employed, which makes Inverse Gamma noise a natural choice for this scenario.}
    
    %Taking $h=\frac{1}{2}\|x\|_2^2$, then \eqref{bregman score denoiser} is the well-known Gradient step denoiser \cite{Hurault1}.
   Within this framework, two different algorithms can be considered:
    \begin{itemize}
        \item[$\bullet$] \textbf{B-RED}:
    \begin{equation}
        x^{k+1}=\Pi_{[0,1]^n}(\nabla h^*(\nabla h(x^k))-\tau(\lambda\nabla \text{KL}(y,Ax^k)+\nabla g_{\gamma}(x^k))).\label{B-RED eq}
    \end{equation}
  Here, the step-size $\tau>0$ is chosen by employing the backtracking strategy in Algorithm \ref{BT}, $\lambda>0$ is the regularization parameter and $\gamma$ is the level of Inverse Gamma noise to be tuned.
    \item[$\bullet$] \textbf{B-PnP}:
    \begin{equation}
        x^{k+1}=\mathcal{B}_{\gamma}(\nabla h^*(\nabla h(x^k)-\tau\nabla \text{KL}(y,Ax^k)),
    \end{equation}
  where the step-size is chosen such that $\tau L<1$, with $L$ given by Lemma \ref{NoLip constant of KL}.
    \end{itemize}
    %The numerical implementation of these two algorithms is based on \url{https://github.com/samuro95}.  
\end{itemize}

%We show now our results on the reference task.

\begin{figure}[h!]
    \centering
    \begin{subfigure}[b]{0.24\textwidth} % Modificato la larghezza per 4 immagini
        \begin{overpic}[width=\linewidth]{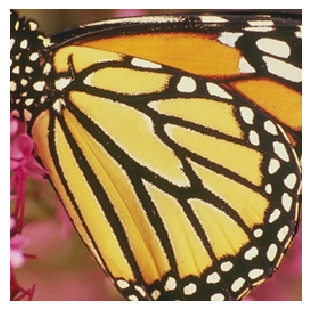}

        \put(2,2){\includegraphics[width=1cm]{images/gaussian_kernel.png}}
        \put(35.5,29){%
            \setlength{\fboxrule}{2.5pt}%
            \color{blue}\framebox(13,13){}%
        }
        \put(57,3.5){%
            \setlength{\fboxrule}{2.5pt}%
            \color{blue}\framebox(40,40){}%
        }
        %sx bass dx alto
        \put(57,3.5){%
            \includegraphics[width=1.5cm,
            trim= 77 67 117 127,  clip]{images/butterfly_gauss_alpha_40/GT.png}%
        }
    \end{overpic}
        \caption{Clean \\\mbox{\footnotesize(PSNR$\uparrow$, SSIM$\uparrow$, LPIPS$\downarrow$)}}
    \end{subfigure}
    \hfill
    \begin{subfigure}[b]{0.24\textwidth} % Modificato la larghezza per 4 immagini
                \begin{overpic}[width=\linewidth]{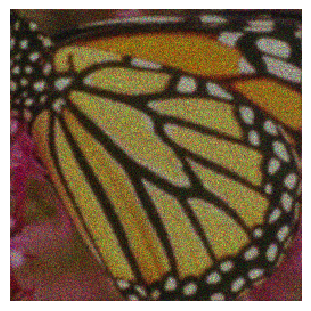}
        \put(35.5,29){%
            \setlength{\fboxrule}{2.5pt}%
            \color{blue}\framebox(13,13){}%
        }
        \put(57,3.5){%
            \setlength{\fboxrule}{2.5pt}%
            \color{blue}\framebox(40,40){}%
        }
        %sx bass dx alto
        \put(57,3.5){%
            \includegraphics[width=1.5cm,
            trim=77 67 117 127, clip]{images/butterfly_gauss_alpha_40/mes.png}%
        }
    \end{overpic}
        \caption{Observed \\
        (14.047,
        0.409,
        0.452) }
    \end{subfigure}
    \hfill
    \begin{subfigure}[b]{0.24\textwidth} % Modificato la larghezza per 4 immagini
                               \begin{overpic}[width=\linewidth]{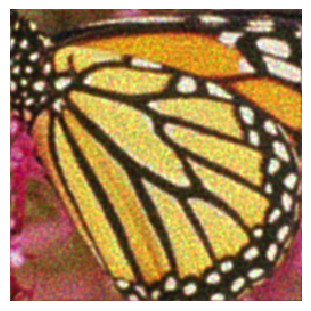}
        \put(35.5,29){%
            \setlength{\fboxrule}{2.5pt}%
            \color{blue}\framebox(13,13){}%
        }
        \put(57,3.5){%
            \setlength{\fboxrule}{2.5pt}%
            \color{blue}\framebox(40,40){}%
        }
        %sx bass dx alto
        \put(57,3.5){%
            \includegraphics[width=1.5cm,
            trim=77 67 117 127, clip]{images/butterfly_gauss_alpha_40/RL.png}%
        }
    \end{overpic}
        \caption{RL \\
        (22.149, 0.638,
        0.355)}
    \end{subfigure}
    \hfill
    \begin{subfigure}[b]{0.24\textwidth} % Modificato la larghezza per 4 immagini
                                 \begin{overpic}[width=\linewidth]{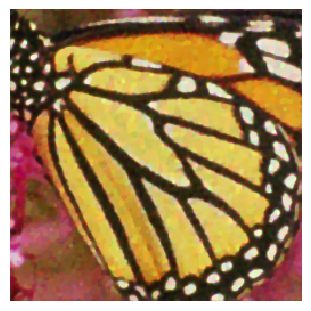}
        \put(35.5,29){%
            \setlength{\fboxrule}{2.5pt}%
            \color{blue}\framebox(13,13){}%
        }
        \put(57,3.5){%
            \setlength{\fboxrule}{2.5pt}%
            \color{blue}\framebox(40,40){}%
        }
        %sx bass dx alto
        \put(57,3.5){%
            \includegraphics[width=1.5cm,
            trim=77 67 117 127, clip]{images/butterfly_gauss_alpha_40/TV.png}%
        }
    \end{overpic}
        \caption{KL+TV$_\epsilon$ \\
        (23.583, 0.761, 0.174)}
    \end{subfigure}
    \vskip\baselineskip
    \begin{subfigure}[b]{0.24\textwidth} % Modificato la larghezza per 4 immagini
                                 \begin{overpic}[width=\linewidth]{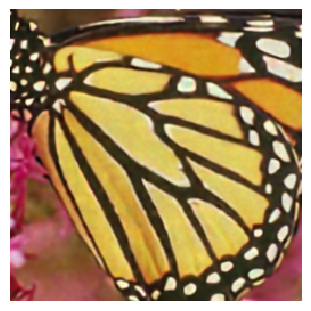}
        \put(35.5,29){%
            \setlength{\fboxrule}{2.5pt}%
            \color{blue}\framebox(13,13){}%
        }
        \put(57,3.5){%
            \setlength{\fboxrule}{2.5pt}%
            \color{blue}\framebox(40,40){}%
        }
        %sx bass dx alto
        \put(57,3.5){%
            \includegraphics[width=1.5cm,
            trim=77 67 117 127, clip]{images/butterfly_gauss_alpha_40/B-RED.png}%
        }
    \end{overpic}
        \caption{B-RED \\
        (\underline{25.894},0.855,0.161)
        }
    \end{subfigure}
    \hfill
    \begin{subfigure}[b]{0.24\textwidth} % Modificato la larghezza per 4 immagini
                                 \begin{overpic}[width=\linewidth]{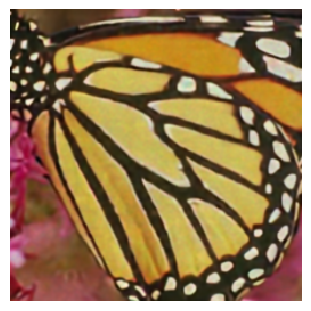}
        \put(35.5,29){%
            \setlength{\fboxrule}{2.5pt}%
            \color{blue}\framebox(13,13){}%
        }
        \put(57,3.5){%
            \setlength{\fboxrule}{2.5pt}%
            \color{blue}\framebox(40,40){}%
        }
        %sx bass dx alto
        \put(57,3.5){%
            \includegraphics[width=1.5cm,
            trim=77 67 117 127, clip]{images/butterfly_gauss_alpha_40/BPNP.png}%
        }
    \end{overpic}
        \caption{B-PnP \\
        (25.740, \underline{0.851},0.183)}
    \end{subfigure}
    \hfill
    \begin{subfigure}[b]{0.24\textwidth} % Modificato la larghezza per 4 immagini
                                 \begin{overpic}[width=\linewidth]{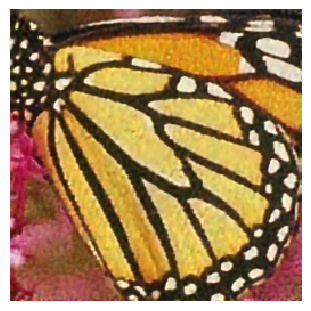}
        \put(35.5,29){%
            \setlength{\fboxrule}{2.5pt}%
            \color{blue}\framebox(13,13){}%
        }
        \put(57,3.5){%
            \setlength{\fboxrule}{2.5pt}%
            \color{blue}\framebox(40,40){}%
        }
        %sx bass dx alto
        \put(57,3.5){%
            \includegraphics[width=1.5cm,
            trim=77 67 117 127, clip]{images/butterfly_gauss_alpha_40/ICNN.png}%
        }
    \end{overpic}
        \caption{DEQ-S \\
    (24.566, 0.754, \underline{0.157})}
    \end{subfigure}
    \hfill
\begin{subfigure}[b]{0.24\textwidth}
    \centering
                             \begin{overpic}[width=\linewidth]{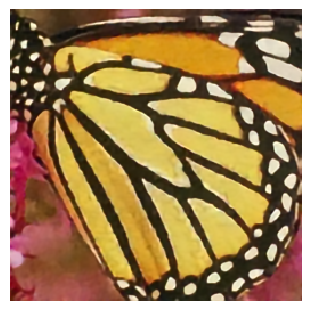}
        \put(35.5,29){%
            \setlength{\fboxrule}{2.5pt}%
            \color{blue}\framebox(13,13){}%
        }
        \put(57,3.5){%
            \setlength{\fboxrule}{2.5pt}%
            \color{blue}\framebox(40,40){}%
        }
        %sx bass dx alto
        \put(57,3.5){%
            \includegraphics[width=1.5cm,
            trim=77 67 117 127, clip]{images/butterfly_gauss_alpha_40/DNcnn.png}%
        }
    \end{overpic}
    \caption{DEQ-RED \\ 
    \mbox{(\textbf{26.417}, \textbf{0.864}, \textbf{0.083})}}
\end{subfigure}
    \caption{Reconstruction comparison in terms of PSNR, SSIM and LPIPS for image deblurring with Gaussian blur and Poisson noise with $\alpha=40$.}
    \label{gaussian deblurring}
\end{figure}

In Figure \ref{gaussian deblurring}, we evaluate the results obtained by the different methods on a Gaussian deblurring task (Fig.~\ref{gaussian kernel}) at test time with noise level $\alpha=40$.
Notice that the absence of explicit regularization for Richardson-Lucy (RL) yields poor reconstruction results.
In contrast, the reconstruction computed using the KL+TV$_\epsilon$ model provides good reconstruction quality, though it remains inferior to both Bregman PnP and DEQ-RED. We observed that DEQ-S struggles to effectively remove image noise. This aligns with findings in \cite{Salimans}, where it is observed that directly parameterizing regularization functional through a neural network can lead to suboptimal results.
Regarding B-RED and B-PnP, we observe that their performance is remarkably similar, with B-RED performing slightly better in terms of the three metrics. Finally, DEQ-RED provides the best reconstruction, achieving a significantly higher PSNR value.

\begin{figure}[h]
    \centering
\begin{subfigure}[b]{0.24\textwidth}
    \begin{overpic}[width=\linewidth]{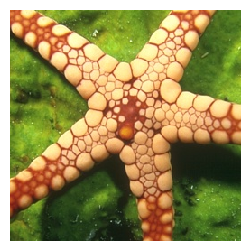}

        \put(2,2){\includegraphics[width=1cm]{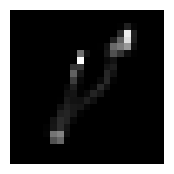}}
        \put(35,50){%
            \setlength{\fboxrule}{2.5pt}%
            \color{blue}\framebox(20,20){}%
        }
        \put(55.5,4.3){%
            \setlength{\fboxrule}{2.5pt}%
            \color{blue}\framebox(40,40){}%
        }
        \put(55.5,4.3){%
            \includegraphics[width=1.5cm, trim=90 130 110 70, clip]{images/starfish_motion_7_alpha_60/gt.png}%
        }
    \end{overpic}
    \caption{Clean \\ \mbox{\footnotesize(PSNR$\uparrow$, SSIM$\uparrow$,  LPIPS$\downarrow$)}}
\end{subfigure}
    \hfill
    \begin{subfigure}[b]{0.24\textwidth} % Modificato la larghezza per 4 immagini
            \begin{overpic}[width=\linewidth]{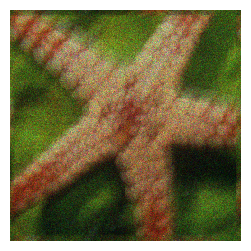}
        \put(35,50){%
            \setlength{\fboxrule}{2.5pt}%
            \color{blue}\framebox(20,20){}%
        }
        \put(55.5,4.3){%
            \setlength{\fboxrule}{2.5pt}%
            \color{blue}\framebox(40,40){}%
        }
        \put(55.5,4.3){%
            \includegraphics[width=1.5cm, trim=90 130 110 70, clip]{images/starfish_motion_7_alpha_60/mes.png}%
        }
    \end{overpic}
        \caption{Observed \\ 
        (13.306,
        0.1664, 
        0.573)}
    \end{subfigure}
    \hfill
    \begin{subfigure}[b]{0.24\textwidth} % Modificato la larghezza per 4 immagini
                    \begin{overpic}[width=\linewidth]{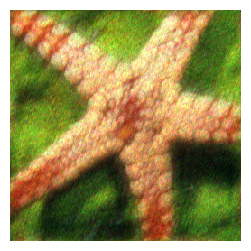}
        \put(35,50){%
            \setlength{\fboxrule}{2.5pt}%
            \color{blue}\framebox(20,20){}%
        }
        \put(55.5,4.3){%
            \setlength{\fboxrule}{2.5pt}%
            \color{blue}\framebox(40,40){}%
        }
        \put(55.5,4.3){%
            \includegraphics[width=1.5cm, trim=90 130 110 70, clip]{images/starfish_motion_7_alpha_60/rl.png}%
        }
    \end{overpic}
        \caption{RL \\ 
        (20.458, 
        0.474, 
        0.444)}
    \end{subfigure}
    \hfill
    \begin{subfigure}[b]{0.24\textwidth} % Modificato la larghezza per 4 immagini
                        \begin{overpic}[width=\linewidth]{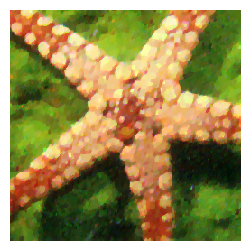}
        \put(35,50){%
            \setlength{\fboxrule}{2.5pt}%
            \color{blue}\framebox(20,20){}%
        }
        \put(55.5,4.3){%
            \setlength{\fboxrule}{2.5pt}%
            \color{blue}\framebox(40,40){}%
        }
        \put(55.5,4.3){%
            \includegraphics[width=1.5cm, trim=90 130 110 70, clip]{images/starfish_motion_7_alpha_60/tv.png}%
        }
    \end{overpic}
        \caption{KL+TV$_\epsilon$ \\ 
        (22.386,
        0.582,
        0.292)}
    \end{subfigure}
    \vskip\baselineskip
    \begin{subfigure}[b]{0.24\textwidth} % Modificato la larghezza per 4 immagini
                        \begin{overpic}[width=\linewidth]{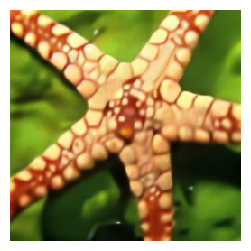}
        \put(35,50){%
            \setlength{\fboxrule}{2.5pt}%
            \color{blue}\framebox(20,20){}%
        }
        \put(55.5,4.3){%
            \setlength{\fboxrule}{2.5pt}%
            \color{blue}\framebox(40,40){}%
        }
        \put(55.5,4.3){%
            \includegraphics[width=1.5cm, trim=90 130 110 70, clip]{images/starfish_motion_7_alpha_60/rec_RED.png}%
        }
    \end{overpic}
        \caption{B-RED \\ (\underline{24.201}, 
        \underline{0.704}, 
        0.335)}
    \end{subfigure}
    \hfill
    \begin{subfigure}[b]{0.24\textwidth} % Modificato la larghezza per 4 immagini
                        \begin{overpic}[width=\linewidth]{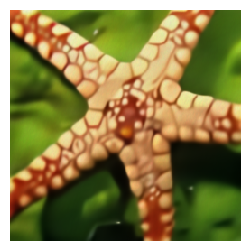}
        \put(35,50){%
            \setlength{\fboxrule}{2.5pt}%
            \color{blue}\framebox(20,20){}%
        }
        \put(55.5,4.3){%
            \setlength{\fboxrule}{2.5pt}%
            \color{blue}\framebox(40,40){}%
        }
        \put(55.5,4.3){%
            \includegraphics[width=1.5cm, trim=90 130 110 70, clip]{images/starfish_motion_7_alpha_60/rec_BPG.png}%
        }
    \end{overpic}
        \caption{B-PnP \\    (24.150, 
        0.703,
        0.336)}
    \end{subfigure}
    \hfill
    \begin{subfigure}[b]{0.24\textwidth} % Modificato la larghezza per 4 immagini
                        \begin{overpic}[width=\linewidth]{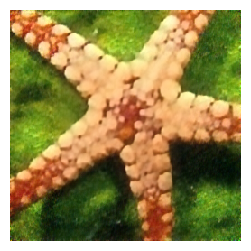}
        \put(35,50){%
            \setlength{\fboxrule}{2.5pt}%
            \color{blue}\framebox(20,20){}%
        }
        \put(55.5,4.3){%
            \setlength{\fboxrule}{2.5pt}%
            \color{blue}\framebox(40,40){}%
        }
        \put(55.5,4.3){%
            \includegraphics[width=1.5cm, trim=90 130 110 70, clip]{images/starfish_motion_7_alpha_60/rec_icnn.png}%
        }
    \end{overpic}
        \caption{DEQ-S \\    (23.216,
        0.637,
        \underline{0.267})}
    \end{subfigure}
    \hfill
    \begin{subfigure}[b]{0.2375\textwidth} % Modificato la larghezza per 4 immagini
                        \begin{overpic}[width=\linewidth]{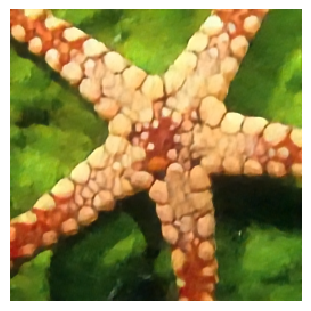}
        \put(35,50){%
            \setlength{\fboxrule}{2.5pt}%
            \color{blue}\framebox(20,20){}%
        }
        \put(57,3.5){%
            \setlength{\fboxrule}{2.5pt}%
            \color{blue}\framebox(40,40){}%
        }
        \put(57,3.5){%
            \includegraphics[width=1.5cm, trim=80 115 95 60, clip]{images/starfish_motion_7_alpha_60/rec_dncnn.png}%
        }
    \end{overpic}
        \caption{DEQ-RED \\ (\textbf{24.473},
        \textbf{0.717},
        \textbf{0.197})}
    \end{subfigure}
    \caption{Reconstruction comparison in terms of PSNR, SSIM and LPIPS for image deblurring task with motion blur and Poisson noise with $\alpha=60$.}
    \label{fig:motion blur deblurring}
\end{figure}

Figure \ref{fig:motion blur deblurring} illustrates a more challenging scenario involving motion blur. The detailed features of the 'Clean' image are hardly discernible in the 'Observed' image due to the action of the blur. The RL method proves particularly ineffective on this particular task, while the KL+TV$_\epsilon$ solution introduces excessive smoothing.
Both B-RED and B-PnP reconstructions are visually appealing; however, artifacts are noticeable, especially in the background where most details appear lost with a general tendency towards over-smoothing. 
%On the other hand, these models effectively remove Poisson noise.
Note that the DEQ-S reconstruction over-smooths details  leaving somewhat noisy the background.
DEQ-RED achieves the best overall reconstruction metrics, with texture-preserving smoothing enforced both in the background and in the textured areas.
%, especially considering the task's complexity. 
This example highlights the robustness of DEQ methods, and in particular of DEQ-RED, to artifacts. On the other hand, PnP approaches, while providing visually appealing denoised reconstructions tend to significantly over-smooth images.

\begin{table}[t!]
\centering
\caption{Performance comparison of different methods under levels of Poisson noise ranging from low to high $\alpha\in\left\{100,60,40\right\}$ averaged on the three convolution kernels considered. Evaluations are made on the test set.}
\label{table 1}
\footnotesize % Riduci il font
\resizebox{\textwidth}{!}{% <-- Inizia qui il ridimensionamento, larga quanto il testo
\begin{tabular}{@{}l *{9}{c} @{}}
\toprule
 & \multicolumn{3}{c}{$\alpha = 100$} & \multicolumn{3}{c}{$\alpha = 60$} & \multicolumn{3}{c}{$\alpha = 40$} \\
\cmidrule(lr){2-4} \cmidrule(lr){5-7} \cmidrule(lr){8-10}
Method & PSNR $\uparrow$ & SSIM $\uparrow$ & LPIPS $\downarrow$ & PSNR $\uparrow$ & SSIM $\uparrow$ & LPIPS $\downarrow$ & PSNR $\uparrow$ & SSIM $\uparrow$ & LPIPS $\downarrow$ \\
\midrule
RL & 22.7381 & 0.5379 & 0.5570 & 22.2968 & 0.5111 & 0.6017 & 21.9496 & 0.4901 & 0.6373 \\
KL+TV$_\epsilon$ & 23.9879 & 0.6304 & \underline{0.3888} & 23.5203 & 0.6045 & \underline{0.4254} & 23.1438 & 0.5844 & 0.4547 \\
B-RED & \underline{25.0710} & \underline{0.6939} & 0.4120 & \underline{24.6566} & \underline{0.6723} & 0.4384 & \textbf{24.4603} & 0.6530 & \underline{0.4533} \\
B-PnP & 24.8976 & 0.6847 & 0.4417 & 24.5818 & 0.6674 & 0.4401 & 24.2912 & \underline{0.6612} & 0.4534 \\
DEQ-S & 24.2080 & 0.6730 & 0.3935 & 23.7252 & 0.6221 & 0.4346 & 23.4423 & 0.5912 & 0.4591 \\
DEQ-RED & \textbf{25.4245} & \textbf{0.7091} & \textbf{0.3058} & \textbf{24.9087} & \textbf{0.6826} & \textbf{0.3354} & \underline{24.4586} & \textbf{0.6626} & \textbf{0.3571} \\
\bottomrule
\end{tabular}% <-- Il '%' qui è importante per evitare spazi extra
} % <-- Finisce qui il ridimensionamento
\end{table}

Table \ref{table 1} presents a comparative analysis of the methods on a test set composed by 46 images, aggregating results across three noise levels and three kernels.
For the low and medium noise scenario ($\alpha\in\left\{100,60\right\}$), DEQ-RED consistently achieved the highest values across all three metrics. 
At $\alpha=60$, B-RED edged out DEQ-RED slightly in terms of PSNR and SSIM, though their overall metric values remained very close.
In the high noise case ($\alpha=40$), B-RED  delivered a slightly superior PSNR. Conversely, DEQ-RED consistently proved to be the best performer for LPIPS and SSIM.
%In summary, DEQ-RED \textcolor{red}{provides a superior performance for the two lower level of noise, while in the highest noise condition the results are close to PnP} 

In Table \ref{tab:reconstruction_parameters} we present the number of hyperparameters required to obtain the reconstructions using all methods considered. RL and KL+TV$_\epsilon$ need each one only one parameter to be tuned per image. In contrast, B-RED and B-PnP require the tuning of two pairs of parameters: $(\lambda_{\text{in}},\gamma_{\text{in}})$ and $(\lambda_{\text{algo}},\gamma_{\text{algo}})$. This requirement stems from the methodology proposed in \cite{Hurault}, where an initial fixed number of iterations is performed using $(\lambda_{\text{in}},\gamma_{\text{in}})$ as regularization parameter/denoiser hyper-parameter to achieve a robust initialization. Subsequently, the primary algorithm is run using the second set of parameters $(\lambda_{\text{algo}},\gamma_{\text{algo}})$. Note however, that despite the theoretical need for parameter tuning for each reconstruction, in \cite{Hurault}  optimal parameter values depending solely on the noise level $\alpha$ are provided.  Contrarily to the previous methods, the proposed algorithms, DEQ-S, DEQ-RED are instead parameter-free, = having been trained for a specific forward operator and noise levels. Note that, however, both DEQ-MD approaches showed good generalization properties w.r.t.~different task and noise level, which will be discussed in Section \ref{gen: diff level of noise}.

\begin{table}[h]
    \centering
    \caption{Number of hyperparameters needed at inference time for each reconstruction method considered.}
    \label{tab:reconstruction_parameters}
\begin{tabular}{|l|c|p{10.5cm}|}
        \hline
        \textbf{Method} & \textbf{Number} & \textbf{Parameters type} \\
        \hline
        RL                 & 1               & $N^{stop}$: Optimal early stopping \\
        KL+TV$_\epsilon$              & 1               & $\lambda$: Optimal regularization parameter \\
        B-RED              & 4               & $(\lambda_{\text{in}},\gamma_{\text{in}},\lambda_{
        \text{algo}},\gamma_{\text{algo}})$: Initialization/reconstruction parameters \\
        B-PnP              & 4               & $(\lambda_{in},\gamma_{in},\lambda_{algo},\gamma_{algo})$: Initialization/reconstruction parameters \\
        DEQ-S              & 0               & - \\
        DEQ-RED            & 0               & - \\
        \hline
    \end{tabular}
\end{table}

We conclude this section presenting in Table \ref{tab:computational_times_FP} the training times, when needed and the CPU times required to conclude a forward pass, i.e.~to perform the inference phase and reconstruct the image. 
The test was executed on the entire test dataset simultaneously, whenever the code implementation supports batch processing.

\begin{table}[h]
\centering
\caption{Computational times for the training procedure and inference. The first column reports the total training time for the data-driven approaches. The second column shows the time required to perform inference on the entire test set, while the third column indicates the average time per image on the test set.}
\label{tab:computational_times_FP}
\resizebox{\textwidth}{!}{%
\begin{tabular}{lccc}
\toprule
\textbf{Method} & \textbf{Training time (hours)} & \textbf{Time on test set (s)} & \textbf{Average time per image (s)} \\
\midrule
RL & -- & 0.55 & 0.019 \\
KL+TV & -- & 146.44 & 3.18 \\
B-RED & 216 & 1900.15 & 41.30 \\
B-PnP & 216 & 860.11 & 18.69 \\
DEQ-S & 105 & 1050.34 & 22.83 \\
DEQ-RED & 50 & 408.94 & 8.89 \\
\bottomrule
\end{tabular}%
}
\end{table}

 The computational times for training the DEQ-MD model were obtained following the modalities described in Section \ref{loss function}. Values are averaged over the number of noise levels and tasks described in Section \ref{numerical experiments}. For both DEQ-S and DEQ-RED, time per epoch is higher in early epochs and then decreases until it stabilizes. An explanation of this behavior could be that in early epochs a lot of iterations are needed to get convergence in the forward pass (see the Adjoint-based case in Figure \ref{number iter for different init}) while later this effect is mitigated. This could be further improved using  a pre-training strategy (see Appendix \ref{pre-training}). For both B-RED and B-PnP, training times are reported w.r.t.~the training of the Bregman Denoiser, whose training modalities can be found in \cite{Hurault}.

 Regarding inference time, model-based RL stands out as the fastest algorithm, exhibiting a significant speed advantage over other methods. Among the learned approaches, DEQ-RED offers competitive CPU time. This efficiency stems from the use of the backtracking strategy, which is effectively implemented to leverage GPU parallelization and batch processing. Note that, however, such time gain is not observed using the DEQ-S method as the CNN employed requires more computation time within its fully-connected layers.

\subsection{Generalization to different tasks/noise levels}
In this section we discuss the generalization property of the DEQ-MD approaches presented in two respects:
\begin{itemize}
    \item Does a DEQ-MD model trained on a level of noise $\alpha$ generalize on a level of noise $\alpha'\neq \alpha$?
    \item Does a DEQ-MD model trained on a deblurring task using a certain convolution kernel generalize to a different task, such as, e.g.~ a super resolution problem?
\end{itemize}
For the following discussion, we focus only on the DEQ-RED method, which shows superior performances w.r.t DEQ-S as extensively discussed before.
\subsubsection{Generalization to different noise levels} \label{gen: diff level of noise}

%We recall that the forward pass of the DEQ-MD model is equivalent to the minimization of the functional $\Psi$ \eqref{Psi fun}. 
In order to generalize a DEQ-MD model trained on a noise level $\alpha$ to a different noise level $\alpha'\neq \alpha$ not seen during training, we consider the following $\lambda$-penalized instance of the DEQ-MD approach:
\begin{equation}
     f^\lambda _{\theta}(x;y)=\Pi_{[0,a]^n}\left(\frac{x}{1+\tau x (\nabla \text{KL}(y,Ax)+\lambda\nabla R_{\theta}(x))}\right), \label{eq:f_theta_with_lam}
\end{equation}
which is associated to the following minimization problem:
\begin{equation}
    \Psi_\lambda (x)= \text{KL}(y, Ax) + \lambda R_{\theta}(x) + \iota_{[0,a]^n}(x) 
\end{equation}
for which all convergence guarantees discussed in Section \ref{sec: conv result} hold.
The parameter $\lambda$ plays here the role of an actual regularization parameter whose magnitude should weight more/less the regularization against the KL data term depending on the ration $\alpha/\alpha'$. A natural choice thus consists in considering as $\lambda= \frac{\alpha}{\alpha'}$. To evaluate the effect of this choice, we tested all considered kernels by evaluating DEQ-MD models trained on a single noise level $\alpha \in {100, 60, 40}$ at different test noise levels $\alpha'\neq \alpha$, using the rule above to select $\lambda$ in \eqref{eq:f_theta_with_lam}. Results are reported in Table \ref{table: gen to diff levels of noise}.

\begin{table}[h!]
\centering
\caption{Generalization of DEQ-MD w.r.t.~different noise levels.}
\label{table: gen to diff levels of noise}
\resizebox{\textwidth}{!}{%
\begin{tabular}{c|ccc|ccc|ccc}
\toprule
\multirow{2}{*}{$\alpha$ training} & \multicolumn{3}{c|}{$\alpha$ inference = 100} & \multicolumn{3}{c|}{$\alpha$ inference = 60} & \multicolumn{3}{c}{$\alpha$ inference = 40} \\
\cmidrule(lr){2-4} \cmidrule(lr){5-7} \cmidrule(lr){8-10}
 & PSNR $\uparrow$ & SSIM $\uparrow$ & LPIPS $\downarrow$ & PSNR $\uparrow$ & SSIM $\uparrow$ & LPIPS $\downarrow$ & PSNR $\uparrow$  & SSIM $\uparrow$ & LPIPS $\downarrow$ \\
\midrule
100 & \underline{25.4245} & \underline{0.7091} & \underline{0.3058} & 24.7776 & 0.6778 & 0.3465 & 24.2627 & 0.6513 & 0.3816 \\
60  & \textbf{25.4828} & \textbf{0.7112} & \textbf{0.3020} & \textbf{24.9087} & \underline{0.6826} & \underline{0.3354} & \underline{24.3854} & \underline{0.6568} & \underline{0.3689} \\
40  & 25.4006 & 0.7079 & 0.3077 & \underline{24.8800} & \textbf{0.6841} & \textbf{0.3335} & \textbf{24.4586} & \textbf{0.6626} & \textbf{0.3571} \\
\bottomrule
\end{tabular}%
}
\end{table}

\begin{figure}[h!]
    \centering
    % Prima riga: 2 immagini
    \begin{subfigure}{0.25\textwidth}
    \begin{overpic}[width=\linewidth]{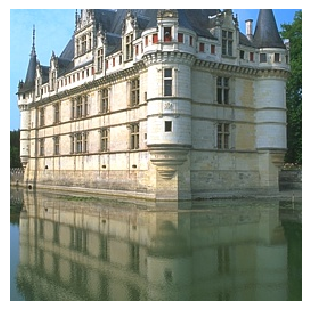}

        \put(2,2){\includegraphics[width=1cm]{images/starfish_motion_7_alpha_60/kernel_motion_7}}
    \end{overpic}
        \caption{Clean\\
        \mbox{\footnotesize(PSNR$\uparrow$,SSIM$\uparrow$,LPIPS$\downarrow$)}}
    \end{subfigure}
    \hspace{0.05\textwidth} % spazio tra le immagini
    \begin{subfigure}{0.25\textwidth}
        \includegraphics[width=\linewidth]{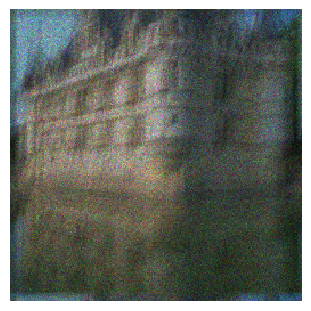}
        \caption{Observed\\
        (14.564
        0.186,
        0.514)}
    \end{subfigure}
    
    \vspace{0.5em} % spazio tra le righe

    % Seconda riga: 3 immagini
    \begin{subfigure}{0.25\textwidth}
        \includegraphics[width=\linewidth]{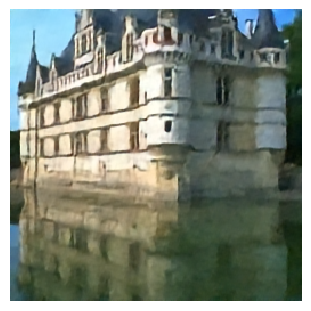}
        \caption{$\alpha^{\text{train}}=100$,
        $\lambda =1 $ \\
        (\underline{23.353}, \textbf{0.709}, \underline{0.312})}
    \end{subfigure}
    \hspace{0.02\textwidth} % piccolo spazio tra le immagini
    \begin{subfigure}{0.25\textwidth}
        \includegraphics[width=\linewidth]{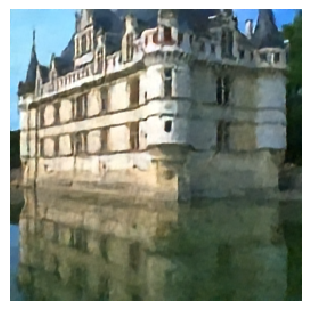}
        \caption{$\alpha^{\text{train}}=60$,
        $\lambda = 0.6$\\
        (23.317,
        0.704,
        0.313)}
    \end{subfigure}
    \hspace{0.02\textwidth}
    \begin{subfigure}{0.25\textwidth}
        \includegraphics[width=\linewidth]{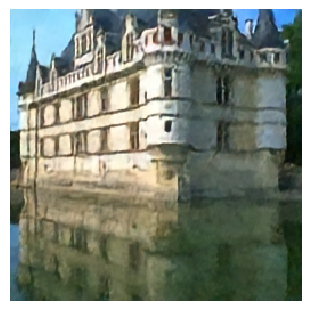}
        \caption{$\alpha^{\text{train}}=40$, $\lambda = 0.4$\\
        (\textbf{23.395},
        \underline{0.709},
        \textbf{0.305})}
    \end{subfigure}

    \caption{Generalization study on different levels of noise. Test is performed  on a deblurring task with motion blur with $\alpha^{\text{test}}=100$.}
    \label{fig: gen on level of noise with motion blur}
\end{figure}

\begin{figure}[h!]
    \centering
    % Prima riga: 2 immagini
     \begin{subfigure}{0.25\textwidth}
    \begin{overpic}[width=\linewidth]{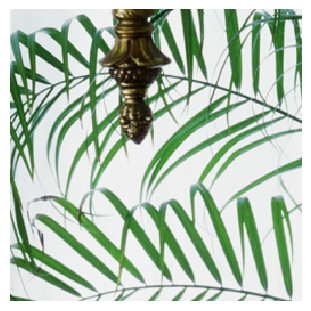}

        \put(2,2){\includegraphics[width=1cm]{images/gaussian_kernel.png}}
    \end{overpic}
        \caption{Clean\\
\mbox{\footnotesize(PSNR$\uparrow$,SSIM$\uparrow$,LPIPS$\downarrow$)}}
    \end{subfigure}
    \hspace{0.05\textwidth} % spazio tra le immagini
    \begin{subfigure}{0.25\textwidth}
        \includegraphics[width=\linewidth]{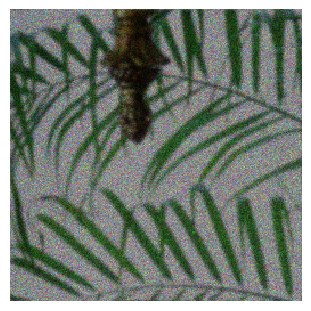}
        \caption{Observed\\
        (17.586,
        0.503,
        0.5435)}
    \end{subfigure}
    
    \vspace{0.5em} % spazio tra le righe

    % Seconda riga: 3 immagini
    \begin{subfigure}{0.25\textwidth}
        \includegraphics[width=\linewidth]{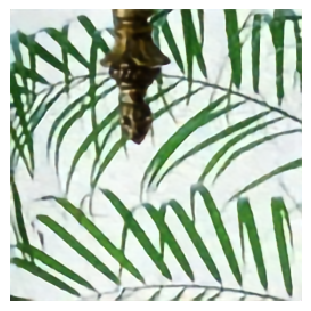}
        \caption{$\alpha^{\text{train}}=100$, $\lambda =5/3 $ \\
        (24.514,
        \underline{0.896},
        0.090)}
    \end{subfigure}
    \hspace{0.02\textwidth} % piccolo spazio tra le immagini
    \begin{subfigure}{0.25\textwidth}
        \includegraphics[width=\linewidth]{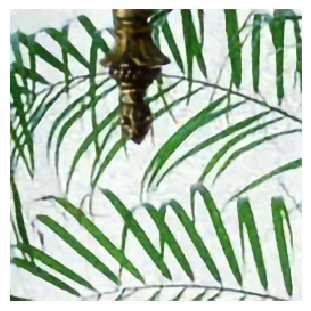}
        \caption{$\alpha^{\text{train}}=60$,
        $\lambda = 1$\\
        (\underline{24.793},
        0.889,
        \textbf{0.074})}
    \end{subfigure}
    \hspace{0.02\textwidth}
    \begin{subfigure}{0.25\textwidth}
        \includegraphics[width=\linewidth]{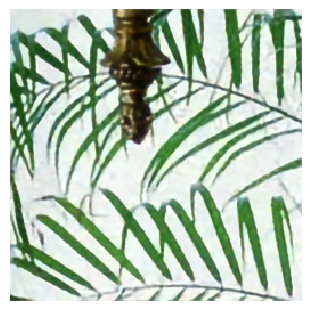}
        \caption{$\alpha^{\text{train}}=40$, 
        $\lambda =2/3$\\
        (\textbf{24.885},
        \textbf{0.896},
        \underline{0.083})}
    \end{subfigure}

    \caption{Generalization study on different levels of noise. Test is performed  on a deblurring task with Gaussian blur with  $\alpha^{\text{test}}=60$.}
    \label{fig: gen on diff level of noise with gauss blur}
\end{figure}

For each noise level, we refer to the metrics value on the diagonal as the reference ones, since in those cases the model has been tested on the same level of noise seen during training. We note that upon the choice of $\lambda$ discussed above, DEQ-MD appear to be easily generalizable to unseen noise levels, especially when the inference is performed on a lower intensity of Poisson noise w.r.t the one in the training. In particular, for $\alpha=100$, the best model appears to be the one trained on $\alpha=60$ with $\alpha=40$ providing performances close to the reference values. The same consideration holds for the inference with $\alpha=60$ and the model trained on $\alpha=40$. 
%This analysis suggests to that it might be possible to train only one DEQ-MD for all level noise, choosing in the training a sufficient high level of Poisson noise. 
While the proposed selection rule for $\lambda$ might not be the optimal one, it sill provides a fairly simple criterion, based only on the knowledge of the noise levels, which streamlines hyperparameter selection. Figure \ref{fig: gen on level of noise with motion blur} and \ref{fig: gen on diff level of noise with gauss blur} visually demonstrate the effectiveness of the proposed approach.

\newpage
\subsubsection{Generalization to super-resolution}

\begin{figure}[h!] 
    \centering 
    % Prima riga: 2 immagini
    \begin{subfigure}{0.25\textwidth}
    \begin{overpic}[width=\linewidth]{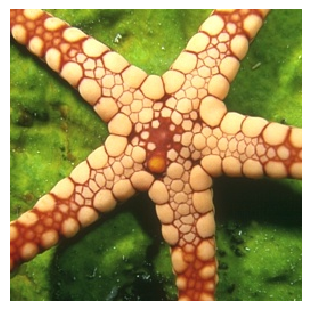}

        \put(2,2){\includegraphics[width=1cm]{images/gaussian_kernel.png}}
    \end{overpic}
        \caption{Clean\\
        \mbox{\footnotesize (PSNR$\uparrow$,
        SSIM$\uparrow$,
        LPIPS$\downarrow$)}}
    \end{subfigure}
    \hspace{0.02\textwidth} % spazio tra le immagini
   \begin{subfigure}{0.25\textwidth}
    \centering
        \raisebox{0.2\textwidth}{\includegraphics[width=0.5\linewidth]{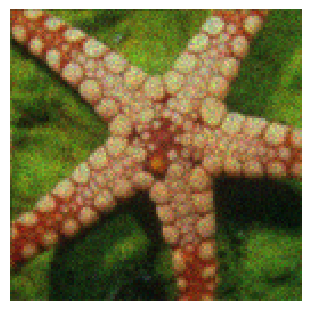}}
        \caption{Observed\\
        (18.989, 0.635,
        0.445)}
    \end{subfigure}
    \hspace{0.02\textwidth} % spazio tra le immagini
     \begin{subfigure}{0.25\textwidth}
        \includegraphics[width=\linewidth]{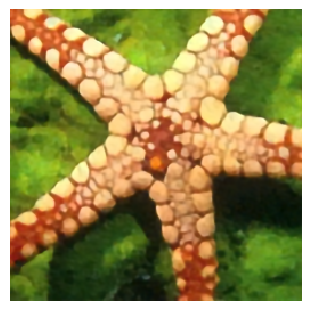}
        \caption{DEQ-RED\\
        (25.897,
        0.769,
        0.204)}
    \end{subfigure}
    
    \vspace{0.5em} % spazio tra le righe

    % Seconda riga: 3 immagini
    \begin{subfigure}{0.25\textwidth}
    \begin{overpic}[width=\linewidth]{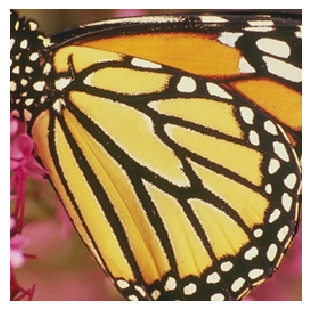}

        \put(2,2){\includegraphics[width=1cm]{images/motion_kernel_1.png}}
    \end{overpic}
\caption{Clean\\
        \mbox{\footnotesize (PSNR$\uparrow$,
        SSIM$\uparrow$,
        LPIPS$\downarrow$)}}
    \end{subfigure}
    \hspace{0.02\textwidth} % piccolo spazio tra le immagini
    \begin{subfigure}{0.25\textwidth}
    \centering
        \raisebox{0.2\textwidth}{\includegraphics[width=0.5\linewidth]{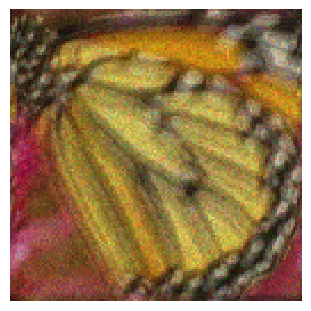}}
        \caption{Observed\\
        (14.495,
        0.255,
        0.704)}
    \end{subfigure}
    \hspace{0.02\textwidth}
    \begin{subfigure}{0.25\textwidth}
        \includegraphics[width=\linewidth]{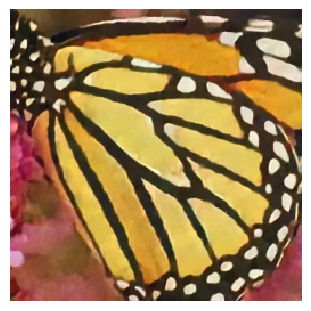}
        \caption{DEQ-RED\\
        (24.208,
        0.809,
        0.113)}
    \end{subfigure}

    \caption{Generalization test to $\times 2$ super-resolution task.
    First row: Gaussian blur and $\alpha=60$.
    Second row: Motion blur and $\alpha=100$.
    The value of the metrics for the observed images have been calculated between the clean images and image $(SA)^\top y$.}
    \label{figure: SR}
\end{figure}

In this section, we study the generalization ability of DEQ-MD models trained on a specific image deblurring task to a related but distinct task: Super-Resolution under Poisson noise. In this setting, an image $y\in\mathbb{R}^m$ is obtained as a degraded version of an image $x\in\mathbb{R}^n$ with $n=s^2m$, according to:
\begin{equation}
    y=\mathrm{Poiss}(\alpha S A x),
\end{equation}
where $S \in \mathbb{R}^{m \times n}$ is a down-sampling operator, $s$ the scale factor and $A\in\mathbb{R}^{n\times n}$ is a convolution operator associated with the kernels shown in Figure~\ref{fig:kernels}. The operator $S$ extracts the top-left pixel from each $s\times s$ block.  Figure~\ref{figure: SR} presents results for the case $s=2$. Even without task-specific training, DEQ-RED exhibits notable generalization to this Super-Resolution setting.

We observed a degradation of the results for more challenging scale factors $s\geq 3$. 

\begin{remark}[The Importance of normalization]
   We stress that in order to get good super-resolution results as the ones in Figure \ref{figure: SR}, it is important to normalize the forward operator $T=SA$ so as to have $\|T^{\top}T\|_2=1$. In practice, it is numerically important to normalize the operator $T$, in order to ensure that its overall energy is consistent across different resolutions. This normalization prevents the inversion process from being biased or ill-conditioned due to scale discrepancies introduced by the blur-and-subsample operations. Without such normalization, the reconstructed results tend to be of lower quality, with a noticeable loss of image details and fine structures. The normalization used is available at \url{https://github.com/christiandaniele/DEQ-MD}.
\end{remark}

\section{Conclusions and Outlook}
\label{sec:conclusions}

We proposed a novel Deep Equilibrium Model (DEQ) framework, named DEQ-MD, for solving linear Poisson inverse problems, leveraging Mirror Descent as a geometry-aware optimization backbone. Unlike classical Plug-and-Play or unrolled architectures, DEQ-MD learns regularization functionals end-to-end through the fixed point of the MD algorithm with a geometry tailored to the Kullback–Leibler divergence. This approach allows us to bypass standard smoothness and contractivity assumptions typically required in DEQ theory, while still guaranteeing convergence of the reconstruction procedure. Our analysis avoids operator contractivity and instead relies on non-convex optimization principles, enabling greater flexibility in the network design and interpretability through the minimization of an explicit functional. In this context, we established refined convergence results extending \cite[Theorem 2]{Bolte2007} within the Kurdyka–Łojasiewicz framework, providing a broadly applicable tool for handling Kullback–Leibler divergences in the presence of non-convex regularization. Beyond the theoretical framework, our method yields competitive reconstructions using significantly lighter architectures, while offering faster inference and full parameter-free operation at test time. Though developed for the Poisson setting, our formulation is general and potentially applicable to other problems and geometries via suitable choices of the underlying Bregman potential (e.g., phase retrieval \cite{godeme2024stablephaseretrievalmirror}).

Looking ahead, we aim to extend convergence guarantees for the Jacobian-Free Backpropagation (JFB) scheme in more general, non-contractive and non-smooth settings, and to integrate self-supervised training strategies. This would further increase the applicability of our approach in domains like microscopy and medical imaging, where ground-truth data is scarce or unavailable.

\section*{Acknowledgments}
The work of C. Daniele and L. Calatroni was supported by the ANR JCJC project TASKABILE ANR-22-CE48-0010. The work of C. Daniele, L. Calatroni and S. Villa was supported by the funding received from the European Research Council (ERC) Starting project MALIN under the European Union’s Horizon Europe programme (grant agreement No. 101117133).\\
S. Villa acknowledges the support of the European Commission (grant TraDE-OPT 861137), of the European Research Council (grant SLING 819789), the US Air Force Office of Scientific Research (FA8655-22-1-7034), the Ministry of Education, University and Research (PRIN 202244A7YL project ``Gradient Flows and Non-Smooth Geometric Structures with Applications to Optimization and Machine Learning’’), and the project PNRR FAIR PE0000013 - SPOKE 10. The research by S. Villa. has been supported by the MIUR Excellence Department Project awarded to Dipartimento di Matematica, Università di Genova, CUP D33C23001110001. S. Villa is a member of the Gruppo Nazionale per l’Analisi Matematica, la Probabilità e le loro Applicazioni (GNAMPA) of the Istituto Nazionale di Alta Matematica (INdAM). \\
S. Vaiter thanks PEPR PDE-AI (ANR-23-PEIA-0004), the chair 3IA BOGL through the 3IA Cote d’Azur Investments (ANR-23-IACL-0001) and the ANR PRC MAD (ANR-24-CE23-1529).\\
This work represents only the view of the authors.
 The European Commission and the other organizations are not responsible for any use that may be made of the information it contains.

\bibliographystyle{plain} % Puoi cambiare 'plain' con lo stile desiderato
\bibliography{bibliografia} % Assicurati che il tuo file .bib si chiami 'bibliografia.bib'

@article{ICNNpaper,
  title={{W}hat's in a {P}rior? {L}earned {P}roximal {N}etworks for {I}nverse {P}roblems},
  author={Fang, Zhenghan and Buchanan, Sam and Sulam, Jeremias},
  journal={arXiv preprint arXiv:2310.14344},
  year={2023}
}

@inproceedings{BolteLePauwelsSilveti2021,
 author = {Bolte, J\'{e}r\^{o}me and Le, Tam and Pauwels, Edouard and Silveti-Falls, Tony},
 booktitle = {Advances in Neural Information Processing Systems},
 editor = {M. Ranzato and A. Beygelzimer and Y. Dauphin and P.S. Liang and J. Wortman Vaughan},
 pages = {13537--13549},
 publisher = {Curran Associates, Inc.},
 title = {Nonsmooth Implicit Differentiation for Machine-Learning and Optimization},
volume = {34},
 year = {2021}
}

@article{Bertero_2009,
  title={{I}mage {D}eblurring with {P}oisson {D}ata: {F}rom {C}ells to {G}alaxies},
  author={Bertero, Mario and Boccacci, Patrizia and Desider{\`a}, Gabriele and Vicidomini, Giuseppe},
  journal={Inverse Problems},
  volume={25},
  number={12},
  pages={123006},
  year={2009},
  publisher={IOP Publishing}
}

@article{Dupe2009,
  title={{A} {P}roximal {I}teration for {D}econvolving {P}oisson {N}oisy {I}mages {U}sing {S}parse {R}epresentations},
  author={Dup{\'e}, Fran{\c{c}}ois-Xavier and Fadili, Jalal M. and Starck, Jean-Luc},
  journal={IEEE Transactions on Image Processing},
  volume={18},
  number={2},
  pages={310--321},
  year={2009},
  publisher={IEEE}
}

@article{conservativederivatives,
  title={{A}utomatic {d}ifferentiation of {n}onsmooth {i}terative {a}lgorithms},
  author={Bolte, J{\'e}r{\^o}me and Pauwels, Edouard and Vaiter, Samuel},
  journal={Advances in Neural Information Processing Systems},
  volume={35},
  pages={26404--26417},
  year={2022}
}

@book{Beck,
  author={Beck, A.},
  title={{F}irst-{O}rder {O}ptimisation {M}ethods},
  publisher={Society for Industrial and Applied Mathematics (SIAM)},
  year={2017},
}

@inproceedings{JFB,
  title={{J}{F}{B}: {J}acobian-{F}ree {B}ackpropagation for {I}mplicit {N}etworks},
  author={Fung, Samy Wu and Heaton, Howard and Li, Qiuwei and McKenzie, Daniel and Osher, Stanley and Yin, Wotao},
  booktitle={{P}roceedings of the {A}{A}{A}{I} {C}onference on {A}rtificial {I}ntelligence},
  volume={36},
  number={6},
  pages={6648--6656},
  year={2022}
}

@article{Attouch2013,
  author={Attouch, H. and Bolte, J. and Svaiter, B. F.},
  title={{C}onvergence of {D}escent {M}ethods for {S}emi-{A}lgebraic and {T}ame {P}roblems: {P}roximal {A}lgorithms, {F}orward--{B}ackward {S}plitting, and {R}egularized {G}auss--{S}eidel {M}ethods},
  journal={Mathematical Programming},
  volume={137},
  number={1-2},
  pages={91--129},
  year={2013},
}

@article{Bolte2018,
  title={{F}irst {O}rder {M}ethods {B}eyond {C}onvexity and {L}ipschitz {G}radient {C}ontinuity with {A}pplications to {Q}uadratic {I}nverse {P}roblems},
  author={Bolte, J{\'e}r{\^o}me and Sabach, Shoham and Teboulle, Marc and Vaisbourd, Yakov},
  journal={SIAM Journal on Optimization},
  volume={28},
  number={3},
  pages={2131--2151},
  year={2018},
  publisher={SIAM}
}

@article{Figuereido2010,
  title={{R}estoration of {P}oissonian {I}mages {U}sing {A}lternating {D}irection {O}ptimization},
  author={Figueiredo, M{\'a}rio AT and Bioucas-Dias, Jos{\'e} M},
  journal={IEEE Transactions on Image Processing},
  volume={19},
  number={12},
  pages={3133--3145},
  year={2010},
  publisher={IEEE}
}

@article{Bonettini2018a,
  title={{I}nertial {V}ariable {M}etric {T}echniques for the {I}nexact {F}orward--{B}ackward {A}lgorithm},
  author={Bonettini, Silvia and Rebegoldi, Simone and Ruggiero, Valeria},
  journal={SIAM Journal on Scientific Computing},
  volume={40},
  number={5},
  pages={A3180--A3210},
  year={2018},
  publisher={SIAM}
}

@article{DiSerafino2021,
  title={{D}irectional {T}{G}{V}-{B}ased {I}mage {R}estoration {U}nder {P}oisson {N}oise},
  author={di Serafino, Daniela and Landi, Germana and Viola, Marco},
  journal={Journal of Imaging},
  volume={7},
  number={6},
  pages={99},
  year={2021},
  publisher={MDPI}
}

@article{Chambolle2018,
  title={{S}tochastic {P}rimal-{D}ual {H}ybrid {G}radient {A}lgorithm {W}ith {A}rbitrary {S}ampling and {I}maging {A}pplications},
  author={Chambolle, Antonin and Ehrhardt, Matthias J and Richt{\'a}rik, Peter and Schonlieb, Carola-Bibiane},
  journal={SIAM Journal on Optimization},
  volume={28},
  number={4},
  pages={2783--2808},
  year={2018},
  publisher={SIAM}
}

@article{Lucy1974,
  title={{A}n {I}terative {T}echnique for the {R}ectification of {O}bserved {D}istributions},
  author={Lucy, Leon B.},
  journal={Astronomical Journal},
  volume={79},
  pages={745},
  year={1974}
}

@article{Byrne1993,
  title={{I}terative {I}mage {R}econstruction {A}lgorithms {B}ased on {C}ross-{E}ntropy {M}inimization},
  author={Byrne, Charles L.},
  journal={IEEE Transactions on Image Processing},
  volume={2},
  number={1},
  pages={96--103},
  year={1993},
  publisher={IEEE}
}

@article{Attouch2010,
  author={Attouch, H. and Bolte, J. and Redont, P. and Soubeyran, A.},
  title={{P}roximal {A}lternating {M}inimisation and {P}rojection {M}ethods for {N}onconvex {P}roblems: {A}n {A}pproach {B}ased on the {K}urdyka-{Ł}ojasiewicz Inequality},
  journal={Mathematics of Operations Research},
  volume={35},
  number={2},
  pages={438--457},
  year={2010},
}

@article{Loja,
  author={{\L}ojasiewicz, S.},
  title={{T}riangulation of {S}emi-{A}nalytic {S}ets},
  journal={Annali della Scuola Normale Superiore di Pisa, Classe di Scienze},
  series={3e série},
  volume={18},
  number={4},
  pages={449--474},
  year={1964},
}

@phdthesis{Perrin2020,
  author={Perrin, M. S.},
  title={{S}emialgebraic {G}eometry},
  school={University of Sydney},
  year={2020},
}

@article{Bolte2007,
  author={Bolte, J. and Daniilidis, A. and Lewis, A.},
  title={{T}he {Ł}ojasiewicz Inequality for Nonsmooth Subanalytic Functions With Applications to Subgradient Dynamical Systems},
  journal={SIAM Journal on Optimization},
  volume={17},
  number={4},
  pages={1205--1223},
  year={2007},
}

@book{Shiota2012,
  author={Shiota, M.},
  title={{G}eometry of {S}ubanalytic and {S}emialgebraic {S}ets},
  volume={150},
  publisher={Springer Science \& Business Media},
  year={2012},
}

@article{nolip,
  title={{A} {D}escent {L}emma {B}eyond {L}ipschitz Gradient Continuity: First-Order Methods Revisited and Applications},
  author={Bauschke, Heinz H. and Bolte, J{\'e}r{\^o}me and Teboulle, Marc},
  journal={Mathematics of Operations Research},
  volume={42},
  number={2},
  pages={330--348},
  year={2017},
  publisher={Informs}
}

@book{nemirovskii_yudin_1983,
  author={Nemirovski{\u\i}, Arkadiĭ Semenovich and I{\u U}din, David Berkovich},
  title={{P}roblem {C}omplexity and {M}ethod {E}fficiency in {O}ptimization},
  series={Wiley-Interscience Series in Discrete Mathematics},
  publisher={John Wiley \& Sons},
  address={Chichester / New York},
  year={1983},
  isbn={0-471-10345-4}
}

@article{Hurault,
  title={{Convergent {B}regman {P}lug-and-{P}lay {I}mage {R}estoration for {P}oisson {I}nverse {P}roblems}},
  author={Hurault, Samuel and Kamilov, Ulugbek and Leclaire, Arthur and Papadakis, Nicolas},
  journal={Advances in Neural Information Processing Systems},
  volume={36},
  pages={27251--27280},
  year={2023}
}

@book{BauschkeCombettes2011,
  title     = {{Convex Analysis and Monotone Operator Theory in Hilbert Spaces}},
  author    = {Bauschke, Heinz H. and Combettes, Patrick L.},
  year      = {2011},
  publisher = {Springer},
  address   = {New York},
  series    = {CMS Books in Mathematics},
  isbn      = {978-1-4419-9466-0}
}

@article{Teboulle2018,
	Author = {Teboulle, Marc},
	Journal = {Mathematical Programming},
	Number = {1},
	Pages = {67--96},
	Title = {A simplified view of first order methods for optimization},
	Volume = {170},
	Year = {2018}}

@article{Kamilov2023,
  title={{P}lug-and-{P}lay {M}ethods for {I}ntegrating {P}hysical and {L}earned {M}odels in {C}omputational {I}maging: {T}heory, {A}lgorithms, and {A}pplications},
  author={Kamilov, Ulugbek S and Bouman, Charles A and Buzzard, Gregery T and Wohlberg, Brendt},
  journal={IEEE Signal Processing Magazine},
  volume={40},
  number={1},
  pages={85--97},
  year={2023},
  publisher={IEEE}
}

@inproceedings{meinhardt2017learning,
  title={{L}earning {P}roximal {O}perators: {U}sing {D}enoising {N}etworks for {R}egularizing {I}nverse {I}maging {P}roblems},
  author={Meinhardt, Tim and Moller, Michael and Hazirbas, Caner and Cremers, Daniel},
  booktitle={{P}roceedings of the {I}{E}{E}{E} {I}nternational {C}onference on {C}omputer {V}ision},
  pages={1781--1790},
  year={2017}
}

@article{laumont2022bayesian,
  title={{B}ayesian {I}maging {U}sing {P}lug \& {P}lay {P}riors: {W}hen {L}angevin {M}eets {T}weedie},
  author={Laumont, R{\'e}mi and Bortoli, Valentin De and Almansa, Andr{\'e}s and Delon, Julie and Durmus, Alain and Pereyra, Marcelo},
  journal={SIAM Journal on Imaging Sciences},
  volume={15},
  number={2},
  pages={701--737},
  year={2022},
  publisher={SIAM}
}

@article{Hurault1,
  title={{G}radient {S}tep {D}enoiser for {C}onvergent {P}lug-and-{P}lay},
  author={Hurault, Samuel and Leclaire, Arthur and Papadakis, Nicolas},
  journal={arXiv preprint arXiv:2110.03220},
  year={2021}
}

@article{Unrolling,
  title={{A}lgorithm {U}nrolling: {I}nterpretable, {E}fficient {D}eep {L}earning for {S}ignal and {I}mage {P}rocessing},
  author={Monga, Vishal and Li, Yuelong and Eldar, Yonina C},
  journal={IEEE Signal Processing Magazine},
  volume={38},
  number={2},
  pages={18--44},
  year={2021},
  publisher={IEEE}
}

@article{AdlerOktem2018,
  title={{L}earned {P}rimal-{D}ual {R}econstruction},
  author={Adler, Jonas and {\"O}ktem, Ozan},
  journal={IEEE Transactions on Medical Imaging},
  volume={37},
  number={6},
  pages={1322--1332},
  year={2018},
  publisher={IEEE}
}

@article{Willett,
  title={{D}eep {E}quilibrium {A}rchitectures for {I}nverse {P}roblems in {I}maging},
  author={Gilton, Davis and Ongie, Gregory and Willett, Rebecca},
  journal={IEEE Transactions on Computational Imaging},
  volume={7},
  pages={1123--1133},
  year={2021},
  publisher={IEEE}
}

@article{primi_DEQ,
  title={{D}eep {E}quilibrium {M}odels},
  author={Bai, Shaojie and Kolter, J. Zico and Koltun, Vladlen},
  journal={Advances in Neural Information Processing Systems},
  volume={32},
  year={2019}
}

@article{Bonettini_2009,
  title={{A} {S}caled {G}radient {P}rojection {M}ethod for {C}onstrained {I}mage {D}eblurring},
  author={Bonettini, Silvia and Zanella, Riccardo and Zanni, Luca},
  journal={Inverse Problems},
  volume={25},
  number={1},
  pages={015002},
  year={2008},
  publisher={IOP Publishing}
}

@article{godeme2024stablephaseretrievalmirror,
  title={{S}table {P}hase {R}etrieval {W}ith {M}irror {D}escent},
  author={Godeme, Jean-Jacques and Fadili, Jalal and Amra, Claude and Zerrad, Myriam},
  journal={arXiv preprint arXiv:2405.10754},
  year={2024}
}

@inproceedings{Shabili2022,
  title={{Bregman {P}lug-and-{P}lay {P}riors}},
  author={Al-Shabili, Abdullah H. and Xu, Xiaojian and Selesnick, Ivan and Kamilov, Ulugbek S.},
  booktitle={2022 {I}{E}{E}{E} {I}nternational {C}onference on {I}mage {P}rocessing ({I}{C}{I}{P})},
  pages={241--245},
  year={2022},
  organization={IEEE}
}

@article{klatzer2025,
  title={{E}fficient {B}ayesian {C}omputation {U}sing {P}lug-and-{P}lay {P}riors for {P}oisson {I}nverse {P}roblems},
  author={Klatzer, Teresa and Melidonis, Savvas and Pereyra, Marcelo and Zygalakis, Konstantinos C},
  journal={arXiv preprint arXiv:2503.16222},
  year={2025}
}

@article{KLsmoothed,
  title={{T}his {I}s {S}{P}{I}{R}{A}{L}-{T}{A}{P}: {S}parse {P}oisson {I}ntensity {R}econstruction {A}lgorithms—{T}heory and {P}ractice},
  author={Harmany, Zachary T and Marcia, Roummel F and Willett, Rebecca M},
  journal={IEEE Transactions on Image Processing},
  volume={21},
  number={3},
  pages={1084--1096},
  year={2011},
  publisher={IEEE}
}

@article{Samuel,
  title={{O}ne-{S}tep {D}ifferentiation of {I}terative {A}lgorithms},
  author={Bolte, J{\'e}r{\^o}me and Pauwels, Edouard and Vaiter, Samuel},
  journal={Advances in Neural Information Processing Systems},
  volume={36},
  pages={77089--77103},
  year={2023}
}

@inproceedings{Levin,
  title={{U}nderstanding and {E}valuating {B}lind {D}econvolution {A}lgorithms},
  author={Levin, Anat and Weiss, Yair and Durand, Fredo and Freeman, William T.},
  booktitle={2009 {I}{E}{E}{E} {C}onference on {C}omputer {V}ision and {P}attern {R}ecognition},
  pages={1964--1971},
  year={2009},
  organization={IEEE}
}

@article{Adam,
  title={{A}dam: {A} {M}ethod for {S}tochastic {O}ptimization},
  author={Kingma, Diederik P.},
  journal={arXiv preprint arXiv:1412.6980},
  year={2014}
}

@article{RelativeSmoothness,
  title={{R}elatively {S}mooth {C}onvex {O}ptimization by {F}irst-{O}rder {M}ethods, and {A}pplications},
  author={Lu, Haihao and Freund, Robert M and Nesterov, Yurii},
  journal={SIAM Journal on Optimization},
  volume={28},
  number={1},
  pages={333--354},
  year={2018},
  publisher={SIAM}
}

@article{RebegoldiCalatroni2022,
  title={{S}caled, {I}nexact, and {A}daptive {G}eneralized {F}{I}{S}{T}{A} for {S}trongly {C}onvex {O}ptimization},
  author={Rebegoldi, Simone and Calatroni, Luca},
  journal={SIAM Journal on Optimization},
  volume={32},
  number={3},
  pages={2428--2459},
  year={2022},
  publisher={SIAM}
}

@article{Richardson72,
  title={{B}ayesian-{B}ased {I}terative {M}ethod of {I}mage {R}estoration},
  author={Richardson, William Hadley},
  journal={Journal of the Optical Society of America},
  volume={62},
  number={1},
  pages={55--59},
  year={1972},
  publisher={Optical Society of America}
}

@book{Rockafellar,
  title={{C}onvex {A}nalysis},
  author={Rockafellar, R. Tyrrell},
  volume={28},
  year={1997},
  publisher={Princeton University Press}
}

@article{Romano2017,
  title={{T}he {L}ittle {E}ngine {T}hat {C}ould: {R}egularization by {D}enoising ({R}{E}{D})},
  author={Romano, Yaniv and Elad, Michael and Milanfar, Peyman},
  journal={SIAM Journal on Imaging Sciences},
  volume={10},
  number={4},
  pages={1804--1844},
  year={2017},
  publisher={SIAM}
}

@article{TV,
  title={{N}onlinear {T}otal {V}ariation {B}ased {N}oise {R}emoval {A}lgorithms},
  author={Rudin, Leonid I and Osher, Stanley and Fatemi, Emad},
  journal={Physica D: Nonlinear Phenomena},
  volume={60},
  number={1-4},
  pages={259--268},
  year={1992},
  publisher={Elsevier}
}

@inproceedings{Salimans,
  title={{S}hould {E}{B}{M}s {M}odel the {E}nergy or the {S}core?},
  author={Salimans, Tim and Ho, Jonathan},
  booktitle={{E}nergy {B}ased {M}odels {W}orkshop-{I}{C}{L}{R} 2021},
  year={2021}
}

@article{Villa,
  title={{P}roximal {G}radient {M}ethods for {M}achine {L}earning and {I}maging},
  author={Salzo, Saverio and Villa, Silvia},
  journal={Harmonic and Applied Analysis: From Radon Transforms to Machine Learning},
  pages={149--244},
  year={2021},
  publisher={Springer}
}

@inproceedings{Sawatzky2009,
  title={{T}otal {V}ariation {P}rocessing of {I}mages {W}ith {P}oisson {S}tatistics},
  author={Sawatzky, Alex and Brune, Christoph and M{\"u}ller, Jahn and Burger, Martin},
  booktitle={{I}nternational {C}onference on {C}omputer {A}nalysis of {I}mages and {P}atterns},
  pages={533--540},
  year={2009},
  organization={Springer}
}

@article{Shepp1982,
  title={{M}aximum {L}ikelihood {R}econstruction for {E}mission {T}omography},
  author={Shepp, Lawrence A and Vardi, Yehuda},
  journal={IEEE Transactions on Medical Imaging},
  volume={1},
  number={2},
  pages={113--122},
  year={2007},
  publisher={IEEE}
}

@article{sun2018online,
  title={{A}n {O}nline {P}lug-and-{P}lay {A}lgorithm for {R}egularized {I}mage {R}econstruction},
  author={Sun, Yu and Wohlberg, Brendt and Kamilov, Ulugbek S},
  journal={IEEE Transactions on Computational Imaging},
  volume={5},
  number={3},
  pages={395--408},
  year={2019},
  publisher={IEEE}
}

@article{sun2021scalable,
  title={{S}calable {P}lug-and-{P}lay {A}{D}{M}{M} {W}ith {C}onvergence {G}uarantees},
  author={Sun, Yu and Wu, Zihui and Xu, Xiaojian and Wohlberg, Brendt and Kamilov, Ulugbek S},
  journal={IEEE Transactions on Computational Imaging},
  volume={7},
  pages={849--863},
  year={2021},
  publisher={IEEE}
}

@inproceedings{Terris,
  title={{B}uilding {F}irmly {N}onexpansive {C}onvolutional {N}eural {N}etworks},
  author={Terris, Matthieu and Repetti, Audrey and Pesquet, Jean-Christophe and Wiaux, Yves},
  booktitle={{I}{C}{A}{S}{S}{P} 2020-2020 {I}{E}{E}{E} {I}nternational {C}onference on {A}coustics, {S}peech and {S}ignal {P}rocessing ({I}{C}{A}{S}{S}{P})},
  pages={8658--8662},
  year={2020},
  organization={IEEE}
}

@inproceedings{Venkatakrishnan2013,
  title={{P}lug-and-{P}lay {P}riors for {M}odel {B}ased {R}econstruction},
  author={Venkatakrishnan, Singanallur V and Bouman, Charles A and Wohlberg, Brendt},
  booktitle={2013 {I}{E}{E}{E} {G}lobal {C}onference on {S}ignal and {I}nformation {P}rocessing},
  pages={945--948},
  year={2013},
  organization={IEEE}
}

@article{Anderson,
  title={{A}nderson {A}cceleration for {F}ixed-{P}oint {I}terations},
  author={Walker, Homer F and Ni, Peng},
  journal={SIAM Journal on Numerical Analysis},
  volume={49},
  number={4},
  pages={1715--1735},
  year={2011},
  publisher={SIAM}
}

@article{SSIM,
  title={{I}mage {Q}uality {A}ssessment: {F}rom {E}rror {V}isibility to {S}tructural {S}imilarity},
  author={Wang, Zhou and Bovik, Alan C and Sheikh, Hamid R and Simoncelli, Eero P},
  journal={IEEE Transactions on Image Processing},
  volume={13},
  number={4},
  pages={600--612},
  year={2004},
  publisher={IEEE}
}

@article{DnCNN,
  title={{B}eyond a {G}aussian {D}enoiser: {R}esidual {L}earning of {D}eep {C}{N}{N} for {I}mage {D}enoising},
  author={Zhang, Kai and Zuo, Wangmeng and Chen, Yunjin and Meng, Deyu and Zhang, Lei},
  journal={IEEE Transactions on Image Processing},
  volume={26},
  number={7},
  pages={3142--3155},
  year={2017},
  publisher={IEEE}
}

@ARTICLE{ByrneAccEMML1998,
  author={Byrne, C.L.},
  journal={IEEE Transactions on Image Processing}, 
  title={Accelerating the EMML algorithm and related iterative algorithms by rescaled block-iterative methods}, 
  year={1998},
  volume={7},
  number={1},
  pages={100-109},
  doi={10.1109/83.650854}}

@inproceedings{LPIPS,
  title={{T}he {U}nreasonable {E}ffectiveness of {D}eep {F}eatures as a {P}erceptual {M}etric},
  author={Zhang, Richard and Isola, Phillip and Efros, Alexei A and Shechtman, Eli and Wang, Oliver},
  booktitle={{P}roceedings of the {I}{E}{E}{E} {C}onference on {C}omputer {V}ision and {P}attern {R}ecognition},
  pages={586--595},
  year={2018}
}

@article{Zunino2023_ISMinverse,
  title={{R}econstructing the {I}mage {S}canning {M}icroscopy {D}ataset: {A}n {I}nverse {P}roblem},
  author={Zunino, Alessandro and Castello, Marco and Vicidomini, Giuseppe},
  journal={Inverse Problems},
  volume={39},
  number={6},
  pages={064004},
  year={2023},
  publisher={IOP Publishing}
}

@article{propertysubanalsets,
  author    = {E. Bierstone and P. D. Milman},
  title     = {Semianalytic and Subanalytic Sets},
  journal   = {Publications Mathématiques de l'IHÉS},
  volume    = {67},
  year      = {1988},
  pages     = {5--42},
}

@book{Lojasiewicz:1999,
  author    = {Stanis{\l}aw {\L}ojasiewicz},
  title     = {On Semi-analytic and Subanalytic Geometry},
  year      = {1999},
  publisher = {Instytut Matematyki, Uniwersytet Jagielloński},
  address   = {Kraków, Poland},
}

@article{DEQKamilov,
  title={Deep equilibrium learning of explicit regularization functionals for imaging inverse problems},
  author={Zou, Zihao and Liu, Jiaming and Wohlberg, Brendt and Kamilov, Ulugbek S},
  journal={IEEE Open Journal of Signal Processing},
  volume={4},
  pages={390--398},
  year={2023},
  publisher={IEEE}
}

\appendix 
\section{Definition and results in convex/non-convex optimization} \label{appendix: definitions and results}

In this section, we present the main definitions and results of convex and non-convex optimization, necessary for the understanding of this work. We refer to \cite{BauschkeCombettes2011,Villa,Beck,RelativeSmoothness,Teboulle2018} as standard references on the following concepts.

\begin{notation}[Interior and Closure]
Let $S\subset \mathbb{R}^n$. We denote, respectively, $\operatorname{int}(S)$ the interior of $S$ and $\overline{S}$ the closure of S (w.r.t Euclidean topology).     
\end{notation}

\begin{definition}[Domain of a function] \label{def: dom}
    Let $f:\mathbb{R}^n \to \mathbb{R} \cup \{+\infty\}$.\\ Then 
    $\operatorname{dom}(f)=\{x \in \mathbb{R}^n : f(x) < +\infty\}.$
\end{definition}

\begin{definition}[$L$-smoothness]\label{L-smothness}
    Let $f:\mathbb{R}^n \to \mathbb{R} \cup \{+\infty\}$. $f$ is said to be $L$-smooth if:
    \begin{itemize}
        \item $f$ is differentiable in $\operatorname{int(dom}(f))$;
        \item $\exists L>0 : \|\nabla f(x)-\nabla f (y)\| \leq L\|x-y\|$ $\forall x,y \in \operatorname{dom}(f)$.
    \end{itemize}
\end{definition}

\begin{definition}[Legendre function]\label{def:Legendre function}
    Let \( h: \mathbb{R}^n \to \mathbb{R} \cup \{+ \infty\} \) be a lower semicontinuous proper convex function. Then $h$ is:
\begin{enumerate}
    \item[(i)] essentially smooth, if \( h \) is differentiable on $ \operatorname{int}(\operatorname{dom}(h))$, with moreover $\| \nabla h(x^k) \| \to + \infty$ \ for every sequence \( \{x^k\}_{k \in \mathbb{N}} \subset \operatorname{dom}(h) \) converging to a boundary point of \( \operatorname{dom}(h)\) as \( k \to +\infty \);
    \item[(ii)] of Legendre type if \( h \) is essentially smooth and strictly convex on \( \operatorname{int}(\operatorname{dom}(h))\).
\end{enumerate}

\end{definition}

\begin{definition}[Bregman divergence] \label{Bregman divergence}
    Given a Legendre function $h: \mathbb{R}^n \to \mathbb{R} \cup \{+ \infty\}$ its associated Bregman distance $D_h:\mathbb{R}^n \times \operatorname{int}(\operatorname{dom}(h)) \longrightarrow [0,+\infty]$ is defined in the following way:
\begin{equation}
    D_h(x,y)=h(x) - h(y) - \langle \nabla h (y), x-y  \rangle  \label{bregman distance}
\end{equation}
with $x \in \mathbb{R}^n$ and $y \in \operatorname{int}(\operatorname{dom}(h))$. If $x \notin \operatorname{dom}(h)$, then $D_h(x,y)=+\infty \ \forall y \in \operatorname{int}(\operatorname{dom}(h))$. 
\end{definition}
\begin{definition}[Bregman Proximal Operator] \label{BProx}
Given a Legendre function $h:\mathbb{R}^n \to \mathbb{R}$ and $\mathcal{R}:\mathbb{R}^n \to \mathbb{R} $ lower semicontinuous and proper (possibly non-convex),  the Bregman Proximal operator of $\mathcal{R}$ (w.r.t $h$) is defined as follows: 
\begin{equation}
    \operatorname{Prox}_{\tau \mathcal{R}}^h(x)=\operatorname*{argmin}_{u\in \mathbb{R}^n}~  \mathcal{R}(u) + \frac{1}{\tau}D_h(u,x). \label{eq of Bprox}
\end{equation}
\end{definition}
Notice that the operator in \eqref{eq of Bprox} may be set-valued. A sufficient condition on $\mathcal{R}$ to have \eqref{eq of Bprox} single-valued is $\mathcal{R}$ being convex and coercive.

\begin{example}[Bregman Proximal operator of $\iota_{[0,a]^n}$ w.r.t $h(x)=-\sum_{i=1}^n \log(x_i)$]\label{Example}
We show here an example of computation of a Bregman Proximal gradient, useful for the aim of this work.
    By definition we need to solve the following problem:

\[\operatorname{Prox}_{\tau [0,a]^n}^h(x)=\operatorname*{argmin}_{u\in \mathbb{R}^n}  \tau \iota_{[0,a]^n}(u) + D_h(u,x),\]
for some $a>0$.
Now using the definition of $D_h$ and observing that \(\iota_{[0,a]^n}(u)=\sum_{j=1}^{n}\iota_{[0,a]}(u_j)\)
we obtain the explicit problem: 
\[\operatorname*{argmin}_{u\in \mathbb{R}^n} \sum_{j=1}^{n}\tau \iota_{[0,a]}(u_j) -\log(u_j)+\log(x_j)+\frac{1}{x_j}(u_j-x_j)  \]

Since the function we are minimizing is separable we can solve the problem reasoning on each component. Thus it suffices to solve:
\begin{equation}
    \operatorname*{argmin}_{u \in \mathbb{R}} \tau \iota_{[0,a]}(u) -\log(u)+\log(x)+ \frac{1}{x}(u-x) \label{Breg prox of projection via burg's}
\end{equation}
We write the problem in a constrained form as: 
\begin{equation}
    u^*_x\in\operatorname*{argmin}_{u \in [0,a]} \left\{ D_h(u,x):= -\log(u)+\log(x)+ \frac{1}{x}(u-x) \right\}
\end{equation}
We have to analyze two cases:
\begin{itemize}
    \item[$\bullet$]
    \item [$\bullet$] if $x>a$, compute $D'_{h}(\cdot,x)$ the derivative of $D_h(u,x)$ w.r.t the first argument. We get 
    \[D'_{h}(u,x)=-\frac{1}{u}+\frac{1}{x}=\frac{u-x}{ux}\]
    This derivative is negative in $(0,x)$ and positive in $(x,+\infty)$ and $u=x$ is the global minimum. Thus for $x>a$ the minimum in $[0,a]$ is attained in $u^*_x=a$.
\end{itemize}
Summarizing
\[\operatorname{Prox}_{\tau [0,a]}^h(x)=\begin{cases}
x & \text{if}\hspace{2mm} x \in [0,a] \\
a & \text{if} \hspace{2mm} x > a
\end{cases}\] Finally, we have that:
\begin{equation}
    \operatorname{Prox}_{\tau [0,a]^n}^h(x)=\Pi_{[0,a]^n}(x)
\end{equation}
where $\Pi$ is the standard Euclidean projection on $[0,a]^n$.
\end{example}
We now recall some definitions and results which will be used to prove convergence of the Mirror Descent algorithm in non-convex setting.
\begin{definition}[Kurdyka-Łojasiewicz property \cite{Attouch2013}]\label{def:KŁ}
 A function $f : \mathbb{R}^n \to \mathbb{R} \cup \{ +\infty \}$ is said to have the Kurdyka-Łojasiewicz property at $x^* \in \operatorname{dom}(f)$ if there exists $\eta \in (0, +\infty)$, a neighborhood $U$ of $x^*$ and a continuous concave function $\psi : [0, \eta) \to \mathbb{R}_+$ such that $\psi(0) = 0$, $\psi$ is $C^1$ on $(0, \eta)$, $\psi' > 0$ on $(0, \eta)$ and $\forall x \in U \cap \{ x \mid f(x^*) < f(x) < f(x^*) + \eta \}$, the Kurdyka-Łojasiewicz inequality holds:

\begin{equation}
  \psi'(f(x) - f(x^*)) \operatorname{dist}(0, \partial f(x)) \geq 1. \label{KurdykaL ineq}  
\end{equation}

Proper lower semicontinuous functions which satisfy the Kurdyka-Łojasiewicz inequality \eqref{KurdykaL ineq} at each point of $\operatorname{dom}(\partial f)$ are called KŁ functions.
\label{Kurdyka-Lojasiewicz}

\end{definition}
\begin{theorem} [\cite{Bolte2018}]\label{KL then convergence of iterates}
    Let $f : \mathbb{R}^n \to \mathbb{R} \cup \{ +\infty \}$ be a proper lower semicontinuous, lower bounded function and $a,b,>0$. Consider a bounded sequence $\{x^n\}_{n \in \mathbb{N}}$ satisfying the following conditions:

\begin{itemize}
    \item [\textbf{H1}] Sufficient decrease condition: $\forall n \in \mathbb{N},$
    \begin{equation}
        f(x^{n+1}) + a \| x^{n+1} - x^n \|^2 \leq f(x^n). 
    \end{equation}
    
    \item [\textbf{H2}] Relative error condition: $\forall n \in \mathbb{N}$, there exists $\omega_{n+1} \in \partial f(x_{n+1})$ such that
    \begin{equation}
        \| \omega^{n+1} \| \leq b \| x^{n+1} - x^n \|. 
    \end{equation}
    
    \item [\textbf{H3}] Continuity condition: Any subsequence $(x^{k_i})_{i \in \mathbb{N}}$ converging towards $\bar x$ verifies
    \begin{equation}
        f(x^{k_i}) \to f(\bar x) \quad \text{as} \quad i \to +\infty. 
    \end{equation}
\end{itemize}

If $f$ verifies the Kurdyka–Łojasiewicz property \eqref{Kurdyka-Lojasiewicz}, then the sequence $\{x^k\}_{k \in \mathbb{N}}$ has finite length, i.e., 
\[
\sum_{k=1}^\infty \|x^{k+1} - x^k\| < \infty,
\]
and it converges to $x^*$, with $x^*$ a critical point of f.
\end{theorem}

The preceding theorem does not require the objective function to be convex. For this reason, the Kurdyka-Łojasiewicz property \eqref{KurdykaL ineq} is a crucial tool for proving the convergence of schemes within non-convex frameworks. Although the class of functions satisfying this property is broad, direct verification can be challenging. We therefore introduce a class of functions known to satisfy the property of interest, the subanalytic functions.
\begin{definition}[Subanalytic sets and functions] \label{def: subanal}

\begin{itemize}
  \item A subset \( S \) of \( \mathbb{R}^n \) is a semianalytic set if each point of \( \mathbb{R}^n \) admits a neighborhood \( V \) for which there exists a finite number of real analytic functions \( f_{i,j}, g_{i,j}: \mathbb{R}^n \to \mathbb{R} \) such that
  \[
  S \cap V = \bigcup_{j=1}^p \bigcap_{i=1}^{q_j} \left\{ x \in \mathbb{R}^n \mid f_{i,j}(x) = 0, g_{i,j}(x) < 0 \right\}.
  \]
  \item A subset \( S \) of \( \mathbb{R}^n \) is a subanalytic set if each point of \( \mathbb{R}^n \) admits a neighborhood \( V \) for which
  \[
  S \cap V = \{ x \in \mathbb{R}^n \mid (x, y) \in U \},
  \]
  where \( U \) is a bounded semianalytic subset of \( \mathbb{R}^n \times \mathbb{R}^p \) with \( p \geq 1 \).
  \item A function \( f : \mathbb{R}^n \to \mathbb{R} \cup \{ +\infty \} \) is called subanalytic if its graph
  \[
  \{(x, y) \in \mathbb{R}^n \times \mathbb{R} \mid y = f(x) \}
  \]
  is a subanalytic subset of \( \mathbb{R}^n \times \mathbb{R} \).
  \end{itemize}
\end{definition}

\begin{lemma}[Stability of subanalytic functions \cite{Shiota2012}] \label{Lemma:stab of subanal}
The sum of two subanalytic functions is subanalytic if
at least one of the two functions maps a bounded set onto a bounded set or if both
functions are lower bounded.
\end{lemma}

\begin{theorem}[KŁ property and subanalytic functions \cite{Bolte2007}] \label{thm:KL and subanal}
Let \( f : \mathbb{R}^n \to \mathbb{R} \cup \{ +\infty \} \) be a proper, subanalytic function with closed domain and assume that \( f \) is continuous on its domain, then \( f \) is a KŁ function.
\end{theorem}
We conclude this section presenting two classes of functions which are subanalytic: the analytic and the semialgebraic functions.

\begin{definition}[Real analytic function]
    A function \( f : A\to \mathbb{R} \ \) with $A \subset  \mathbb{R}^n $ an open set is said to be (real) analytic at \( x_0 \in A \) if there exists a neighborhood \( U \) of \( x_0 \) such that \( f \) can be represented by a convergent power series for all $x\in U$. The function $f$ is said to be analytic if it is analytic at each point of $A$.
\end{definition}

 \begin{definition}[Semialgebraic sets and functions \cite{Attouch2010}]\label{def:Semialgebraic}
\begin{itemize}

  \item A subset \( S \) of \( \mathbb{R}^n \) is a real semialgebraic set if there exists a finite number of real polynomial functions \( P_{i,j}, Q_{i,j} : \mathbb{R}^n \to \mathbb{R} \) such that
  \[
  S = \bigcup_{j=1}^p \bigcap_{i=1}^{q_j} \left\{ x \in \mathbb{R}^n \mid P_{i,j}(x) = 0, Q_{i,j}(x) < 0 \right\}
  \]
  \item A function \( f : \mathbb{R}^n \to \mathbb{R} \cup \{ +\infty \} \) is called semialgebraic if its graph 
  \[
  \{(x, y) \in \mathbb{R}^n \times \mathbb{R} \mid y = f(x) \}
  \]
  is a semialgebraic subset of \( \mathbb{R}^n \times \mathbb{R} \).
\end{itemize}
\begin{lemma}[Elementary properties of semi-algebraic sets and functions\cite{Loja}] \label{Elementary properties}
The closure, the interior, and the boundary of a semi-algebraic set are semi-algebraic.
Polynomials and indicator function of a semialgebraic set are semialgebraic.
\end{lemma}
\end{definition}
\begin{lemma}[Stability of semialgebraic class \cite{Perrin2020}]
    The sum, product, and composition of semialgebraic
functions are semialgebraic.
\end{lemma}

\begin{lemma} [\cite{Shiota2012}] \label{Lemma:Semialg e anal are subanal}
    Semialgebraic functions and analytic functions with full domain are subanalytic.
\end{lemma}

Note, however, that the trivial  extension 
$$
\tilde f(x) =
\begin{cases}
f(x) & \text{if } x \in A,\\[6pt]
+\infty & \text{if } x \notin A.
\end{cases}
$$
to $\mathbb{R}^n$  of an analytic function $f: A \subset \mathbb{R}^n \to \mathbb{R}$, defined on $A \neq \mathbb{R}^n$, may fail to be subanalytic.
As an example, consider $f:(0,+\infty) \to \mathbb{R}$, defined by $f(x)=\sin\left(\frac{1}{x}\right)$.
% }
% \textcolor{red}{
% \[
% \tilde f(x) =
% \begin{cases}
% \sin\left(\frac{1}{x}\right) & \text{if } x > 0,\\[6pt]
% +\infty & \text{if } x \le 0.
% \end{cases}
% \]
% }

An example of function which is analytic on a proper open subset of $\mathbb{R}^n$ with a  subanalytic extension is given in Section SM2.2. This shows that the assumption of full domain is sufficient, but not necessary for subanalyticity.
% \begin{lemma}
%     \textcolor{red}{Let $f:\mathbb{R}^n \to \mathbb{R} \cup  \{+ \infty\}$ be an analytic function. If the graph of $f$ is closed, then $f$ is subanalytic.} 
% \end{lemma}
% \begin{proof}
% We call $U = \operatorname{dom}(f)$, which is open by definition. 
% We now consider the graph of $f$:

% \begin{equation*}
% \Gamma_f = \{(x, y) \in \mathbb{R}^n \times \mathbb{R} \mid y = f(x)\}
% = \{(x, y) \in U \times \mathbb{R} \mid y = f(x)\}.
% \end{equation*}

% We want to show that $\Gamma_f$ is a subanalytic subset of $\mathbb{R}^n \times \mathbb{R}$. In particular we show that $\Gamma_f$ is semianalytc. We consider the following two cases, fixing $x \in \mathbb{R}^n \times R$:
% \begin{itemize}
%     \item $x \notin \Gamma_f:$ in this case, since $\Gamma_f$ is closed, then there exists a neighborhood $V$ of $x$, s.t. $V \cap \Gamma_f = \emptyset$. In this case since the intersection is empty the condition is trivial.
%     \item $x \in \Gamma_f:$ In this case it is possible to choose a neighborhood $V$ of $x$ s.t. \\  
%     $\Gamma_f \cap V = \{x \in \mathbb{R}^n \mid g(x)=0\}$. In particular we consider $g: \mathbb{R}^n \times \mathbb{R}\to \mathbb{R}$, defined as $g(x,y)=y- f(x)$.
% \end{itemize}

% \end{proof}

\subsubsection{Indicator function of $[0,a]^n$ is K\L} \label{indicator fun is KL}
Consider, for any $a>0$, the following function $\iota_{[0,a]^n}$. 
We start proving that $\iota_{[0,a]^n}(x) \hspace{1mm} \text{for} \hspace{1mm}a\hspace{1mm}>0 \hspace{1mm}\text{is semialgebraic}$. From Lemma \ref{Elementary properties} it suffices to show that the set $[0,a]^n$ is semialgebraic. We will show that $(0,a)^n$ is semialgebraic and since $[0,a]^n=\overline{(0,a)^n}$ by Lemma \ref{Elementary properties} we conclude. For definition \ref{def:Semialgebraic} to hold we need to exhibit the polynomials $P_{i,j}$ and $Q_{i,j}$ satisfying the definition. With the following choices:
\[P_{i,j}(x)=0 \hspace{3mm}\forall i  \hspace{2mm}\forall j,\]
\[Q_{ij}(x)=x_j(x_j-a)\hspace{5mm}j=1,...,n,\hspace{1mm}i=1,\]
we can write:
\[(0,a)^n=\cap_{j=1}^n \{x \in \mathbb{R}^n \mid P_{i,j}(x)=0,\hspace{1mm}Q_{i,j}(x)<0\},\]
To conclude that $\iota_{[0,a]^n}$ is a KŁ function, we use Lemma \ref{Lemma:Semialg e anal are subanal}. This lemma implies that the function is subanalytic. We then apply Theorem \ref{thm:KL and subanal}, as $\iota_{[0,a]^n}$ is continuous on its domain, as per Definition \ref{def: dom}. Since $\iota_{[0,a]^n} = 0 \hspace{1mm}\forall x \in [0,a]^n$, the conclusion follows.

\section{Proofs of Section 4}
In this Section, we report the proofs of Propositions 4.3 and 4.5, and Corollary 4.4. Many of the results presented here are technical extensions of well-known results in non-convex optimization, addressing the validity of the Kurdyka–Łojasiewicz property for subanalytic functions with non-closed domains, as  for the KL-dependent functional under consideration. Surprisingly, despite the extensive literature on tailored algorithms for KL minimization in non-convex settings, we could not find any satisfactory result addressing this subtle but crucial technical detail. Clarifying this point not only fills a gap in the existing theory of convergence of non-convex optimization algorithms for such functions, but also strengthens the mathematical foundation of the results presented in the main text.

\subsection{Proof of Proposition
\ref{KL with not-closed domain}} \label{proof of prop 4.4}

Before proving our main result, we illustrate the need for it with two examples.  
Consider the functions $\Psi_1, \Psi_2: \mathbb{R}^2 \to \mathbb{R} \cup \{+\infty\}$, defined by
\[
\Psi_i(x) := \operatorname{KL}(y, A_i x) + \iota_{[0,+\infty)^2}(x), \quad i=1,2,
\]
where $\iota_{[0,+\infty)^2}$ denotes the indicator function of the nonnegative orthant and the matrices $A_i$ are defined by
\[
A_1 = \frac{1}{4} I, \quad 
A_2 = 
\begin{bmatrix}
\frac{1}{2} & 0 \\
0 & \frac{1}{2}
\end{bmatrix}.
\]
In the first case, it is easy to verify that
\[
\operatorname{dom}(\Psi_1) = [0,+\infty)^2 \setminus \{(0,0)\},
\]
which is neither closed nor open. In contrast, for the second case,
\[
\operatorname{dom}(\Psi_2) = (0,+\infty)^2,
\]
which is open. The domain of functions defined in terms of the KL divergence may thus fail to be closed, highlighting the need for Proposition \ref{KL with not-closed domain}.

We recall the following  result.
\begin{theorem}[Theorem 3.1 \cite{Bolte2007}] \label{strong KL}
Let $\Psi : \mathbb{R}^n \to \mathbb{R} \cup \{+\infty\}$ be a subanalytic function with closed
domain, and assume that $\Psi \!\mid_{\operatorname{dom} \Psi}$ is continuous. Let $x^* \in \mathbb{R}^n$ be a
critical point of $\Psi$. Then, there exists an exponent $\theta \in [0,1)$ such that the function
\begin{equation} \label{eq:8}
   x\mapsto \frac{|\Psi(x) - \Psi(x^*)|^{\theta}}{m_\Psi(x)}
\end{equation}
is bounded around $x^*$, where $m_{\Psi}$ is defined as
$$
m_\Psi(x) := \inf \{ \|x^*\| : x^* \in \partial \Psi(x) \}=\text{dist}(0,\partial\Psi(x)),
$$
under the convention $m_\Psi(x)=+\infty$ whenever $\partial \Psi(x) = \emptyset.$
\end{theorem}
In addition we state the following lemma, useful for the proof.
\begin{lemma} \label{Lemma: preimage of a compact set}
    Let $\Psi:\mathbb{R}^n \to \mathbb{R} \cup \{+ \infty\}$ be a subanalytic and coercive function and $K \subset \mathbb{R}$ a compact and subanalytic set. Then $\Psi^{-1}(K)$ is a subanalytic set of $\mathbb{R}^n$.
\end{lemma}
\begin{proof}
Before starting we recall standard results on subanalytic sets \cite{propertysubanalsets,Lojasiewicz:1999}. In particular given $E,F$ two subanalytic sets of $\mathbb{R}^n$, then $E \cap F, E \times F$ are subanalytic.\\
    Now we write:
    \[\Psi^{-1}(K)=\Pi_{\mathbb{R}^n}(\Gamma_\Psi \cap (\mathbb{R}^n \times K)),\]
    where 
$\Pi_{\mathbb{R}^n}: \mathbb{R}^n \times \mathbb{R} \to \mathbb{R}^n,
    $    with $\Pi_{\mathbb{R}^n}(x_1,...,x_n,x_{n+1})=(x_1,...,x_n)$.
    Since $\Psi$ is subanalytic, then $\Gamma_\Psi$ is a subanalytic set of $\mathbb{R}^n \times \mathbb{R}$.
     Similarly, $\mathbb{R}^n \times K$ is subanalytic, as the Cartesian product of two subanalytic sets. Moreover, $\Gamma_\Psi \cap (\mathbb{R}^n \times K)=\{(x, \Psi(x))\mid x \in \mathbb{R}^n, \Psi(x) \in K)\}$ is bounded since $\Psi$ is coercive, thus it is relatively compact. Finally, note that the function $\Pi_{\mathbb{R}^n}$ is subanalytic, as  it is analytic with full domain. To conclude, we use that the image of a relatively compact set by a subanalytic mapping is subanalytic \cite{propertysubanalsets,Lojasiewicz:1999}, which concludes the proof.
\end{proof}

We can now prove Proposition \ref{KL with not-closed domain}.\\
\begin{proof}
We fix $x^*\in\text{crit}(\Psi)$. By definition %\cite{rockafellar1998variational}% 
    $\Psi(x^*)=M < +\infty$. 
    We now divide the proof in two parts:
    \begin{enumerate}
        \item \textbf{Part I}: We suppose there exists a compact set $D'$ such that:
        \begin{enumerate}[label=\roman*.]
            \item $x^* \in D' \subset \operatorname{dom}(\Psi)$;
            \item there exists a neighborhood $U$ of $x^*$, s.t.~$U \cap \operatorname{dom}( \Psi )\subset D'$; 
            \item $D'$ is a subanalytic set.
        \end{enumerate}
        We can thus define $\tilde \Psi:= \Psi+ \iota_{D'}$ such that $\text{dom}(\tilde{\Psi}) = D'$. By construction, we show that $\tilde \Psi$ satisfies Theorem \ref{strong KL} at $x^*$ and that there exists a neighborhood of $x^*$ where $\tilde \Psi\equiv \Psi$, by which the thesis follows.
        \item \textbf{Part II}: We show that such $D'$ always exists, providing a constructive way to define it by considering several exhaustive cases.
    \end{enumerate}
    \paragraph{Part I} In order to apply Theorem \ref{strong KL} we need to show that $\tilde \Psi$ is subanalytic and continuous when restricted to its domain. 
    %, in fact by construction $\Psi=\tilde \Psi$ on $D'$, which is contained in $\operatorname{dom}(\Psi)$. $D'$ being compact is thus closed. Similarly, it holds that 
    $\tilde \Psi \!\mid_{D'}$ is continuous, since it holds that $\tilde \Psi \!\mid_{D'}=\Psi \!\mid_{D'}$, where the second function is continuous by hypothesis. In addition $D'$ is closed, since compact. 
    We now show that $\tilde{\Psi}$ is subanalytic. For that, we first show that $\iota_{D'}$ is subanalytic.  Since the graph of $\iota_{D'}$  is $\Gamma_{\iota_{D'} } = D' \times \{0\}$ and since cartesian product of subanalytic sets is subanalytic \cite{propertysubanalsets,Lojasiewicz:1999}, the claim follows. Since $\Psi$ is subanalyic, it  follows from Lemma \ref{Lemma:stab of subanal} that the sum $\Psi + \iota_{D'}$ remains subanalytic. We show now that: 
    \begin{equation}
         x\mapsto \frac{|\Psi(x) - \Psi(x^*)|^{\theta}}{m_\Psi(x)} \label{eq: KL ratio for Psi}
    \end{equation}
    is bounded around $x^*$.
    By ii., there exists $U$ neighborhood of $x^*$, s.t. $U \cap \operatorname{dom} \Psi\subset D'$. Notice that $\forall y \in U \setminus  D'$, there holds $\Psi(y)=\tilde{\Psi}(y)=+\infty$. Furthermore, since $U=  (U \setminus D')  \cup (U \cap D') $ we conclude that $\Psi$ and $\tilde \Psi$ coincide on $U$. By applying Theorem \eqref{strong KL} to $\tilde \Psi$, we conclude that there exists $U'$ neighborhood of $x^*, \theta \in [0,1)$ and $c>0$ s.t. $$
    \frac{|\tilde \Psi(x) - \tilde \Psi(x^*)|^{\theta}}{m_{\tilde \Psi}(x)} \leq c
    $$ for all $x \in U'$. Since $\Psi\equiv \tilde \Psi$  on $U$, then the function defined by \eqref{eq: KL ratio for Psi} is bounded in $U \cap U'$. Thus, choosing $\psi:=\frac{\cdot ^ {1-\theta}}{c(1-\theta)}$ and choosing $U \cap U'$ as a neighborhood for $x^*$, we conclude that $\Psi $ satisfies \eqref{KurdykaL ineq} at $x^*$.

    \paragraph{Part 2}
    We start defining $D:= \Psi^{-1}((-\infty,M])$, which satisfies i. by definition. Notice that since $\Psi$ is l.s.c., then $D$ is closed and since $\Psi$ is coercive, then $D$ is bounded. Thus $D$ is compact and $x^* \in D$. By definition, there holds $D \cap S = \emptyset$, from which we conclude that $\bar{d}=\text{dist}(D,S)>0$, where $\text{dist}(D,S)$ denotes the set distance between $D$ and $S$. \\
    We now aim at finding the required set $D'$. We observe that there are two possible scenarios:
    \begin{itemize}
        \item $x^* \in \operatorname{int}(D)$: in this case it is possible to find a neighborhood $U$ of $x^*$ s.t. $U \subset D$ and thus trivially $U \cap \operatorname{dom}(\Psi) \subset D.$ It thus suffices to pick $D'=D$.
        \item $x^* \in \partial D$: in this case the existence of such neighborhood is in general false. For instance, if $x^{*}$ is a global minimizer, then  $D=\Psi^{-1}((-\infty,x^*])=\{x^*\}$. The idea is to define $D'$ by adding points which will allow for finding the desired neighborhood $U$. To this end, we firstly define the following set:
        $$
         D^{\epsilon}:=\{x \in \mathbb{R}^n \mid \exists \ y \in D \text{ s.t. } \|x-y\|_2 \leq \epsilon\}.
        $$
        Let us consider $B^{\epsilon}: = D^{\epsilon} \cap \operatorname{dom}\Psi$. The set $B^{\epsilon}$  is bounded and contained in $\operatorname{dom} (\Psi)$, although, in general, not closed. To bypass the problem, we proceed as follows: we observe that
        $\text{dist}(B^{\epsilon},S) \geq \bar{d} - \epsilon$, thus, taking $\epsilon<\bar d$ it holds that $\text{dist}(B^\epsilon,S)>0.$ Thus it follows that $\underset{x \in B^{\epsilon}}{\operatorname{sup}} \Psi(x) = \tilde{M}_\epsilon < +\infty$.
    %We can always find a strictly positive $\epsilon$, s.t. $\underset{x \in B^{\epsilon}}{\operatorname{sup}} \Psi(x) < +\infty.$ 
This leads to $D'=\Psi^{-1}(-\infty,\tilde M_\epsilon]$, which is a compact contained in $\operatorname{dom} (\Psi)$, hence $B^{\epsilon} \subset D'$. It thus suffices to consider $U=B_{\frac{\epsilon}{2}}(x^*)$, which by construction is such that $U \cap \operatorname{dom}(\Psi) \subset D'$.
    \end{itemize}
    To conclude, it remains to show that in both cases $D'$ is a subanalytic set. Let us observe that $D'$ is given by the following rule $D' = \Psi^{-1}((-\infty,M])$ for some $M \in \mathbb{R}$. Using the fact that $\Psi$ is lower bounded (since it is l.s.c. and coercive), we actually get more precisely that $D'=\Psi^{-1}([-m,M])$, where $m = \inf \Psi$. 
    To conclude, we thus simply apply Lemma \ref{Lemma: preimage of a compact set} to $\Psi$ and $K=[-m,M]$.
    \end{proof}

\subsection{Proof of Corollary \ref{Kl corollary}} \label{proof of corollary 4.5}
\begin{proof}
The proof consists of checking the validity of the assumptions of Proposition \ref{KL with not-closed domain}.
    \paragraph{1. $\Psi$ is subanalytic}  
    The expression of $\text{KL}(y,Ax)$ is:
  \begin{equation}
    \text{KL}(y,Ax) = \sum_{i=1}^{m} y_i \log\left( \frac{y_i}{(Ax)_i} \right) + (Ax)_i - y_i= \sum_{i=1}^{m} {\text{kl}_i}(x), \label{KL equation appendix}
\end{equation}
where $\text{kl}_i$ is the $i$-th term in the sum defining $\text{KL}(y,A\cdot)$. The idea is to first show that $\widetilde{\text{KL}}_i$ is subanalytic $\forall \ i=1,...,m$ and then using the stability of the sum of subanalytic functions. 
% The considered functions are analytic, but there holds: $$\operatorname{dom}(\text{kl}_i)=\{x \in \mathbb{R}^n \mid (Ax)_i >0\},$$ which is not the full space $\mathbb{R}^n$. To bypass this issue, we show the desired subanalyticity property using Definition \ref{def: subanal}. 
Our goal is showing that:
$$
\Gamma_{\text{kl}_i}= \{(x,\text{kl}_i(x)) \mid (Ax)_i > 0\} \subset \mathbb{R}^n \times \mathbb{R} 
$$
is a subanalytic set. In particular, we will prove that it is semianalytic, which implies that it is subanalytic. We firstly show that $\Gamma_{\text{kl}_i}$ is closed. This is true since:
$$\lim_{x \to \ z \in \partial \operatorname{dom}({\text{kl}_i})}
{\text{kl}_i}(x)= +\infty.$$ 
Such property implies that all accumulations points of $\Gamma_{\text{kl}_i}$ belong to $\Gamma_{\text{kl}_i}$. Once we have this we consider two cases:
\begin{itemize}
    \item $x \notin \Gamma_{\text{kl}_i}$: in this case it possible to find a neighborhood $V$ of $x$ such that $V \cap \Gamma_{\text{kl}_i}= \emptyset$. Thus, the condition is trivially satisfied.
    \item $x\in \Gamma_{\text{kl}_i}$:
    in this case we consider $V= \mathbb{R}^n \times \mathbb{R}$ and we consider $f,g: \mathbb{R}^n \times \mathbb{R} \to \mathbb{R}$ defined by: 
    \begin{align}
f(x,t) &:= \left(\frac{(Ax)_i}{y_i} \right)^{y_i} \exp(y_i - (Ax)_i) - \exp(t), \\
g(x,t) &:= -(Ax)_i
\end{align}
under the convention $\left(\frac{(Ax)_i}{0} \right)^0=1$. The functions $f,g$ are real-analytic as sum, composition and product of analytic functions. In addition, it holds that:
$$
\Gamma_{\text{kl}_i} \cap V= \Gamma_{\text{kl}_i} = \{x \in \mathbb{R}^n \times \mathbb{R} \mid f(x,t)=0, g(x,t)<0\}.
$$
\end{itemize}
This shows the semianaliticity of $\Gamma_{\text{kl}_i}$ and thus the claim. To conclude we observe that ${\text{kl}_i}$ is non-negative for all $i=1,...,m$ and thus the sum is still subanalytic (Lemma \ref{Lemma:stab of subanal}).\\
% We observe that \eqref{KL equation} is in fact the composition of $\text{KL}(y,\cdot)$ and the linear operator $A$. 
% $\text{KL}(y,\cdot)$ is a separable function that is the sum of the following component-wise functions:
% \begin{equation}
%     t(x)=y \log \left(\frac{y}{x}\right) + x - y
%     \label{KL without T}
% \end{equation}
% Since \(x \mapsto  \log(x),x \mapsto  \frac{1}{x}\) and $x \mapsto  x$ are analytic and the sum and composition of analytic functions is still analytic, then \eqref{KL without T} is analytic. Since both $\text{KL}(y,\cdot)$ and $A$ are analytic, then their composition, i.e. $\text{KL}(y,A\cdot)$ is analytic, thus, in particular, subanalytic  as stated in Lemma \ref{Lemma:Semialg e anal are subanal}.\\ 
Next, note that given $a>0$, the function $\iota_{[0,a]^n}(x)$ is semialgebraic, as we showed in Section~\ref{indicator fun is KL}.\\
To conclude we draw some considerations on $R_\theta
$. A neural network consists in the composition of linear functions (which are analytic) and non-linear activation functions. We focus on the particular case where the activation function is the $\operatorname{Softplus}_{\beta}(x)=\frac{1}{\beta} \log (1+ \exp (\beta x))$, with $\beta>0$.
This function is analytic for any $\beta>0$, being it the composition and sum of analytic functions. Thus, a network with softplus as activation functions is analytic (in particular subanalytic as stated in Lemma \ref{Lemma:Semialg e anal are subanal}), as composition and sum of analytic functions with full domain. Hence $R_{\theta}$ is subanalytic and all the three functions considered are subanalytic. To conclude, we need to check if the sum of these three functions is subanalytic, which in general is not true. We proceed as follows:
\begin{itemize}
    \item [(1)] $\text{KL}(y,A\cdot)+ \iota_{[0,a]^n}$ is subanalytic by Lemma \ref{Lemma:stab of subanal} since both functions are subanalytic and non-negative.
    \item [(2)] $R_{\theta}(K)$ is compact for every compact $K \subset \mathbb{R}^n$ since $R_{\theta}(x)$ is differentiable on $\mathbb{R}^n$, hence in particular continuous. Thus $R_{\theta}$ maps bounded sets in bounded sets.
    \item[(3)] by Lemma \ref{Lemma:stab of subanal}, $(\text{KL}(y,A\cdot)+\iota_{[0,a]^n})+R_{\theta}$ is subanalytic as sum of two subanalytic with at least one mapping bounded sets to bounded sets.
\end{itemize}
\paragraph{2. $\Psi$ is l.s.c. and continuous on its domain} The function $\Psi$ is the sum of three l.s.c.~functions. Moreover,  $\Psi=\text{KL}(y,A\cdot) + R_{\theta}$ on its domain, which is continuous. 
\paragraph{3. $\Psi$ is coercive} This holds since $\Psi= +\infty$ in $\mathbb{R}^n\setminus [0,a]^n$.
\paragraph{4. S is closed} The set $S$ can be explicitly defined in this case. Let us start by observing that $\overline{\operatorname{dom}(\Psi)}=[0,a]^n$. Given $x \in [0,a]^n$, we trivially observe that $\Psi(x)=+\infty$ if $Ax$ has at least one $0$ component. By now considering the operator $A_i : \mathbb{R}^n \to \mathbb{R}$, defined by $A_i(x):=\langle a_i,x \rangle$, with $a_i$ being the $i$-th row of $A$, we consider $K_i = \text{ker}(A_i)$. We observe that
\[
Ax \text{ has at least a $0$ component }\quad \iff \quad x \in \underset{i=1,..,n}{\cup} K_i.
\] 
Thus, $S$ can be written as $S= \left( \underset{i=1,..,n}{\cup} K_i \right) \cap [0,a]^n$. We have that $S$ is closed as the intersection of a closed set with a finite union of closed sets.
\end{proof}

\begin{remark}
    Although $\text{KL}(y,A\cdot)$ is analytic, there holds $$\operatorname{dom}(\text{kl}_i)=\{x \in \mathbb{R}^n \mid (Ax)_i >0\},$$ which is not the full space $\mathbb{R}^n$, so this property cannot be used to prove subanalyticity as a direct consequence.  This is why in the proof we directly rely on Definition \ref{def: subanal}  
    showing that its graph is semianalytic, which is not straightforward. 
A natural attempt would be to consider an analytic function whose zero level set represents the graph of such function, for instance
\[
f(x,t) := \mathrm{KL}(y, Ax) - t.
\]
However, this function is analytic only on a proper open subset of \(\mathbb{R}^n \times \mathbb{R}\), since the KL divergence is not defined everywhere. 
Therefore, \(f\) does not have full domain, whereas Definition~\ref{def: subanal} requires the defining function to be analytic on all of \(\mathbb{R}^n \times \mathbb{R}\). 
A different argument is thus needed to establish the semianalyticity of the graph.

\end{remark}

\subsection{Proposition \ref{prop: convergence}} \label{proof of prop 4.6}

\begin{proof}

Proposition \ref{prop: convergence} is a corollary of Theorem 1 of \cite{Hurault}.
In particular we need to prove the following facts:
\begin{itemize}
    \item[(i)] $h: \mathbb{R}^n \rightarrow \mathbb{R} \ \cup \{+\infty\}$ defined in \eqref{burg's entropy} is of Legendre-type and strongly convex on any bounded convex subset of its domain.

    \item[(ii)] $F := \text{KL}(y,A\cdot)+R_{\theta}$ is proper, $C^1$ on $\operatorname{int}(\operatorname{dom}(h))$, with $\operatorname{dom}(h) \subset \operatorname{dom}(F)$.
    
    \item[(iii)] $\mathcal{R} := \iota_{[0,a]^n}$ is proper, lower semicontinuous with $\operatorname{dom}(\mathcal{R}) \ \cap \ \operatorname{int}(\operatorname{dom}(h))  \neq \emptyset$.
    
    \item[(iv)] $\Psi:=F+\mathcal{R}$ is lower-bounded, coercive and verifies the Kurdyka-Łojasiewicz property.
    
    \item[(v)] For $x \in \operatorname{int}(\operatorname{dom}(h))$, $$T_\tau(x)=\operatorname*{argmin}_{x\in \mathbb{R}^n} \mathcal{R}(x)+ \langle x - x^k, \nabla F(x^k) \rangle + \frac{1}{\tau} D_{h}(x,x^k)$$ is nonempty and included in $\operatorname{int}(\operatorname{dom}(h))$.
    \item [(vi)] There is $L>0$ such that, $Lh-F$ is convex on $\operatorname{Conv}(\operatorname{dom}(\mathcal{R})) \cap \operatorname{int}(\operatorname{dom}(h))$. \label{(vi)}
    
    \item[(vii)] \textit{For all} $\alpha > 0$, $\nabla h$ and $\nabla F$ are Lipschitz continuous on $\{ \Psi(x) \leq \alpha \}$.
\end{itemize}
We show only the points that are not straightforward. This proof is based on Section E.3 of \cite{Hurault}. However, in Section \ref{proof on nolip}, we establish the assumption using a general argument, rather than one tailored to the specific choice of $h$.
   \subsubsection{Assumption (iv): properties on $\Psi=\text{KL}
   (y,A\cdot)+R_{\theta}+\iota_{[0,a]^n}$}
These properties were already proved in  \ref{proof of corollary 4.5}.
\subsubsection{ Assumption (v): well-posedness}
We recall the definition of $T_{\tau}$, for $k\geq 0$ and for a given $x^k \in \operatorname{int}(\operatorname{dom}(h))$:
\begin{equation}
    T_{\tau}(x^k)\in\operatorname*{argmin}_{x\in \mathbb{R}^n} \mathcal{R}(x)+ \langle x - x^k, \nabla F(x^k) \rangle + \frac{1}{\gamma} D_{h}(x,x^k).
\end{equation}
Using the definition of $D_h(\cdot,\cdot)$ \eqref{Bregman divergence} and omitting all terms that do not depend on $x$, we consider the following function:
\[\beta(x)=\mathcal{R}(x)+\langle x,\nabla F(x^k) \rangle+ \frac{1}{\gamma}(h(x)-\langle \nabla h(x^k),x \rangle)\]
Since $\mathcal{R}=\iota_{[0,a]^n}$ is convex and $h$ is stricly convex then $\beta$ is strictly convex, thus its set of minimizers is empty or single-valued. In addition, $\beta$ is equal to $+\infty$ in $\mathbb{R}^n \setminus[0,a]^n$, hence it is coercive. Thus $\beta$ admits a unique global minimiser that belongs to $[0,a]^n$. To conclude observe that $\beta=+\infty$ on the boundary of $[0,a]^n$, since $h(x)=+\infty$ if $x \notin \mathbb{R}^n_{++}$ and the other terms are bounded in $[0,a]^n$. Thus the minimiser must belong to $(0,a)^n \subset \operatorname{int}(\operatorname{dom}(h))$, as $\beta$ is proper.
\subsubsection{Assumption (vi): NoLip} \label{proof on nolip} $\text{KL}(y,A\cdot)$ satisfies NoLip property w.r.t Burg's entropy \eqref{burg's entropy} with $L\geq\|y\|_1$ \cite{nolip}. Thus it remains to prove that: 
\[\exists L'>0 \hspace{1mm} \text{such that} \hspace{1mm}L'h-R_{\theta} \    
 \text{is convex on} \operatorname{Conv}(\operatorname{dom(\mathcal{R})}) \cap \operatorname{int(\operatorname{dom}( h))}=(0,a]^n.\]
This is indeed a sufficient condition showing the existence of an $ \tilde{L}>0$ satisfying NoLip property (Definition 4.1 of the main text) for the composite functional, simply by choosing $\tilde{L}=L'+L$. In fact, $\tilde{L}h-(\text{KL}(y,A\cdot)+R_{\theta})=(Lh-\text{KL}(y,A\cdot))+ (L'h-R_{\theta})$. Since the sum of convex functions is convex, this proves the desired property.
To show the existence of such a $L'$, given the regularity of $h$ and $R_{\theta}$ it is enough to show that \cite{nolip}: 
\[\hspace{1mm} \exists L'>0 \ \text{s.t.} \ L'\nabla^2 h(x)-\nabla^2 R_{\theta}(x) \succeq 0, \forall x \in \operatorname{Conv}(\operatorname{dom (\mathcal{R})) \cap \operatorname{int(\operatorname{dom(h)})}}.\]
We start observing that both $h$ and $R_{\theta}$ are $C^2$ functions, thus their Hessians are symmetric (on every point) for Schwarz's Theorem. In addition, for the real spectral Theorem $\nabla^2h(x)$ and $\nabla^2 R_{\theta}(x)$ have real eigenvalues $\forall x \in \mathbb{R}^n_{++}$. We want to show then that:
\[\exists L'>0 \hspace{1mm}\text{such that} \hspace{1mm} \inf_{x \in (0,a]^n } \min \lambda (L'\nabla^2 h(x)-\nabla^2R_{\theta}(x)) \geq0\]
where by $\lambda(\cdot)$ we denote the spectrum of a matrix. 
To show such property, we need the following result:
\begin{lemma}
Let $C,D$ be real symmetric matrices. Let $\lambda^{-}$ and $\lambda^+$ the minimum and maximum eigenvalue of $C$, and $\mu^-$ and $\mu^+$ respectively for $D$. Then, calling $\gamma^-$ and $\gamma^+$ the minimum and maximum eigenvalue of $C+D$ there holds:
\[\lambda^-+\mu^- \leq \gamma^-\leq \gamma^+\leq \lambda^++\mu^+.\]
\begin{proof}
We prove only: $\lambda^-+\mu^- \leq \gamma^-$, since the proof of the other inequality is similar. 
We know:
\[\lambda(C) \subset [\lambda^-,\lambda^+], \quad \lambda(D) \subset [\mu^-,\mu^+].\]
Thus
\[\lambda(C-\lambda^-I)\subset[0,\lambda^+ -\lambda^-], \quad \lambda(D-\mu^-I) \subset [0,\mu^+-\mu^-].\]
Therefore $C-\lambda^-I$ and $D-\mu^-I$ are positive semi-definite.
Thus \((C-\lambda^-I)+(D-\mu^-I)=(C+D)-(\lambda^-+\mu^-)I \hspace{3mm}\text{is positive semi-definite}\), since sum of positive semi-definite matrices is positive semi-definite.\\
The latter implies:
\[\lambda((C+D)-(\lambda^-+\mu^-)I)\subset[0,+\infty].\]
We conclude observing that 
\[C+D=(C+D)-(\lambda^-+\mu^-)I+(\lambda^-+\mu^-)I,\]
by which we have that for the spectrum, there holds: 
\[\lambda(C+D)\subset [\lambda^-+\mu^-,+\infty.] \quad \square\]
\end{proof}
\end{lemma}

To conclude, we have to pay attention to the fact that the Hessian matrices of $h$ and $R_{\theta}$ depend on $x \in \operatorname{int}(\operatorname{dom}(h))$, while, ideally, we would like to get a uniform bound.
Since, by assumption, $R_{\theta}$ is analytic, then $R_{\theta}$ has a Lipschitz gradient (over a compact set), which implies that for $x \in (0,a]^n:$ 
\[s(x)=\|\nabla^2 R_{\theta}(x)\|=\max{\{|\min\lambda(\nabla^2 R_{\theta}(x))|,|\max \lambda (\nabla^2 R_{\theta}(x))|}\} \leq M.\]
On the other side, since $h$ is strongly convex on any bounded and convex set, we have:
\[\inf_{(0,a]^n} \min \lambda(\nabla^2h(x))\geq m >0.\]
Hence, it is thus possible to find $L'>0$ s.t
\[ \inf_{x \in (0,a]^n } \min \lambda (L'\nabla^2 h(x)-\nabla^2R_{\theta}(x)) \geq L'm-M\geq0.\]
\end{proof}
Observe that this way a proceeding is still true for every Legendre function $h$ that is strongly convex on bounded and convex set and do not use the analytic expression of Burg's entropy.
Now the thesis is a direct application of Theorem 1 of \cite{Hurault}.\\
For the sake of completeness we show that, $h$ chosen as the Burg's entropy is strongly convex on bounded and convex set of its domain. Let $B \subset \mathbb{R}^n_{++}$ be a bounded and convex set. By \eqref{burg's entropy}, it follows that: \[\nabla^2 h(x)=\operatorname{diag}\left(\frac{1}{x^2}\right), \text{ with }\frac{1}{x^2}=\left(\frac{1}{x_1^2},...,\frac{1}{x^2_n}\right).\]
Since $h \in C^2{(\mathbb{R}^n)}$, it suffices to show that. there exists $m>0$ such that:
\[\inf_{x\in B} \min \lambda(\nabla^2h(x))\geq m >0.\]
Since $B \subset \mathbb{R}^{n}_{++}$ is bounded, then there exists $M>0$ s.t $B \subset (0,M]^n$ and to conclude we observe:
\[\inf_{x\in B} \min \lambda(\nabla^2h(x))\geq \inf_{x\in (0,M]^n} \min \lambda(\nabla^2h(x))= \inf_{u \in (0,M]}\frac{1}{u^2}=\frac{1}{M^2}>0.\]
\iffalse
\section{More numerical tests}
\subsection{Checking the symmetry of $B_{\epsilon}$} \label{checking symmetry}
To check the symmetry of the Jacobian $B_{\epsilon}$, that we remind below:
\begin{equation}
    B_{\epsilon}=\left((1+\epsilon)I-\frac{\partial f_\theta(x^{\infty})}{\partial x}\right),
\end{equation}
we need only to check the symmetry of $J=\frac{\partial f_\theta(x^{\infty})}{\partial x}$.
Observe that $J$ is symmetric if and only if $Jv=J^{\top}v$, for each vector $v$ (when this makes sense, i.e. if $J$ is a squared matrix). Since this is not easy theoretically, we check this property numerically as follows. 
\begin{itemize}
    \item we generate $v$ such that$v_i \sim \mathcal{U}((0,1))$
    \item we compute $\|Jv-J^{\top}v\|_2$
\end{itemize}

Since the Jacobian depends on $x^{\infty}$, i.e. the reconstruction image associated to a measurement $y$, we considered one Jacobian for each image in the test set and we considered two different operator $A$. 
In the hystograms below we considered the relative difference between $Jv$ and $J^{\top}v$ and we took a mean over the test set.

\begin{figure}[h]
    \centering
    \includegraphics[width=0.85\linewidth]{images/hystograms.png}
    \caption{Study on the symmetry of $J$, generating 10000 random vectors $v$.}
    \label{fig:hysto}
\end{figure}

As shown in figure \ref{fig:hysto}, in both cases the maximum of the relative difference is really small, meaning that we can consider $J$ symmetric on the subspace generated by all the vectors $v$.
\fi

\section{Evaluation Metrics} \label{evaluation metrics}

In the main text we consider the three following well-known metrics to evaluate the quality of the reconstructions:
\begin{itemize}
    \item \textbf{PSNR} (Peak Signal-to-Noise Ratio): This metric quantifies the reconstruction quality of an image in terms of pixel-wise similarity. It takes values in decibels (dB), typically ranging from 20 dB to 40 dB for natural images.
    \item \textbf{SSIM}  (Structural Similarity Index Measure) \cite{SSIM}: SSIM assesses the structural similarity between two images, considering luminance, contrast, and structure. Its values range between 0 and 1, where 1 signifies perfect similarity and 0 indicates no similarity.
    \item \textbf{LPIPS} (Learned Perceptual Image Patch Similarity)  \cite{LPIPS}: This metric measures the perceptual distance between two images by utilizing a pre-trained neural network model, which is designed to correlate better with human perception of similarity. It takes positive values, generally close to 0 for very similar images and increasing as dissimilarity grows. 
\end{itemize}

\section{Pre-training}
\label{pre-training}
We consider here DEQ-RED method. We describe in the following a pre-training strategy which could be beneficial for two reasons.  Firstly, given that each evaluation of the Deep Equilibrium models considered necessitates convergence, we empirically observed that randomly initialized weights often required a substantial number of iterations to reach a fixed point, translating to increased GPU computational time during the training process. Secondly, considering the non-convex nature of the training objective, a judicious initialization of the network weights can play a crucial role in guiding the optimization process towards more favorable local minima, potentially leading to improved overall performance and a more robust final model. For the reasons, we employed a pre-training strategy. In particular, given $R_{\theta}(x)=\frac{1}{2}\|x-N_{\theta}(x)\|^2$ we considered its associated Gradient step denoiser \cite{Hurault1} $D_{\theta}(x):=x-\nabla R_{\theta}(x)$ and we trained it on Gaussian denoising task. To see the benefits of this pre-training strategy, we considered a small training problem where we simultaneously trained two models: one pre-trained and the other not. The dataset considered was a subset of the dataset proposed in Section \ref{dataset}, and for the training, we followed the same modalities as in Section \ref{loss function}
\begin{figure}[h]
    \centering
    \includegraphics[width=1\linewidth]{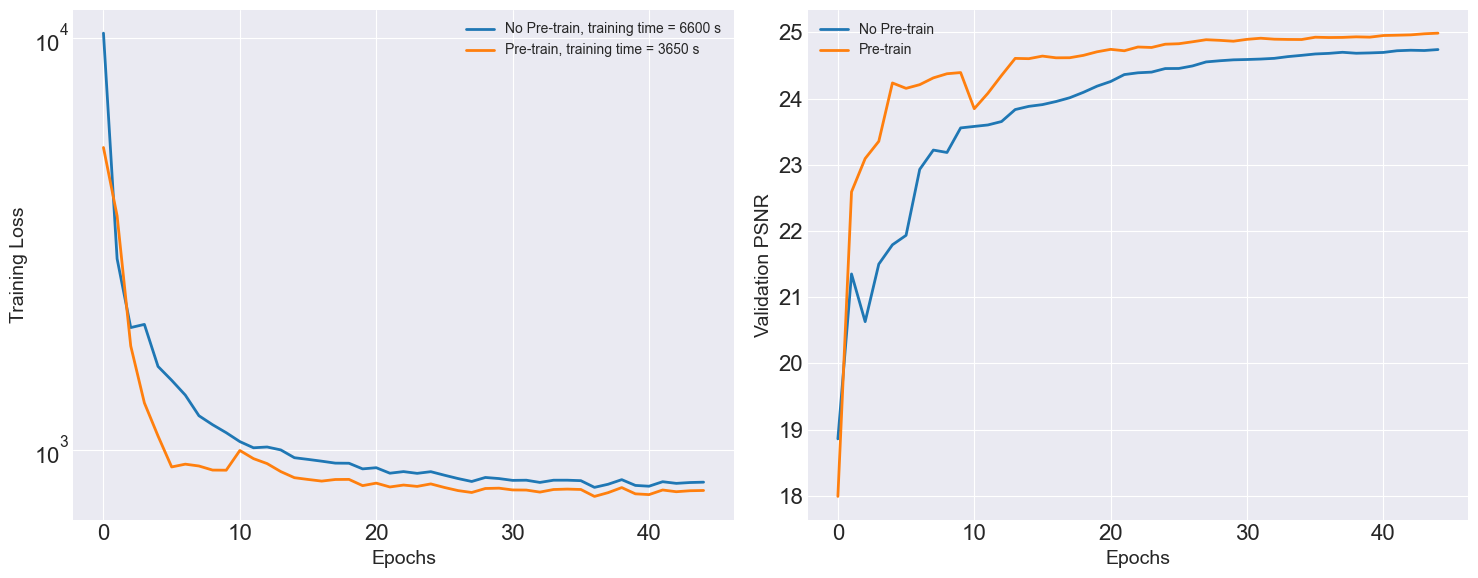}
    \caption{Effect of pre-training for \text{DEQ-RED.}}
    \label{fig:pre-training}
\end{figure}
As illustrated in Figure \ref{fig:pre-training}, the empirical findings unequivocally highlight the advantageous impact of pre-training. This benefit is manifested both in a marked decrease in the computational resources required for training, leading to shorter training durations, and a discernible enhancement in the PSNR metric when evaluated on the validation dataset.

\end{document}